%% file: wsindy_convergence.tex
\documentclass[review,onefignum,onetabnum]{siamart190516}

\input{ex_shared}

\ifpdf
\hypersetup{
  pdftitle={Asymptotic consistency of the WSINDy algorithm in the limit of continuum data},
  pdfauthor={D. A. Messenger D. M. Bortz}
}
\fi

\begin{document}

\maketitle

\begin{abstract}
In this work we study the asymptotic consistency of the weak-form sparse identification of nonlinear dynamics algorithm (WSINDy) in the identification of differential equations from noisy samples of solutions. We prove that the WSINDy estimator is unconditionally asymptotically consistent for a wide class of models which includes the Navier-Stokes equations and the Kuramoto-Sivashinsky equation. We thus provide a mathematically rigorous explanation for the observed robustness to noise of weak-form equation learning. Conversely, we also show that in general the WSINDy estimator is  only conditionally asymptotically consistent, yielding discovery of spurious terms with probability one if the noise level is above some critical threshold and the nonlinearities exhibit sufficiently fast growth. We derive explicit bounds on the critical noise threshold in the case of Gaussian white noise and provide an explicit characterization of these spurious terms in the case of trigonometric and/or polynomial model nonlinearities. However, a silver lining to this negative result is that if the data is suitably denoised (a simple moving average filter is sufficient), then we recover unconditional asymptotic consistency on the class of models with locally-Lipschitz nonlinearities. Altogether, our results reveal several important aspects of weak-form equation learning which may be used to improve future algorithms. We demonstrate our results numerically using the Lorenz system, the cubic oscillator, a viscous Burgers growth model, and a Kuramoto-Sivashinsky-type higher-order PDE.
\end{abstract}


\section{Introduction}\label{sec:intro}

\subsection{Overview}

A widespread challenge in the natural sciences is to create a mathematical model which makes accurate predictions, is mathematically analyzable, and amenable to parameter estimation from data. Typically, parameters exhibit a nonlinear relationship with the observed data, which explains the widespread use of nonlinear least-squares methods to fit parameters by minimizing the difference between the experimental data and numerical simulations. In many cases, parameter uncertainty is layered on top of uncertainty in the model itself, which has led to the field of {\it model selection} \cite{Akaike1974IEEETransAutomControl,Akaike1977Applicationsofstatistics}, whereby criteria are devised to select an appropriate (e.g.\ parsimonious) model from a host of candidate models.

A subject of fervent study in recent years has been the use of data to learn the correct model equations for a given phenomenon, along with the model parameters, without resorting to the laborious nonlinear least-squares approach which often involves many expensive forward solves of candidate models. With equations in hand, one can conduct further experiments computationally which may not be immediately accessible in laboratory conditions. The weak-form sparse identification of nonlinear dynamics algorithm (WSINDy), introduced in \cite{messenger2020weak,messenger2020weakpde}, is one such equation learning algorithm that has been observed to offer several advantages over other equation learning methods, including robustness to noise, high accuracy, low computational cost, and extensibility across many modeling paradigms \cite{messenger2022learning,messenger2022online,messenger2022learninganisotropic}.
WSINDy is based off of the SINDy algorithm \cite{brunton2016discovering,rudy2017data}, which popularized the sparse equation learning technique and demonstrated its viability in discovering ordinary differential equations, partial differential equations, normal forms, and reduced-order models, to name a few. For other important works on sparse equation learning see \cite{xu2008spatiotemporal,schmidt2009distilling,schaeffer2017learning,kaheman2022automatic} and the references therein.

The key difference between SINDy and WSINDy is that the latter discretizes an integral representation of the dynamics, avoiding computation of pointwise derivatives. This leads to implicit noise filtering effects and relaxed constraints on the smoothness of the underlying ground-truth data. In addition to WSINDy as introduced in \cite{messenger2020weak,messenger2020weakpde}, these and other advantages of weak-form equation learning have now been reported in several other studies \cite{schaeffer2017sparse,rosenfeld2019occupation,pantazis2019unified,wang2019variational,wang2021variational,gurevich2019robust,reinbold2020using,chen2020methods,tang2022weakident,alves2022data}, leading to a consensus that exploiting integration can provide significant improvements in the discovery of equations from data. Nevertheless, other promising methods for reducing the deleterious effects of noise do exist, see in particular works that incorporate automatic differentiation \cite{kaheman2022automatic,LagergrenNardiniMichaelLavigneEtAl2020ProcRSocA}, ensembling \cite{fasel2022ensemble}, and improved optimization methods \cite{wentz2022derivative,hokanson2022simultaneous}. However, each of these techniques can easily be combined with a weak-form approach, as demonstrated in \cite{fasel2022ensemble,bertsimas2022learning}.

The goal of this paper is to provide rigorous justification for the observed robustness to noise of WSINDy and other weak form methods. We accomplish this by determining conditions under which a WSINDy model is asymptotically consistent with the true underlying model in the limit of continuum data. We consider the performance of WSINDy on noisy samples $\Ubf^{(m)}$ of some solution $u$ of a differential equation on grids $(\Xbf^{(m)},\tbf^{(m)})$ of increasing resolution within an underlying spatiotemporal domain $\Omega\times(0,T)$. This is cast as a support recovery problem for the support $\supp{\wstar}$ of the true weight vector $\wbf^\star$, which is nonzero only on the subset of terms that exist in the true model. It should be noted that throughout we assume that the true model is contained in the library of models considered, however extension beyond this case is possible (see Section \ref{discussion} discussion point II).

As we will see, for each class of models examined, there exists a critical noise level $\sigma_c>0$ such that for noise levels $\sigma<\sigma_c$, we have $\what^{(m)}$ satisfying $\supp{\what^{(m)}} = \supp{\wstar}$ with probability rapidly approaching one as $m\to \infty$, where $\what^{(m)}$ is the WSINDy model at discretization level $m$, and $m$ is inversely proportional to the volume element of the spatiotemportal grid $(\Xbf^{(m)},\tbf^{(m)})$. Most importantly, we identify a large class of systems for which $\sigma_c=\infty$, leading to unconditional support recovery as $m\to \infty$. We also prove that suitably denoising the data leads to unconditional asymptotic consistency on the class of models with locally-Lipschitz nonlinearities, in particular this is true if a simple moving average filter is used.

\subsection{Related work}

Despite the widespread use and development of sparse data-driven equation learning algorithms applied to dynamical systems since the onset of SINDy, very few studies have examined the performance of these methods from a mathematically rigorous perspective. The authors of \cite{tran2017exact} and \cite{schaeffer2020extracting} report sparse recovery guarantees for chaotic and structured ODE systems (respectively), however each result requires the availability of point-wise derivatives. Other results include convergence of the STLS algorithm used in SINDy \cite{zhang2019convergence} and recovery guarantees for sparse polynomial regression \cite{ho2020recovery}, however the former only asserts convergence to a local minimum of the $\ell_0$-regularized least-squares cost function, and the latter relies on sparsity of the measurement noise. The authors of \cite{russo2022convergence} have recently proved that surrogate models learned using WSINDy converge to the true dynamics as the basis used to represent the dynamics grows in dimension, although this result concerns only continuum noise-free data.

In this work we undertake the problem of analyzing the performance of WSINDy applied to noisy discrete datasets by utilizing a suitable continuum limit, whereby many practical insights can be gleaned. One important related work is \cite{he2021asymptotic}, in which the authors show that with local polynomial differentiation together with filtering kernels that obey certain asymptotics, the LASSO estimator is asymptotically consistent for PDE recovery problems from quadratic libraries of arbitrary spatial differentiation order. This approach appears to be extendable to a wide class of model libraries, but requires the so-called {\it mutual incoherence} condition on the design matrix, which is restrictive in settings with correlated columns (a common paradigm in sparse equation learning from data). In addition, the noise-free data in \cite{he2021asymptotic} is assumed to be classical, which we show in this work is not necessary if one employs the weak form. Nevertheless, we conjecture that results in the present article may be combined with those in \cite{he2021asymptotic} to yield asymptotically consistent algorithms for models with combined weak and classical derivative evaluation. 

\subsection{Summary of results}

As alluded to above, in this work we focus on the performance of WSINDy in the continuum limit as the computational grid $(\Xbf^{(m)},\tbf^{(m)})$ on which the noisy data $\Ubf^{(m)}$ is sampled becomes dense in the domain of definition $\Omega\times (0,T)$ of the dynamical system. Motivated by future algorithmic developments, we aim to provide explicit results whenever possible. These include exponential rates of concentration which we use to quantify the probability of support recovery, and explicit bounds on the critical noise below which we recover the correct model unconditionally in the limit $m\to \infty$ (where $m$ is the number of points that the test function $\psi$ is supported on at the grid resolution $(\Delta x^{(m)},\Delta t^{(m)})$, see \eqref{def_of_m}). The following is a qualitative summary of results with reference to specific theorems and lemmas (see Section \ref{sec:assumptions} for assumptions used in many of these results).

\begin{enumerate}[label=(\Roman*)]
\item We prove asymptotic results that explain the empirically observed robustness of weak-form equation learning methods (Theorem \ref{conv_theorem_nosmooth} and supporting lemmas \ref{lemm:bounds on sigma_c} and \ref{lemm:subset_support_rec}, see also Theorem \ref{full_theorem_nosmooth} and Lemma \ref{lemm:fullsupprec}). 
\item We prove that WSINDy is capable of identifying models from non-smooth data (Theorem \ref{conv_theorem_nosmooth}), which was demonstrated empirically \cite{messenger2020weakpde}, and we quantify the effect of smoothness on convergence (theorems \ref{matrix_concentration} and \ref{smoothconc}, utilizing Lemma \ref{trapconv}).
\item We specify both the class of models and the bounds on the critical noise level below which the weak form is asymptotically consistent (lemmas \ref{lemm:bounds on sigma_c} and \ref{lemm:subset_support_rec}).
\item We prove that suitably denoising the data (e.g.\ with a simple moving average filter) results in unconditional asymptotic consistency of WSINDy over the class of models with locally-Lipschitz nonlinearities (Theorem \ref{smoothconc}). 
\end{enumerate}

\subsection{Outline}
 
 We first cover in Section \ref{sec:prelims} the preliminaries necessary to analyze the limit of large data. These include an overview of the WSINDy algorithm for ordinary and partial differential equations (\ref{sec:wsindyoverview}) which may be skipped for readers familiar with \cite{messenger2020weakpde}, a definition of the continuum problem and related notation (\ref{sec:contprob}), an intuition-building discussion of the bias resulting from taking the continuum limit (\ref{sec:controllablebias}), and lastly the assumptions used throughout and explanations thereof (\ref{sec:assumptions}). In Section \ref{matcon} we prove that the WSINDy discrete linear system concentrates at an exponential rate to the associated continuum linear system under the assumptions in Section \ref{sec:assumptions}. In Section \ref{sec:consistency} we prove that results in Section \ref{matcon} imply {\it conditional} asymptotic consistency for raw data (\ref{nosmooth}) and {\it unconditional} consistency for a wide range of data filtering techniques (\ref{smooth}). Finally, Section \ref{exp} contains numerical examples demonstrating the results from Section \ref{sec:consistency} in practice. In particular, we show that with a simple moving-average filter we achieve stable and accurate recovery of systems for noise levels where WSINDy without filtering fails due to the existence of a critical noise. In the appendix we include a table of notations used (Appendix 
 \ref{app:notation}), supporting lemmas (Appendices \ref{app:trapconv},\ref{app:comblemms},\ref{app:conclemmas}), additional information on numerical examples (Appendices \ref{app:hyperKS},\ref{app:sigma_est}), and extension of several results in Section \ref{sec:consistency} to more general settings (Appendix \ref{app:supprec}).

\section{Preliminaries}\label{sec:prelims}

\subsection{Overview of WSINDy} \label{sec:wsindyoverview}

Let $\Ubf = u(\Xbf, \tbf) + \ep$ be a spatiotemporal dataset defined on the $d$-dimensional spatial grid $\Xbf \subset \Omega\subset \Rbb^d$ over timepoints $\tbf\subset [0, T]$, where $u:\Rbb^d\times\Rbb\to \Rbb^n$ is a weak solution to the PDE
\begin{equation}\label{gen_pde}
\partial^{\alpha^0}u(x,t) = \partial^{\alpha^1}g_1(u(x,t))+\partial^{\alpha^2}g_2(u(x,t))+\dots+\partial^{\alpha^S} g_S(u(x,t)), \qquad x\in \Omega,\ t\in (0,T).
\end{equation}
Here $\ep$ represents i.i.d.\ measurement noise associated with point-wise evaluations of $u$. The WSINDy algorithm uses the weak form of the dynamics \eqref{gen_pde} to identify a PDE model for $\Ubf$ by discovering functional representations of the nonlinear differential operators\footnote{Throughout we use the multi-index notation $\alpha^s =~ (\alpha^s_1,\dots,\alpha^s_d,\alpha^s_{d+1}) \in \Nbb^{d+1}$ to denote partial differentiation with respect to $x = (x_1,\dots, x_d)$ and $t$: 
\[\partial^{\alpha^s}u(x,t) = \frac{\partial^{\alpha^s_1+\cdots+\alpha^s_d+\alpha^s_{d+1}}}{\partial x_1^{\alpha^s_1}\dots \partial x_d^{\alpha^s_d}\partial t^{\alpha^s_{d+1}}}u(x,t).\]
 We will avoid using subscript notation such as $u_x$ to denote partial derivatives, instead using $\partial^\alpha u$ or $\partial_x u$.}%
 $(\partial^{\alpha^s}g_s(\cdot))_{s\in [S]}$ in a computationally efficient sparse regression framework. We adopt a dictionary learning approach and use a basis of $J$ {\it trial functions} $\CalF:=(f_j)_{j\in [J]}$ and {\it differential operators} specified by the set of multi-indices $\pmb{\alpha}:=(\alpha^s)_{s\in [S]}$. We assume that $(g_s)_{s\in [S]} \subset \text{span}(\CalF)$, allowing for the representation of \eqref{gen_pde}:
\begin{equation}\label{diffform}
\partial^{\alpha^0} u  = \sum_{s=1}^S\sum_{j=1}^J \wstar_{(s-1)J+j} \partial^{\alpha^s} f_j(u),
\end{equation}
with the assumption that the true $\wstar\in \Rbb^{SJ}$ is sparse.

To convert \eqref{diffform} to a weak form, we  convolve the equation against a sufficiently smooth \textit{test function} $\psi(x,t)$, compactly-supported in $\Omega\times (0,T)$, arriving at a {\it convolutional weak form} of the equation
\begin{equation}\label{conv_form}
\Big(\partial^{\alpha^0}\psi\Big) * u (\xbf,t) = \sum_{s=1}^S\sum_{j=1}^J \wstar_{(s-1)J+j} \Big(\partial^{\alpha^s}\psi\Big) * f_j(u)(\xbf,t).
\end{equation}
%
Equation \eqref{conv_form} only holds if the support of $\psi$ centered at $(\xbf,t)$ lies inside $\Omega\times(0,T)$, or 
\begin{equation}\label{IBPbcs}
\supp{\psi(\xbf-\cdot,t-\cdot)}\subset \Omega\times (0, T)
\end{equation}
which serves to eliminate any boundary terms that arise during integration by parts. 

To discretize \eqref{conv_form}, we first select a finite set of {\it query points} $\CalQ := \{(\xbf_k,t_k)\}_{k\in[K]}$ satisfying \eqref{IBPbcs}. We then define discrete convolution operators
\[\Psi^s := \partial^{\alpha^s}\psi(\Ybf,\mathfrak{t})(\Delta x)^d\Delta t, \qquad s = 0, \dots, S\]
where $(\Ybf,\mathfrak{t})$ denotes a centered reference grid at the same spatiotemporal resolution $(\Delta x, \Delta t)$ as $(\Xbf,\tbf)$ and the factor $(\Delta x)^d\Delta t$ is the uniform weight of the trapezoidal rule given compact support of $\psi$ which eliminates boundary terms with weight $1/2$. The WSINDy linear system $(\Gbf,\bbf)$ with $\Gbf\in \Rbb^{K\times \mathfrak{J}}$, $\bbf\in \Rbb^{K\times n}$ is then defined by  
\begin{equation}\label{conv_disc}
\begin{dcases} \hspace{1.5cm}\bbf_k := \Psi^0 * \Ubf (\xbf_k,t_k),\\
\Gbf_{k,(s-1)J+j} :=  \Psi^s * f_j(\Ubf) (\xbf_k,t_k),\end{dcases}
\end{equation}
where $K$ is the number of convolution query points, $\mathfrak{J} =SJ$ is the size of the candidate model library, and $n$ is the dimension of the state variable $u$ being observed. The discrete $(d+1)$-dimensional convolution between $\Psi^s$ and $f_j(\Ubf)$ at a point $(\xbf_k,t_k) \in (\Xbf,\tbf)$ is defined by
\[\Psi^s*f_j\left(\Ubf\right)(\xbf_k,t_k) := \sum_{\ell_1=1}^{N_1}\cdots\sum_{\ell_{d+1}=1}^{N_{d+1}} \Psi^s_{k_1-\ell_1,\dots,k_{d+1}-\ell_{d+1}} f_j\left(\Ubf_{\ell_1,\dots,\ell_{d+1}}\right).\]

To arrive at an approximate model we solve for a sparse solution $\what$ such that $\Gbf\what \approx \bbf$ using a {\it modified sequential thresholding least squares} algorithm (proposed in \cite{messenger2020weakpde}) $\what = \text{MSTLS}(\Gbf,\bbf;\, \CalL,\pmb{\lambda})$ that combines the traditional sequential thresholding algorithm 
\begin{equation}\label{STLS}
\text{STLS}(\Gbf,\bbf;\, \lambda)\qquad \begin{dcases} 
\hspace{0.1cm} \wbf^{(\ell+1)} = H_\lambda\left(\argmin_{\supp{\wbf}\subset S^{(\ell)}}\nrm{\Gbf\wbf-\bbf}_2^2\right)\\
S^{(\ell+1)} =\supp{\wbf^{(\ell+1)}}
\end{dcases}
\end{equation}
with a line search for the sparsity threshold $\lambda$:
\begin{equation}\label{MSTLS2}
\text{MSTLS}(\Gbf,\bbf;\, \CalL,\pmb{\lambda})\qquad \begin{dcases} 
\hspace{0.1cm} \widehat{\lambda} = \min\left\{\lambda\in \pmb{\lambda} \ :\ \CalL(\lambda) = \min_{\lambda\in \pmb{\lambda}} \CalL(\lambda)\right\}\\
\widehat{\wbf}  =\text{STLS}(\Gbf,\bbf;\,\widehat{\lambda}).\end{dcases}
\end{equation}
Here $H_\lambda$ is the hard thresholding operator  defined by
\begin{equation}\label{hardthresh}
H_\lambda(\wbf)_i=\begin{cases} \wbf_i, & |\wbf_i|\geq \lambda \\ 0, & \text{otherwise}\end{cases}
\end{equation}
and $\CalL$ is the auxiliary loss function, introduced in \cite{messenger2020weakpde} and defined by
\begin{equation}\label{lossfcn}
\CalL(\lambda) = \frac{\nrm{\Gbf(\wbf^\lambda-\wbf^0)}_2}{\nrm{\Gbf\wbf^0}_2}+\frac{\nrm{\wbf^\lambda}_0}{\mathfrak{J}},
\end{equation}
which is defined on outputs $\wbf^\lambda = $ STLS$(\Gbf,\bbf;\,\lambda)$ of the sequential thresholding least squares algorithm with threshold $\lambda$. The collection of candidate sparsity thresholds $\pmb{\lambda}$ is chosen by the user (a finite set of equally log-spaced values of $\lambda$ is seen to be a successful choice for $\pmb{\lambda}$ in \cite{messenger2020weakpde}).

Altogether, the hyperparameters of the WSINDy algorithm are the reference test function $\psi$, the query points $\CalQ$, the model library composed of trial functions $\CalF$ and differential operator indices $\pmb{\alpha}$, and the sparsity thresholds $\pmb{\lambda}$. These are collected in Table \ref{hypparam}. 

\begin{rmrk}
In the original formulation of MSTLS in \cite{messenger2020weakpde}, the authors employed a different hard thresholding operator than \eqref{hardthresh} to incorporate relative term magnitudes $\nrm{\Gbf_i\wbf_i}_2/\nrm{\bbf}$ and a rescaling step to improve conditioning. Here we have chosen to present results using \eqref{hardthresh} in order to simplify the analysis and identify explicit conditions for convergence, however we conjecture that the original formulation of MSTLS is more computationally advantageous, especially for support recovery from multiscale data. We elaborate on the practical advantages and possible theoretical conclusions regarding the original formulation of MSTLS in Section \ref{discussion}. 
\end{rmrk}

\begin{table}
\centering
\begin{tabular}{|c|c|c|}
\hline Hyperparameter & Domain & Description \\ 
\hline $\psi$ & $C^{|\pmb{\alpha}|}(\Rbb^d+1)$ & test function \\ 
\hline $\CalQ := \{(\xbf_k,t_k)\}_{k\in[K]}$ & $\Rbb^{K\times(d+1)}$ satisfying \eqref{IBPbcs}& convolution query points  \\ \hline $\CalF := (f_j)_{j\in[J]}$ & $C(\Rbb)$ & trial functions\\ 
\hline $\pmb{\alpha} = (\alpha_s)_{s=0,\dots,S}$ & $\Nbb^{(S+1) \times (d+1)}$ & partial derivatives \\ \hline $\pmb{\lambda}$ & $[0,\infty)$ & candidate sparsity thresholds \\ \hline
\end{tabular}
\caption{Hyperparameters for the WSINDy Algorithm. Note that $\psi\in C^{|\pmb{\alpha}|}$ and $\CalQ$ satisfying \eqref{IBPbcs} ensure that the convolutional weak form \eqref{conv_form} is well defined. Here $|\pmb{\alpha}| := \max_{\alpha^s\in \pmb{\alpha}}\nrm{\alpha^s}_\infty$ is the maximum derivative order in library.}
\label{hypparam}
\end{table}

\subsection{WSINDy in the continuum limit: definitions and notation} \label{sec:contprob}

In this work we analyze the performance WSINDy in the limit of continuum data, in the sense that the solution $u$ is sampled on computational grids at finer and finer scales (to be made precise below), all while keeping the spatiotemporal domain $\Omega\times (0,T)$ fixed and while fixing the hyperparameters in Table \ref{hypparam}. As we will see in Section \ref{sec:controllablebias}, this continuum limit leads to a linear system $(\overline{\Gbf},\overline{\bbf})$ (defined in \eqref{contsys}), that is in general {\it biased} from the noise-free continuum linear system, denoted by $(\Gbf^\star,\bbf^\star)$, which has entries given analytically by either side of equation \eqref{conv_form}. By analyzing this biased system, we are able to prove that the bias is controllable, in that we prove conditions under which WSINDy still recovers the true support $S^\star := \supp{\wstar}$ of the true model coefficients $\wstar$ in the limit (see Section \ref{nosmooth}). We let the least squares solution to $(\overline{\Gbf},\overline{\bbf})$ be denoted by $\overline{\wbf}^0$, which under mild assumptions on the measurement noise satisfies $\overline{\Gbf}\overline{\wbf}^0=\overline{\bbf}$.

Specifically, we consider a sequence of samples $\{\Ubf^{(m)}\}_{m=1}^\infty$ on a dense set of successively finer computational grids $\{(\Xbf^{(m)}, \tbf^{(m)})\}_{m=1}^\infty\subset \Omega\times(0, T)$, each of which is equally spaced with resolution $(\Delta x^{(m)}, \Delta t^{(m)})$. As in Section \ref{sec:wsindyoverview}, we assume an i.i.d.\ additive noise model $\Ubf^{(m)} = u(\Xbf^{(m)}, \tbf^{(m)})+\ep$, detailed assumptions  of which are specified in Section \ref{sec:sampling}. With the hyperparameters in Table \ref{hypparam} fixed, for each $m$ we let $(\Gbf^{(m)},\bbf^{(m)})$ denote the linear system associated with WSINDy at discretization level $m$. A notable (if unsurprising) result of this paper is that under the assumptions in Section \ref{sec:assumptions}, there exists a {\it continuum linear system} 
\begin{equation}\label{contsys}
(\overline{\Gbf},\overline{\bbf}):=\lim_{m\to \infty}(\Gbf^{(m)},\bbf^{(m)})
\end{equation}
with convergence in probability and exhibiting an exponential concentration rate. We let $(\Gbf^\star,\bbf^\star)$ be the noise-free continuum linear system such that 
\begin{equation}\label{nzfreesys}
\bbf^\star = \Gbf^\star\wstar
\end{equation}
(i.e.\ the entries of $(\Gbf^\star,\bbf^\star)$ are given on either side of equation \eqref{conv_form}). We will refer to solving \eqref{contsys} for $\overline{\wbf}$ such that $\overline{\Gbf}\overline{\wbf}\approx \overline{\bbf}$ as the {\it continuum problem} and to solving \eqref{nzfreesys} for $\wstar$ as the {\it noise-free problem}. Throughout, weight vectors $\overline{\wbf}^\lambda$ refer to sequential thresholding least squares solutions to \eqref{contsys} with threshold $\lambda$. In other words, 
\begin{equation}\label{wlam}
\overline{\wbf}^\lambda = \text{STLS}(\overline{\Gbf},\overline{\bbf}; \lambda),
\end{equation}
and $\wbf^{(m),\lambda}$ is defined similarly with regard to $(\Gbf^{(m)},\bbf^{(m)})$. In particular, $\overline{\wbf}^0$ is the least-squares solution to the continuum problem and $\wbf^{(m),0}$ is the least-squares solution to the discrete system at level $m$. 

\subsection{Controllable bias}\label{sec:controllablebias}

In this section we informally discuss in what sense the bias in the continuum problem is controllable, such that WSINDy applied to $(\overline{\Gbf},\overline{\bbf})$ still recovers $\supp{\wstar}$. This is crucial to proving that WSINDy applied to the discrete systems $(\Gbf^{(m)},\bbf^{(m)})$ also recovers $\supp{\wstar}$ with probability rapidly approaching one as $m\to \infty$. To derive explicit results, we focus on the case of nonlinearities $\CalF$ containing only polynomial and trigonometric functions, although similar results may be available for more general libraries. Favoring a more intuition-building presentation, we introduce concepts in this section in the setting of ordinary differential equations, saving more general results to sections \ref{matcon} and \ref{sec:consistency}.

In terms of the continuum linear system and noise-free continuum linear system, for a fixed candidate weight vector $\wbf$ it holds that 
\[\lim_{m\to \infty}\left(\bbf^{(m)} - \Gbf^{(m)} \wbf\right) = \overline{\bbf}-\overline{\Gbf}\wbf=\Gbf^\star(\wstar-\wbf)+(\Gbf^\star-\overline{\Gbf})\wbf+(\overline{\bbf} - \bbf^\star).\]
Hence, we see that recovering $\wstar$ asymptotically ultimately depends on the proximity of $(\overline{\Gbf}, \overline{\bbf})$ to $(\Gbf^\star, \bbf^\star)$. A main result of the current manuscript is to show that although the continuum system is in general biased from the noise-free continuum system (i.e.\ $(\overline{\Gbf}, \overline{\bbf})\neq (\Gbf^\star, \bbf^\star)$), asymptotic recovery of supp$(\wstar)$ is still guaranteed for polynomial and trigonometric systems {\it for noise levels falling below a critical noise threshold}. Below we provide an explicit characterization of this threshold in the case of Gaussian noise, which depends on the relative magnitudes of the true coefficients $\wstar$ and the growth rate of nonlinearities present in the true model  (see the bounds \eqref{explicitboundsonsigmac}). Furthermore, we show in Section \ref{smooth} that if we modify the WSINDy algorithm (as presented in Section \ref{sec:wsindyoverview}) to additionally include a filtering step, then we recover the true support supp$(\wstar)$ asymptotically for all problems satisfying the assumptions in Section \ref{sec:assumptions}.

We will now explicitly characterize the continuum linear system $(\overline{\Gbf},\overline{\bbf})$ and its implications. The focus of Sections \ref{matcon} and \ref{sec:consistency} is to derive the following results in the PDE setting, together with exponential concentration bounds. However, for the sake of building intuition, consider data $\Ubf = u(\tbf)+\ep$ where $u:\Rbb\to\Rbb$ is a function of time satisfying some ODE $\frac{d}{dt}u = F(u(t))$, and as before $\ep$ represents i.i.d.\ measurement noise. Entries of the discrete linear systems $(\Gbf^{(m)},\bbf^{(m)})$ (defined in general by equation \eqref{conv_disc}) then consist of discretized integrals of the form
\begin{equation}\label{Im}
T^{(m)} := \sum_{i=1}^{m} \varphi(t_i) f(\Ubf_i)\Delta t^{(m)},
\end{equation}
where $\varphi$ denotes an arbitrary derivative of the reference test function $\psi$ and $f\in \CalF$ is a given function (possibly nonlinear). Letting $\rho$ be the distribution of the measurement noise $\ep$, under mild assumptions on $f$ and $\rho$ it holds that\footnote{We define the cross-correlation $f\star\rho(u) = \int_{\Rbb}f(u+x)d\rho(x)$.
This is equivalent to a convolution when $\rho$ is symmetric.} 
\begin{equation}\label{EIm}
\lim_{m\to \infty}\Ebb[T^{(m)}] = \lim_{m\to \infty}\sum_{i=1}^{m} \varphi(t_i) f\star\rho(u(t_i))\Delta t^{(m)} = \int_{\supp{\varphi}}  \varphi(t) f\star\rho(u(t))\,dt
\end{equation}
and 
\[\lim_{m\to \infty}\Vbb[T^{(m)}] = \lim_{m\to \infty}\Delta t^{(m)}\left(\sum_{i=1}^{m} \left(\varphi(t_i)\right)^2\Big(f^2\star\rho(u(t_i))-(f\star\rho(u(t_i)))^2\Big) \Delta t^{(m)}\right) =  0.\]
In other words, all entries of the linear system $(\Gbf^{(m)},\bbf^{(m)})$ converge to deterministic quantities of the form \eqref{EIm}, where $f$ has been replaced by $f\,\star\,\rho$, representing a bias between the continuum and noise-free problems.

If the trial function library $\CalF$ is such that $\{f\star\rho\ |\ f\in \CalF\} \in \text{span}(\CalF)$, then a linear transformation exists between the continuum least-squares solution $\overline{\wbf}^0$ and true solution $\wstar$, despite the inherent bias. This turns out to be the case for polynomial and trigonometric $f$, in which case for moderate noise levels and reasonably well-behaved dynamics, recovery of the correct model in the limit of large $m$ follows from a single round of thresholding. This partially explains the observed robustness to noise in \cite{messenger2020weakpde}. We now informally derive these linear transformations for polynomial and trigonometric libraries. 

When $f$ is polynomial, the bias $f\star\rho - f$ is a polynomial {\it of lower degree}. Specifically, with $f(x) = x^k$, we have that
\[f\star\rho(x) = \sum_{j=0}^k{k \choose j}M_{k-j}(\rho)x^j = f(x) + \sum_{j=0}^{k-1}{k \choose j}M_{k-j}(\rho)x^j,\]
where $M_j(\rho)$ denotes the $j$th moment of the distribution $\rho$. For example if $\rho$ is a white Gaussian noise distribution, then\footnote{The double factorial for an integer $n$ is defined $n!!:=n(n-2)\cdots(1)$, with $(0)!!=(-1)!!=1$.}
\[M_j(\rho) = \begin{dcases} (j-1)!!\sigma^j, & j \text{ even} \\0, & \text{otherwise.}\end{dcases}\]
The monomial library $P^{(q)} = (1,x,\dots,x^q)$ thus transforms under cross-correlation with Gaussian $\rho$ as
\[P^{(q)}\star\rho = P^{(q)}\Abf^{(q)},\] 
where $\Abf^{(q)}$ is defined
\begin{equation}\label{Agauss}
\Abf^{(q)}_{i,j} = \begin{dcases} {j\choose i}(j-i-1)!!\sigma^{j-i}, & j\geq i, \ (j-i)\ \text{even}  \\ 0, &\text{otherwise.}\end{dcases}
\end{equation}
In words, $\Abf^{(q)}$ is upper triangular with $1$'s along the diagonal and $0$'s along odd superdiagonals. The least-squares solution to the continuum problem is then given by $\overline{\wbf}^0 = (\Abf^{(q)})^{-1}\wstar$, and using that (see Lemma \ref{Ainv} in Appendix \ref{app:comblemms})
\[\left(\Abf^{(q)}\right)^{-1}_{i,j} = (-1)^{\frac{j-i}{2}} \Abf^{(q)}_{i,j},\]
the resulting coefficient error obeys the bound,
\[\nrm{\wstar-\overline{\wbf}^0}_\infty \leq C\sigma^2\nrm{\wbf^\star}_\infty.\]
The constant $C = C(p_{\max},\sigma)$ depends on the maximum degree monomial $p_{\max}$ in the true model, as well as the noise variance $\sigma^2$, and takes modest values for small $p_{\max}$: 
\[C(1,\sigma) = 0,\quad C(2,\sigma) = 1,\quad C(3,\sigma) = 3,\quad C(4,\sigma) = 6+3\sigma^2,\quad C(5,\sigma) = 10+15\sigma^2.\] 
Hence, if $p_{\max}$ is not too large and the noise variance $\sigma^2$ is moderate, then spurious terms in $\overline{\wbf}^0$ will be removed by a single round of thresholding.

The situation is even nicer for trigonometric functions. With $f=e^{i\omega x}$, we have 
\[f\star\rho(x) = \int_\Rbb e^{i\omega(x+y)}\rho(y)dy = e^{i\omega x} \widehat{\rho}(\omega) = \hat{\rho}(\omega) f(x)\]
where $\widehat{\rho}$ is the Fourier transform $\rho$ up to scaling. A library of trigonometric terms $F^{\pmb{\omega}} = (e^{i\omega\cdot(\cdot)})_{\omega\in \pmb{\omega}}$ then satisfies
\[F^{\pmb{\omega}}\star\rho = F^{\pmb{\omega}}\Dbf^{\pmb{\omega}}\]
where $\Dbf^{\pmb{\omega}}$ is diagonal. In this case, under mild restrictions on $\widehat{\rho}$, solving the continuum least squares problem produces $\overline{\wbf}^0$ with supp$(\overline{\wbf}^0) =$ supp$(\wstar)$ and 
\[\nrm{\wstar-\overline{\wbf}^0}_\infty \leq \max_{\omega\in \pmb{\omega}^\star}|1-\hat{\rho}(\omega)|\nrm{\wstar}_\infty\]  
where $\pmb{\omega}^\star\subset \pmb{\omega}$ is the set of trigonometric frequencies present in the true model. In the case of Gaussian white noise we have
\[\widehat{\rho}(\omega) = \exp\left(-\frac{\sigma^2\omega^2}{2}\right),\]
so that for $\sigma\leq \frac{0.14}{\nrm{\pmb{\omega}^\star}_\infty}$, the vector $\overline{\wbf}^0$ will be 99\% accurate. 

\begin{rmrk}\label{rmrkonomega}
For some noise distributions $\rho$ it is possible to have $\supp{\overline{\wbf}^0} \subsetneq \supp{\wstar}$ on trigonometric functions if $\widehat{\rho}(\omega)=0$ at some $\omega\in \pmb{\omega}$. Consider the uniform distribution $\rho = (2a)^{-1}\ind{[-a,a]}$. Then $\widehat{\rho}(n\pi/a) = 0$ for $n\in \Zbb$, which leads to the restriction $\nrm{\pmb{\omega}}_\infty<\frac{\pi}{a}$ in order for $\Dbf^{\pmb{\omega}}$ to be invertible. In this case the maximum allowable frequency is inversely proportional to the standard deviation $\sigma = a/\sqrt{3}$, which is not unreasonable: perturbations of the solution $u$ that are comparable to the periods of trigonometric terms will render their frequencies unobservable.  
\end{rmrk}

\begin{rmrk}
An understandable concern of using a weak-form recovery method is the biased recovery outlined in this section. However, we emphasize that the weak form is still advantageous over many strong-form methods which rely on point-wise derivative estimates. For a finite difference derivative operator $D$, the mean-squared error in the noisy derivative approximation of a one-dimensional derivative $u'(t)$ using data $\Ubf = u(\tbf) +\ep$ sampled at resolution $\Delta t$ satisfies 
\[\Ebb_{\ep \sim \rho}[|u'(t)-D\Ubf(t)|^2]=\CalO\left(\frac{\sigma^2}{\Delta t^2}\right),\]
where $\sigma^2 = \Vbb[\rho]$. In this form, there is no hope of convergence as $\Delta t\to 0$ for any fixed $\sigma>0$ unless the data is suitably denoised as $\Delta t$ is brought to zero. As we will see in this work, weak-form recovery methods do not suffer from instabilities as $\Delta t\to 0$.
\end{rmrk}

\subsection{Assumptions}\label{sec:assumptions}

Here we list the main assumptions used in the results below, along with a brief overview of their motivations. Specific results will specify which of the following assumptions are in effect.

\subsubsection{Regularity of the true solution}\label{sec:regularity}

In order to present results in as general a context as possible with respect to solution regularity, we specify the regularity of the true solution $u$ by defining the following function spaces. For open, bounded $D\subset \Rbb^{d+1}$, and for $p\in[1,\infty]$ and $k>0$, define
\begin{equation}\label{eq:Hspaces}
\CalH^{k,p}(D) := \left\{ f\in  L^p(D)\ :\ \exists \text{ disjoint, open } (D_i)_{i=1}^{\ell}\ \text{ s.t. }\ \overline{D} = \bigcup_{i=1}^\ell\overline{D}_i\,,\ f\big\vert_{D_i}\in H^k(D_i), \ \partial D_i\in C^{0,1}\right\},
\end{equation}
where $H^k(D)$ is the space of functions on $D$ with weak derivatives up to order $k$ in $L^2(D)$. The spaces $\CalH^{k,p}(D)$ are similar to the broken Sobolev spaces used in the analysis of discontinuous Galerkin finite element methods (see e.g.\ \cite{houston2002discontinuous}). We assume that $u\in \CalH^{k,\infty}(\Omega\times (0,T))$ for some $k > (d+1)/2$ is a weak solution to \eqref{diffform} with coefficients $\wbf=\wstar$ and that any points of reduced regularity (e.g.\ discontinuities) are contained to the boundaries between subdomains $D_i$ of $D:=\Omega\times (0,T)$. The restriction $k > (d+1)/2$ ensures by the Sobolev embedding theorem that $u$ is bounded and piecewise H\"older continuous on each $D_i$. We show in Appendix \ref{app:trapconv}
 that this is sufficient regularity for the trapezoidal rule to converge for integrals of smooth functions of $u$.

\subsubsection{Sampling model}\label{sec:sampling}
We assume that every pointwise evaluation\footnote{With $u\in \CalH^{k,\infty}(\Omega\times (0,T))$ such that $k>(d+1)/2$, we have that pointwise evaluations of $u$ are well-defined (apart from a set of measure zero, e.g.\ when considering solutions with shocks) by the Sobolev embedding theorem (see Appendix \ref{app:trapconv} for more details).}  of the solution $u$ (resulting in the samples $\{\Ubf^{(m)}\}_{m=1}^\infty$ for each grid $\{(\Xbf^{(m)},\tbf^{(m)})\}_{m=1}^\infty$) produces a measurement error $\ep \sim \rho$, where $\rho$ is symmetric and sub-Gaussian, that is, 
\[\nrm{\rho}_{\text{SG}}:=\inf\{\lambda >0 \,:\, \Ebb_{\ep\sim \rho}\left[\exp(\ep^2/\lambda^2)\right]\leq 2\}<\infty.\]  
In particular, this includes Gaussian white noise and bounded noise. In addition, we assume the noise at distinct spatiotemporal points is uncorrelated: $\Ebb[\ep(x,t)\ep(y,s)] = \sigma^2\ind{(x,t)=(y,s)}$ for two points $((x,t),(y,s))\in (\Omega\times[0,T])\times(\Omega\times[0,T])$ and fixed variance $\sigma^2$. We refer to $\sigma$ as the {\it noise level} and note that $\sigma\leq \sqrt{2}\nrm{\rho}_{\text{SG}}$ (see the textbook \cite{vershynin2018high} for more details). 

Furthermore, we assume that each grid $(\Xbf^{(m)},\tbf^{(m)})$ has equal resolution $(\Delta x^{(m)},\Delta t^{(m)})$, and throughout the discretization level $m$ is defined as the number of points that the reference test function $\psi$ is supported on at the resolution $(\Delta x^{(m)},\Delta t^{(m)})$, or
\begin{equation}\label{def_of_m}
m := \#\{\supp{\psi}\cap (\Xbf^{(m)},\tbf^{(m)})\}
\end{equation}
where $\#\{\cdot\}$ indicates the set cardinality.

\subsubsection{Model library}\label{sec:modellibrary}

We assume the collection of multi-indices $\pmb{\alpha} = (\alpha^0,\dots,\alpha^S)$ is known and the family of functions $\CalF = (f_j(u))_{j\in[J]}$ consists of $P^{(p_{\max})}$, the space of polynomials of degree at most $p_{\max}$ on $\Rbb^n$, as well as $F^{\pmb{\omega}} = \{\exp(i\omega^T u)\}_{\omega \in \pmb{\omega}}$, a finite collection of Fourier modes on $\Rbb^n$ (i.e. $\pmb{\omega}\subset \Rbb^n$).

\subsubsection{Reference test function}\label{sec:testfcn}
We assume that $\psi\in C^{|\pmb{\alpha}|}(\Omega\times (0,T))$ with compact support in $\Omega\times (0,T)$. When $\psi$ is taken to be separable, $\psi(x) = \phi_1(x_1)\cdots\phi_d(x_d)\phi_{d+1}(t)$, it is assumed that each $\phi_i\in C^{|\pmb{\alpha}|}(\Rbb)$.

\subsubsection{Conditioning of the noise-free system}\label{sec:conditioning}

In this work we ensure that the recovery problem is solvable by assuming that the noise-free continuum matrix $\Gbf^\star$ is full-rank. This implies several assumptions:
\begin{enumerate}[label=(\Roman*)]
\item  The underlying PDE has a unique representation of the form \eqref{diffform} with coefficients $\wstar$ over the library $(\CalF,\pmb{\alpha})$. 
\item The reference test function $\psi$, library $(\CalF,\pmb{\alpha})$, and solution $u$ are such that the set of vectors (with Query points $\CalQ$ arranged into a vector)
\[\left(\partial^{\alpha^s}\psi*f_j(u)(\CalQ)\right)_{s=1,\dots,S,\ j=1,\dots,J}\]
are linearly independent. In particular, the number of convolution query points is not smaller than the library size: $K\geq \mathfrak{J}$, where $K$ is the number of query points, $S$ is the number of differential operators in $\pmb{\alpha}$, $J$ is number of functions $\CalF$, and $\mathfrak{J} = SJ$ is the total number of terms in the model library. 
\end{enumerate}

\begin{rmrk}
The most restrictive assumption above is the full rank assumption on $\Gbf^\star$ in \ref{sec:conditioning}. This implies that there does not exist $\tilde{\wbf}\neq 0$ such that $u$ satisfies an additional constraint of the form
\[0 = \sum_{s,j}\tilde{\wbf}_{s,j}\partial^{\alpha^s}f_j(u).\]   
Such constraints do indeed arise in practice (e.g.\ the divergence-free constraint in the incompressible Navier-Stokes equations). While this cannot be directly checked in practice from noisy data sets, the condition number of the weak-form linear system provides an excellent guide for the conditioning of the underlying noise-free system, and is reflected in the choice of reference test function $\psi$, query points $\CalQ$, and library $(\CalF,\pmb{\alpha})$, all of which the user has control over.

Ultimately, the full-rank assumption is due to our use of the STLS algorithm for performing sparse regression, however this can be relaxed if other algorithms are chosen, and extensions of STLS to the rank-deficient case are possible. However, due to spurious terms arising in the continuum limit, we conjecture that some form of thresholding-based sparse regression is advantageous (over e.g.\ greedy methods or $\ell_1$-regularization). 
\end{rmrk}

\section{Exponential concentration of the WSINDy linear system}\label{matcon}

Using concentration results for sums of heavy-tailed random variables recently proved in \cite{vladimirova2020sub,bakhshizadeh2020sharp}, we prove that $\nrm{\Gbf^{(m)}-\overline{\Gbf}}_\infty\to 0$ with an exponential convergence rate. We will need the following lemmas to connect existing results with the specific form of the entries of the matrix $\Gbf^{(m)}$. Proofs can be found in Appendix \ref{app:conclemmas}.

Using the terminology of \cite{bakhshizadeh2020sharp}, the following lemma establishes that the right tails of summands within each entry of $\Gbf^{(m)}$ can be modelled using a common rate function.

\begin{lemm}\label{lemm1}
Let $f:\Rbb\to \Rbb$ satisfy
\begin{equation}\label{eq:growthcondition}
|f(x)|\leq C_f\left(1+|x|^p\right),
\end{equation}
for some $p\geq 1$ and $C_f>0$, and let $\rho$ be a symmetric sub-Gaussian distribution. For bounded sequences of real numbers $(\alpha_i)_{i\in \Nbb}$ and $(u_i)_{i\in \Nbb}$, and for $\ep_i\sim \rho$ {\normalfont i.i.d.}\, $i\in \Nbb$, define the random variables $Y_i = \alpha_i f(u_i +\ep_i)$. Then for any $\kappa\geq \nrm{\rho}_{\text{SG}}$ and $t>0$, it holds that the right tails of $Y_i$ are captured by a common rate function $I(t)$,
\[\Pbb\left(Y_i > t\right)\leq \exp (-I(t)),\]
where 
\begin{equation}\label{eq:rate_fcn}
I(t):=\ind{(t^*,\infty)}(t)\left[\frac{1}{\kappa^2}\left(\left(\frac{t}{C_f\alpha^*}-1\right)^{1/p}-u^*\right)^2-\log(2)\right] = \frac{t^{2/p}}{\kappa^2(C_f\alpha^*)^{2/p}}I_0(t),
\end{equation}
for $\alpha^* = \sup_i|\alpha_i|$, $u^* = \sup_i|u_i|$, and $t^*:= C_f\alpha^*\left(1+\left(u^*+ \kappa\sqrt{\log(2)}\right)^p\right)$. Moreover, $I_0(t)$ is monotonically increasing from $0$ to $1$ over $t\in(t^*,\infty)$, and is defined in the proof in Appendix \ref{app:conclemmas}.
\end{lemm}
The next lemma allows one to uniformly bound a sequence of independent, non-identically distributed random variables of the form in Lemma \ref{lemm1}.
\begin{lemm}\label{lemm2}
Let $Y_i$ be defined under the same conditions as Lemma \ref{lemm1} and choose $\beta\in(0,1)$. Then there exists $\overline{v}(\beta) < \infty$ such that the sum $S_m = \sum_{i=1}^{m} Y_i$ satisfies
\begin{equation}\label{eq:convergence_of_sums}
\Pbb\left(|S_m-\Ebb S_m|> mt\right) \leq  \begin{dcases} 2 \exp\left(-\frac{\beta}{2} I(mt)\right) + 2m\exp\left(-I(mt)\right), & t \geq t_m(\beta)\\
2 \exp\left(-\frac{mt^2}{2\overline{v}(\beta)}\right) + 2m\exp\left(-\frac{mt_m(\beta)^2}{\overline{v}(\beta)}\right), & 0\leq t < t_m(\beta), \end{dcases}
\end{equation}
where $t_m (\beta):= \sup\{t\geq 0\ :\ t\leq \beta\overline{v}(\beta)\frac{I(mt)}{mt}\} \to 0$ in m.
\end{lemm}

We now present a main result concerning concentration of the WSINDy discrete linear system to its continuum variant. Throughout we use the notation $\nrm{\Gbf}_{\vec{p}}$ to denote the vector $p$ norm of matrix $\Gbf$ stretched into a column vector, and the notation $(\Gbf,\bbf)$ to mean the concatenation of matrix $\Gbf\in \Rbb^{K\times \mathfrak{J}}$ with vector $\bbf\in \Rbb^K$.

\begin{thm}\label{matrix_concentration}
Suppose each function in library $\CalF$ satisfies the growth bound \eqref{eq:growthcondition} for some $p:=p_{\max}$ and Assumptions \ref{sec:regularity}, \ref{sec:sampling}, and \ref{sec:testfcn} hold.  Then it holds that for every $t>\overline{t}(m)$, where $\overline{t}(m)\to 0$ as $m\to \infty$, we have the concentration rates
\begin{equation}\label{eq:convergence_of_G}
\Pbb\Big(\nrm{(\Gbf^{(m)},\bbf^{(m)}) - (\overline{\Gbf},\overline{\bbf})}_{\overrightarrow{\infty}} > t \Big) \leq  \begin{dcases} K\mathfrak{J}\exp\left(-\frac{c}{2} (mt)^{2/p_{\max}}\right) + K\mathfrak{J}m\exp\left(-c (mt)^{2/p_{\max}}\right), & t \geq t_m\\
K\mathfrak{J} \exp\left(-\frac{mt^2}{2\overline{v}}\right) + K\mathfrak{J}m\exp\left(-\frac{mt_m^2}{\overline{v}}\right), & 0\leq t < t_m, \end{dcases}
\end{equation}
where the rate factor $c$ depends on $\nrm{u}_\infty$, $|\Omega\times (0,T)|$, $\pmb{\alpha}$, $\psi$, $\CalF$, and $\nrm{\rho}_{\text{SG}}$, and $\overline{v} = \overline{v}(1/2)$ and $t_m=t_m(1/2)$ from Lemma \ref{lemm2}.\end{thm}
For more details on the rate $c$ in Theorem \ref{matrix_concentration}, see the proof in Appendix \ref{app:conclemmas}. In addition, to make more straight-forward use of these concentration results, we have the following.

\begin{corr}\label{matrix_corr}
Under the assumptions of Theorem \ref{matrix_concentration}, for every $t>0$ and sufficiently large $m$ it holds that 
\begin{equation}\label{eq:convergence_of_G_2}
\Pbb\Big(\nrm{(\Gbf^{(m)},\bbf^{(m)}) - (\overline{\Gbf},\overline{\bbf})}_{\overrightarrow{\infty}} > t \Big) \leq  2K\mathfrak{J}\exp\left(-\frac{c}{2}(mt)^{2/p_{\max}}\right).
\end{equation}
\end{corr}
\begin{proof}
This comes from noting that for every $t$ and sufficiently large $m$ we have 
\[\frac{c}{2}(mt)^{2/p_{\max}} < \min\left(c(mt)^{2/p_{\max}}-\log(m),\ \frac{mt^2}{2\overline{v}}-\log(m)\right).\] 
\end{proof}

\begin{rmrk}
The proof of Theorem \ref{matrix_concentration} reveals several ways that the rate of concentration can be increased. The most effective is to lower the growth-rate $p_{\max}$ and decrease $\nrm{\partial^{\alpha^s}\psi}_\infty$. This implies that in practice, concentration will be determined by the growth rate of nonlinearities and order of differential operators in the true model. In addition, if the true solution $u$ has increased smoothness, $\overline{t}(m)$ goes to zero much faster, leading to faster entry into the regime of exponential concentration. Lastly, $c$ is inversely proportional to $\nrm{\rho}_{\text{SG}}^2$, hence decreasing the variance of the noise directly increases the concentration rate.  
\end{rmrk}

\section{Asymptotic consistency}\label{sec:consistency}

In this section we provide asymptotic consistency results in the form of support recovery of the true model coefficients $\supp{\wstar}$ in the limit as $m\to \infty$ (recall $m$ is the number of points that the reference test function $\psi$ is supported on at the resolution $(\Delta x^{(m)},\Delta t^{(m)})$). 

In Section \ref{nosmooth} we prove that WSINDy with a hyperparameter-free version of the MSTLS algorithm recovers the true model support $S^\star = \supp{\wstar}$ with high-probability from the raw (un-filtered) data as $m\to \infty$, provided $\sigma<\sigma_c$ for some critical noise level $\sigma_c$ and the assumptions in Section \ref{sec:assumptions} are met. This is done by first proving recovery results for the continuum problem, and then combining these with the matrix concentration results in Section \ref{matcon}. In order to demonstrate how explicit bounds on $\sigma_c$ may be derived, in Section \ref{nosmooth} we focus on proving conditions for {\it subset} support recovery, or $\supp{\widehat{\wbf}}\subset S^\star$, under the restricted setting of Gaussian white noise. A proof of full support recovery $\supp{\widehat{\wbf}} = S^\star$, requiring an additional assumption on $(\Gbf^\star,\bbf^\star)$ (see Remark \ref{promise_of_fullsupprec}), is presented in Appendix \ref{app:supprec} and may be extended to the case of arbitrary symmetric sub-Gaussian noise using Lemma \ref{generalmomentmatrix}. Furthermore, Lemma \ref{generalmomentmatrix} provides a general mechanism for deriving bounds on $\sigma_c$ specific to any symmetric sub-Gaussian noise distribution $\rho$. Explicit bounds on $\sigma_c$ may be highly informative for future algorithmic developments.

To summarize, we first prove in Lemma \ref{lemm:bounds on sigma_c} that for noise levels below some critical noise $\sigma_c'$ there exists a feasible $\widehat{\lambda}$ for which the STLS solution of associated continuum problem recovers the true support $S^\star$. We use properties of the Gaussian moment matrix to present explicit bounds $\sigma_c'$. We then show in Lemma \ref{lemm:subset_support_rec} that there exists $\sigma_c\leq \sigma_c'$ below which the one-shot MSTLS algorithm (defined below) applied to the continuum problem produces $\widehat{\wbf}$ satisfying $\supp{\widehat{\wbf}}\subset S^\star$. The utility of Lemma \ref{lemm:subset_support_rec} lies in the fact that the one-shot MSTLS algorithm involves no hyperparameters, hence avoids the task of selecting a feasible $\widehat{\lambda}$. We then combine these results with classical stability of the least-squares problem (Lemma \ref{least_squares_cont}) and the concentration results in Section \ref{matcon} to yield $\supp{\widehat{\wbf}^{(m)}}\subset S^\star$ with high probability for the MSTLS solution $\widehat{\wbf}^{(m)}$ on the linear system $(\Gbf^{(m)},\bbf^{(m)})$, provided $\sigma<\sigma_c$ and $m$ is large enough. This result is extended to yield $\supp{\widehat{\wbf}^{(m)}}= S^\star$ with high probability in Appendix \ref{app:supprec}, and for general symmetric sub-Gaussian noise with the aid of Lemma \ref{generalmomentmatrix}.

In Section \ref{smooth}, we consider the case of filtering the data before applying the WSINDy algorithm. Examining the case of simple moving average filters, which may easily be extended to a wider class of filters, we prove the exponential concentration of $(\Gbf^{(m)},\bbf^{(m)})$ to $(\Gbf^\star,\bbf^\star)$, in contrast to $(\overline{\Gbf},\overline{\bbf})$. This implies that for a wider class of libraries (only requiring $\CalF$ to be locally Lipschitz) and arbitrary symmetric sub-Gaussian noise, that we get $\supp{\widehat{\wbf}^{(m)}}\subset S^\star$ with high probability using the hyperparameter-free MSTLS algorithm. As before, this is strengthened to $\supp{\widehat{\wbf}^{(m)}}= S^\star$ if an additional condition on $(\Gbf^\star,\bbf^\star)$ is satisfied.

\subsection{Asymptotic consistency without filtering}\label{nosmooth}

Recall from Section \ref{sec:controllablebias} that for trigonometric and polynomial libraries we have $\text{span}(\CalF\star\rho) \subset \text{span}(\CalF)$, that is, the cross-correlated library terms are linear combinations of the original library terms. More specifically, under assumptions \ref{sec:regularity}-\ref{sec:testfcn} there exists an upper triangular, block diagonal matrix 
\begin{equation}\label{Ablock}
\Abf = \text{blkdiag}(\underbrace{\Abf^{(p_{\max})},\Dbf^{\pmb{\omega}}}_{s=1},\dots,\underbrace{\Abf^{(p_{\max})},\Dbf^{\pmb{\omega}}}_{s=S})
\end{equation}
where
\[F^{\pmb{\omega}}\star\rho = F^{\pmb{\omega}}\Dbf^{\pmb{\omega}}, \qquad P^{(p_{\max})}\star\rho = P^{(p_{\max})}\Abf^{(p_{\max})}\]
such that $\overline{\Gbf} = \Gbf^\star \Abf$. For all symmetric noise distributions $\rho$ we have 
\[\Dbf^{\pmb{\omega}} = \text{diag}(\hat{\rho}(\pmb{\omega})), \qquad \Abf^{(p_{\max})}_{ij} := \delta_{ij}+\Lbf_{ij}^{(p_{\max})}=\begin{dcases} {j\choose i} M_{j-i}(\rho)& 0\leq i\leq j \\ 0 & \text{otherwise,} \end{dcases}\]
provided the highest moment $M_{p_{\max}}(\rho)$ exists. 

We can decompose $\Abf$ into diagonal and off-diagonal matrices as follows, 
\begin{equation}\label{Ainv_DL}
\Abf = \Dbf + \Lbf
\end{equation}
where $\Dbf :=  \Ibf^{(S)}\otimes \text{blkdiag}(\Ibf^{(p_{\max}+1)},\Dbf^{\pmb{\omega}})$ is diagonal with diagonal entries less than or equal to 1 in magnitude (since $\rho$ is a probability distribution), and $\Lbf := \Ibf^{(S)}\otimes \text{blkdiag}(\Lbf^{(p_{\max})},\Ibf^{(|\pmb{\omega}|)})$ is zero along the diagonal. Furthermore, in the case that $\widehat{\rho}>0$, (e.g.\ when $\rho$ is Gaussian), we have 
\begin{equation}\label{invA}
\Abf^{-1} = \Dbf^{-1} + \tilde{\Lbf}
\end{equation}
where $\tilde{\Lbf} = \Ibf^{(S)}\otimes \text{blkdiag}(\tilde{\Lbf}^{(p_{\max})},\Ibf^{(|\pmb{\omega}|)})$, and $\tilde{\Lbf}^{(p_{\max})}$ is equal to $\Lbf^{(p_{\max})}$ up to sign changes (see Lemma \ref{Ainv} in Appendix \ref{app:comblemms}).

In the following we assume $\rho$ is Gaussian, as it allows for explicit bounds on the critical noise level $\sigma_c$. The polynomial moment matrix $\Abf^{(p_{\max})}$ is then given by \eqref{Agauss}, with inverse given in Lemma \ref{Ainv}, and $\Dbf^{\pmb{\omega}}$ is positive definite. Similar results, but not as explicit, exist for any other symmetric sub-Gaussian noise distribution (see Lemma \ref{generalmomentmatrix}) however restrictions may be needed on $\pmb{\omega}$ to ensure the invertibility of $\Dbf$ (see Remark \ref{rmrkonomega}). First we need the following existing result on the STLS algorithm.

\begin{lemm}\label{lemm:STLS_supp_rec}
    Let $\Abf$ be the corresponding moment matrix such that $\overline{\Gbf} = \Gbf^\star \Abf$ and let $\what =$ STLS$(\overline{\Gbf},\overline{\bbf},\widehat{\lambda})$ be the continuum STLS solution with sparsity threshold $\widehat{\lambda}$. A necessary and sufficient condition for one iteration of STLS to result in ${\normalfont\supp{\what} = \supp{\wstar}}$ is
\begin{equation}\label{necccond1}
\min_{j\in S^\star}\left\vert\left(\Abf^{-1}\wstar\right)_j\right\vert>\widehat{\lambda}> \max_{j\in (S^\star)^c}\left\vert\left(\Abf^{-1}\wstar\right)_j\right\vert.
\end{equation}
Moreover, \eqref{necccond1} is sufficient to ensure $\supp{\what} \subset \supp{\wstar}$ for any number of STLS iterations.
\end{lemm}

\begin{proof}
This is a special case of \cite[Proposition 2.8]{zhang2019convergence}, considering that $\Abf^{-1}\wstar = \overline{\Gbf}^\dagger\overline{\bbf}$.
\end{proof}

In the next lemma we classify which models and noise levels lead to existence of $\widehat{\lambda}$ satisfying \eqref{necccond1} in the case of Gaussian noise.

\begin{lemm}\label{lemm:bounds on sigma_c} 
Let $\rho$ be a Gaussian mean-zero noise distribution and $\Abf$ be the corresponding moment matrix such that $\overline{\Gbf} = \Gbf^\star \Abf$. Let $p$ be the maximum polynomial degree appearing in the true model and define $S^\star:=\supp{\wstar}$. Then we have the following cases:
\begin{enumerate}[label=(\roman*)]
\item If $p \leq 2$ and the true coefficient of $u^2$ is zero, then there exists $\widehat{\lambda}$ satisfying \eqref{necccond1} for any finite noise level $\sigma>0$.
\item If $p\geq 3$, or if $p=2$ and the true coefficient of $u^2$ is nonzero, then there exists a critical noise level $\sigma_c$ satisfying
\begin{equation}\label{explicitboundsonsigmac}
\left(\frac{1}{2{p \choose 2}e}\right)\frac{\min_{j\in S^\star}|\wstar_j|}{\max_{j \in S^\star}|\wbf^\star_j|}\leq \sigma_c^2\leq \frac{1}{{p \choose 2}},
\end{equation}
such that for all $\sigma < \sigma_c$, there exists $\widehat{\lambda}$ such that \eqref{necccond1} holds.
\end{enumerate}
\end{lemm}

\begin{proof}
For $p\leq 2$, if the term $u^2$ does not exist in the true model, in other words all terms $\partial^{\alpha^s}u^2$ in the true model satisfy $|\alpha^s|\geq 1$, then using that 
\begin{equation}\label{u2}
\partial^{\alpha^s}\psi*(u^2\star \rho) = \partial^{\alpha^s}\psi*(\sigma^2 + u^2) = \partial^{\alpha^s}\psi*u^2,
\end{equation}
we see that no spurious terms are generated (recall from Section \ref{sec:controllablebias} that trigonometric terms do not generate spurious terms in the continuum problem). In these cases any $\widehat{\lambda} < \min_{i\in S^\star}|\wstar_i|$ satisfies \eqref{necccond1}, using that $|(\Abf^{-1}\wstar)_i| = |\Dbf^{-1}_{ii}\wstar_i|\geq |\wstar_i|$ for all $i\in\{1,\dots,\mathfrak{J}\}$, since $\nrm{\text{diag}(\Dbf)}_\infty\leq 1$.

Now assume that $p\geq 3$ or $p=2$ with $u^2$ contained in the true model. To derive the upper bound in \eqref{explicitboundsonsigmac}, consider the case where all entries of $\wstar_{S^\star}$ have equal magnitude. Then there exists $\widehat{\lambda}$ such that \eqref{necccond1} holds only if 
\[\min_{j\in S^\star}\left\vert\left(\Abf^{-1}\wstar\right)_j\right\vert{p \choose 2} \sigma^2 < \min_{j\in S^\star}\left\vert\left(\Abf^{-1}\wstar\right)_j\right\vert,\]
since the coefficient of $\partial^{\alpha^s}\psi*u^p$ is nonzero for some $\alpha^s\in \pmb{\alpha}$, which generates the spurious term ${p \choose 2}\sigma^2\partial^{\alpha^s}\psi*u^{p-2}$. Equivalently, such a $\widehat{\lambda}$ exists only if $\sigma^2 < {p \choose 2}^{-1}$, which necessitates $\sigma_c^2\leq {p \choose 2}^{-1}$.

For the sufficient lower bound in \eqref{explicitboundsonsigmac}, using Lemma \ref{app:lemm_bnd_on_AL} to bound $\|\tilde{\Lbf}\|_\infty$, we have that 
\[\max_{j\in (S^\star)^c}\left\vert\left(\Abf^{-1}\wstar\right)_j\right\vert \leq \|\tilde{\Lbf}\|_\infty\nrm{\wbf^\star}_\infty \leq \sigma^2{p \choose 2}\exp\left(\sigma^2{p \choose 2}\right)\nrm{\wbf^\star}_\infty\]
and 
\[\min_{j\in S^\star}\left\vert\left(\Abf^{-1}\wstar\right)_j\right\vert \geq \min_{j \in S^\star}|\Dbf^{-1}\wbf^\star_j|-\|\tilde{\Lbf}\|_\infty\nrm{\wbf^\star}_\infty \geq \min_{j \in S^\star}|\wbf^\star_j| - \sigma^2{p \choose 2}\exp\left(\sigma^2{p \choose 2}\right)\nrm{\wbf^\star}_\infty.\]
A sufficient condition for existence of $\widehat{\lambda}$ satisfying \eqref{necccond1} is thus
\begin{equation}\label{suff_lower_bound}
2\sigma^2{p \choose 2}\exp\left(\sigma^2{p \choose 2}\right)< \frac{\min_{j\in S^\star}|\wstar_j|}{\max_{j \in S^\star}|\wbf^\star_j|}.
\end{equation}
Taking $\sigma \leq \sigma_c$ and using the upper bound $\sigma_c^2 \leq {p \choose 2}^{-1}$, we get that \eqref{suff_lower_bound} is implied by 
\[\sigma^2 < \frac{1}{2{p \choose 2}e}\frac{\min_{j\in S^\star}|\wstar_j|}{\max_{j \in S^\star}|\wbf^\star_j|},\]
hence by the sufficiency of this upper bound we achieve the lower bound on $\sigma_c$ in \eqref{explicitboundsonsigmac}.
\end{proof}

\begin{rmrk}
Lemma \ref{lemm:bounds on sigma_c} provides a rigorous explanation for the robustness to noise of WSINDy observed in \cite{messenger2020weakpde}, as several systems that were shown to yield robust recovery fall into the case $(i)$ for which $\sigma_c=\infty$. These include inviscid Burgers, Korteweg-de Vries, Kuramoto-Sivashinksy, porous medium, Sine-Gordon, and Navier-Stokes. Moreover, {\it all linear differential equations fall into case $(i)$}. Furthermore, several types of nonlinear PDEs including reaction-diffusion and nonlinear Schr\"odinger's (also examined in \cite{messenger2020weakpde}) fall into the case $(ii)$, for which spurious terms will arise in the limit of large data if the noise level is greater than some finite $\sigma_c$.
\end{rmrk}

\begin{rmrk}
Condition \eqref{necccond1} also implies that the STLS algorithm with any number of iterations (at least one) yields a solution $\supp{\what}$  with $\supp{\what} = \supp{\wstar}$, provided the noise level $\sigma^2$ is low enough. This follows from stability of least-squares problems in Lemma \ref{proj_cont} of Appendix \ref{app:supprec}, however we are unable to provide explicit bounds on the critical noise $\sigma_c$ in this case, hence we focus on the one-shot MSTLS algorithm in what follows.
\end{rmrk}

In Lemma \ref{lemm:bounds on sigma_c} we identified a relation in the Gaussian noise case between the noise variance $\sigma^2$ and the existence of a sparsity parameter $\widehat{\lambda}$ such that one step of STLS produces the correct support. In the next two lemmas we identify conditions under which the MSTLS algorithm \eqref{MSTLS2} identifies a feasible $\widehat{\lambda}$ satisfying \eqref{necccond1}, which leads to correct support in one step of STLS. While in practice we often run STLS until it terminates (in a maximum of $\mathfrak{J}$ iterations), for the sake of identifying conditions for convergence we focus on the simpler case of performing only one step of STLS for each candidate $\lambda\in\pmb{\lambda}$.

In the one-step STLS case there are only $\mathfrak{J}+1$ possible sparse solutions. To see this, order $\overline{\wbf}^0$ from least to greatest in absolute value:
\[|\overline{\wbf}^0_{(0)}|\leq |\overline{\wbf}^0_{(1)}|\leq |\overline{\wbf}^0_{(2)}|\leq \cdots\leq |\overline{\wbf}^0_{(\mathfrak{J})}|\]
where $\overline{\wbf}^0_{(0)} = 0$ is inserted. Then for all $\lambda \in (\overline{\wbf}^0_{(i)},\overline{\wbf}^0_{(i+1)})$, one step of STLS produces  the same solution, hence there are at most $\mathfrak{J}$ distinct solutions attainable for $\lambda \in [0,\nrm{\overline{\wbf}^0}_\infty]$, and $\lambda> \nrm{\overline{\wbf}^0}_\infty$ leads to the zero vector. With this in mind, we need only examine the case 
\begin{equation}\label{oneshotpmblambda}
\pmb{\lambda} = \left\{ \frac{|\overline{\wbf}^0_{(i)}|+|\overline{\wbf}^0_{(i+1)}|}{2}\ :\ i=\{0,\mathfrak{J}-1\}\right\}, 
\end{equation}
discarding duplicate values. In what follows, we let $\overline{\wbf}^\lambda = \text{STLS}^{(1)}(\overline{\Gbf},\overline{\bbf},\lambda)$ denote the STLS solution using a single round of thresholding, in other words,
\[\overline{\wbf}^\lambda = \argmin_{\supp{\wbf}\subset \supp{H_\lambda(\overline{\Gbf}^\dagger\overline{\bbf})}} \nrm{\overline{\Gbf}\wbf-\overline{\bbf}}_2^2\]
and $\widehat{\wbf} = \text{MSTLS}^{(1)}(\overline{\Gbf},\overline{\bbf})$ denote the MSTLS solution using one round of thresholding per inner STLS loop over the thresholds \eqref{oneshotpmblambda}.

\begin{lemm}\label{lemm:subset_support_rec}
Let $\rho$ be a Gaussian mean-zero noise distribution and let $p$ be the maximum polynomial degree appearing in the true model. Let $\Gbf^\star_p$ be the restriction of $\Gbf^\star$ to the columns corresponding to polynomial terms of maximum degree $p$. Then there exists a critical noise level $\sigma_c$ such that for any $\sigma\leq \sigma_c$ the estimator $\what = \text{MSTLS}^{(1)}(\overline{\Gbf},\overline{\bbf})$ satisfies 
\begin{equation}\label{subsetsupprec_lemm}
\supp{\what}\subset\supp{\wstar}
\end{equation}
where $\sigma_c$ satisfies
\begin{enumerate}[label=(\roman*)]
\item If $p\leq 2$ and the coefficient of $u^2$ is zero in the true model, then $\sigma_c=\infty$.
\item Otherwise, $\sigma_c$ satisfies the bounds
\begin{equation}\label{boundsMSTLS}
\frac{1}{e{p \choose 2}}\min\left\{\frac{\nrm{\Gbf^\star\wstar}_2}{\mathfrak{J}\nrm{\Gbf^\star_p}_2\nrm{\wstar}_2}, \frac{\min_{j\in S^\star}|\wstar_j|}{2\max_{j\in S^\star}|\wstar_j|}\right\}\leq \sigma_c^2\leq {p \choose 2}^{-1}.
\end{equation}
\end{enumerate}
\end{lemm}
\begin{proof}
The MSTLS loss is defined
\[\CalL(\lambda) = \frac{\nrm{\overline{\Gbf}(\overline{\wbf}^0 -\overline{\wbf}^\lambda)}}{\nrm{\overline{\Gbf} \overline{\wbf}^0}}+\frac{\nrm{\overline{\wbf}^\lambda}_0}{\mathfrak{J}}\]
where $\overline{\wbf}^\lambda=\text{STLS}^{(1)}(\overline{\Gbf},\overline{\bbf},\lambda)$. We first assume that there exists $\widehat{\lambda}$ satisfying \eqref{necccond1}, which is given by $\sigma<\sigma_c$ for $\sigma_c$ satisfying \eqref{explicitboundsonsigmac}. Then it holds that $\widehat{\lambda}\in \pmb{\lambda}$ given by \eqref{oneshotpmblambda}. Also, for $\lambda \in \pmb{\lambda}$ satisfying $\lambda>\min_{j\in S^\star}|(\Abf^{-1}\wstar)_j|$, it holds that $\supp{\overline{\wbf}^\lambda}\subset\supp{\wstar}$, so it suffices to prove that $\CalL(\tilde{\lambda})>\CalL(\widehat{\lambda})$ for all $\tilde{\lambda} \leq \max_{j\in (S^\star)^c}|(\Abf^{-1}\wstar)_j|$.

Indeed, choose $\tilde{\lambda}\leq \max_{j\in (S^\star)^c}|(\Abf^{-1}\wstar)_j|$ whereby
\[S^\star= \supp{\overline{\wbf}^{\widehat{\lambda}}} \subset \supp{\overline{\wbf}^{\tilde{\lambda}}} \subset S^0\]
where $S^\star := \supp{\wstar}$ and $S^0:= \supp{\overline{\wbf}^0}$. Then
\[\CalL(\widehat{\lambda}) - \CalL(\tilde{\lambda}) \leq  \frac{\nrm{\overline{\Gbf}\left(\overline{\wbf}^0 -\overline{\wbf}^{\widehat{\lambda}}\right)}_2-\nrm{\overline{\Gbf}\left(\overline{\wbf}^0- \overline{\wbf}^{\tilde{\lambda}}\right)}_2}{\nrm{\overline{\bbf}}_2} - \frac{1}{\mathfrak{J}} \leq \frac{\nrm{\overline{\Gbf}\left(\overline{\wbf}^0 -\overline{\wbf}^{\widehat{\lambda}}\right)}_2}{\nrm{\overline{\bbf}}_2} - \frac{1}{\mathfrak{J}}.\]
Now, if we are in case $(i)$, it holds that $\overline{\wbf}^0=\overline{\wbf}^{\widehat{\lambda}}$ since no spurious terms are generated. For case $(ii)$, since $\overline{\Gbf}\overline{\wbf}^0 = \Gbf^\star\wstar$, we have
\[\nrm{\overline{\Gbf}\left(\overline{\wbf}^0 -\overline{\wbf}^{\widehat{\lambda}}\right)}_2= \nrm{\Gbf^\star\wstar -\overline{\Gbf}\overline{\wbf}^{\widehat{\lambda}}}_2 = \nrm{\Gbf^\star\wstar - \Gbf^\star\Abf\overline{\wbf}^{\widehat{\lambda}}}_2= \nrm{\Gbf^\star\wstar - \Gbf^\star(\Dbf +\Lbf)\overline{\wbf}^{\widehat{\lambda}}}_2.\]
Recalling that $\overline{\wbf}^{\widehat{\lambda}}_{S^\star}= \overline{\Gbf}_{S^\star}^\dagger\overline{\bbf}$, and using that $\overline{\bbf}=\bbf^\star=\Gbf^\star\wstar$, we see that replacing $\overline{\wbf}^{\widehat{\lambda}}$ in the previous line with any other vector supported on $S^\star$  will increase the norm. Hence, we may define a new vector $\wbf'$ block-wise (see \eqref{Ablock}) by
\[\wbf_B' = (\Dbf^{\pmb{\omega}})^{-1}\wstar_B\]
for blocks $B$ corresponding to trigonometric terms, and 
\[\wbf_B' = \wstar_B\]
for blocks $B$ corresponding to polynomial terms, which leads to
\[\nrm{\overline{\Gbf}\left(\overline{\wbf}^0 -\overline{\wbf}^{\widehat{\lambda}}\right)}_2 \leq \nrm{\Gbf^\star\wstar -\Gbf^\star(\Dbf +\Lbf)\wbf'}_2 = \nrm{\Gbf^\star_{p}\Lbf_p\wstar_p}_2\leq \nrm{\Gbf^\star_p}_2\nrm{\wstar}_2\nrm{\Lbf^{(p)}}_2\] 
where we use the subscript $p$ to denote columns pertaining to polynomial terms of degree at most $p$. Again using the bounds on $\nrm{\Lbf^{(p)}}$ in the appendix, we then get 
\[\nrm{\overline{\Gbf}\left(\overline{\wbf}^0 -\overline{\wbf}^{\widehat{\lambda}}\right)}_2 \leq \sigma^2 e{ p\choose 2}\nrm{\Gbf^\star_p}_2\nrm{\wstar}_2\]
where we also employed the necessary condition $\sigma^2<{p\choose 2}^{-1}$ from the previous lemma.
Hence, for all 
\[\sigma^2<\frac{\nrm{\Gbf^\star\wstar}_2}{\nrm{\Gbf^\star_p}_2\nrm{\wstar}_2} \frac{1}{\mathfrak{J}e{p \choose 2}}\]
it holds that $\CalL(\widehat{\lambda}) < \CalL(\tilde{\lambda})$. We conclude by combining the bounds on $\sigma_c$ given in the previous lemma.
\end{proof}

We now list some classical estimates that will be used to show that support recovery occurs with high probability as the resolution of the data increases. The bound \eqref{xnormbnd} is adapted from classical stability of full-rank linear least squares problems found in \cite[Theorem 5.1]{wedin1973perturbation}, together with elementary norm equivalences.

\begin{lemm}\label{least_squares_cont}
Let $\Abf \in \Rbb^{m\times n}$ have rank $n$ and $\ybf = \Abf\xbf$ for $\xbf\in \Rbb^n$, $\xbf\neq \mathbf{0}$. If a perturbed system $(\Abf',\ybf')$ satisfies \[\nrm{\Abf-\Abf'}_{\overrightarrow{\infty}}<\varepsilon, \qquad \nrm{\ybf-\ybf'}_\infty < \varepsilon\]
where  $\Abf'$ has rank $n$ and\footnote{Note that the constant $C(\Abf,\xbf)$ in general takes the form 
\[C(\Abf,\xbf) = \frac{\sqrt{m}+\sqrt{mn}\nrm{\xbf}_2}{\sigma_n(\Abf)(1-\alpha)}\]
valid for any $\alpha<1$ provided $\varepsilon\leq \alpha \frac{\sigma_n(\Abf)}{\sqrt{mn}}$. For convenience we chose $\alpha = 1/\sqrt{2}$ and the requirement $\varepsilon\leq \frac{\sigma_n(\Abf)}{\sqrt{2mn}}$.}
 $\varepsilon\leq \frac{\sigma_n(\Abf)}{\sqrt{2mn}}$, then the solution $\xbf' = (\Abf')^\dagger \ybf' = ((\Abf')^T\Abf')^{-1}(\Abf')^T\ybf'$ to the perturbed least squares problem satisfies
\begin{equation}\label{xnormbnd}
\nrm{\xbf-\xbf'}_\infty\leq \nrm{\xbf-\xbf'}_2 \leq \left(\frac{\sqrt{2m}+\sqrt{2mn}\nrm{\xbf}_2}{\sigma_n(\Abf)(\sqrt{2}-1)}\right)\varepsilon =: C(\Abf,\xbf)\varepsilon.
\end{equation}
\end{lemm}

We now introduce the main theorem of this section, showing subset support recovery with high probability so long as $\sigma<\sigma_c$ for some critical noise. Results are presented for Gaussian noise distributions, however similar results hold for more general noise distributions with suitable modifications to the proof using Lemma \ref{generalmomentmatrix}. 

\begin{thm}\label{conv_theorem_nosmooth}
Let Assumptions \ref{sec:regularity}-\ref{sec:conditioning} hold with $\rho$ a Gaussian mean-zero noise distribution. There exists a critical noise $\sigma_c>0$ and a stability tolerance $\tau$, both independent of $m$, such that for all $\sigma<\sigma_c$ and $t< \tau$, and for sufficiently large $m$, it holds that 
\begin{equation}\label{subsetsupprec_thm}
\supp{\widehat{\wbf}^{(m)}}\subset\supp{\wstar}
\end{equation}
with probability exceeding $1-2K\mathfrak{J}\exp\left(-\frac{c}{2}\left(mt\right)^{2/p_{\max}}\right)$, where $c$ is from Theorem \ref{matrix_concentration} and $\widehat{\wbf}^{(m)}$ is the one-shot MSTLS solution, $\widehat{\wbf}^{(m)} = \text{MSTLS}^{(1)}(\Gbf^{(m)},\bbf^{(m)})$. Moreover, if $\supp{\widehat{\wbf}^{(m)}}=\supp{\wstar}$, then it holds that
\begin{equation}
\nrm{\widehat{\wbf}^{(m)}-\wstar}_\infty < C'(t+\sigma^2)
\end{equation}
with the same probability, where $C'$ depends only on $(\Gbf^\star_{S^\star},\bbf^\star)$.
\end{thm}
\begin{proof}
Let $\overline{\wbf}^0 =\Abf^{-1}\wstar = \overline{\Gbf}^\dagger \overline{\bbf}$ and define $S = S^0\setminus S^\star$ where $S^0 = \supp{\overline{\wbf}^0}$. By Lemma \ref{lemm:bounds on sigma_c} there exists $\sigma_c'$ such that $\sigma<\sigma_c'$ ensures

\begin{equation}\label{delta1}
\delta_1 := \min_{j\in S^\star}\left\vert \overline{\wbf}^0_j\right\vert-\max_{j\in (S^\star)^c}\left\vert\overline{\wbf}^0_j\right\vert > 0.
\end{equation}
By Lemma \ref{lemm:STLS_supp_rec} this guarantees existence of $\widehat{\lambda}$ such that $\supp{\overline{\wbf}^{\widehat{\lambda}}} = S^\star$, where $\overline{\wbf}^{\widehat{\lambda}} = \text{STLS}^{(1)}(\overline{\Gbf},\overline{\bbf},\widehat{\lambda})$. Lemma \ref{lemm:subset_support_rec} then ensures that there exists $\sigma_c\leq \sigma_c'$ such that for $\sigma<\sigma_c$, 
\begin{equation}\label{delta2}
\delta_2 := \frac{\nrm{\overline{\bbf}}_2}{\mathfrak{J}}-\nrm{\overline{\Gbf}(\overline{\wbf}^0-\overline{\wbf}^{\widehat{\lambda}})}_2 > 0,
\end{equation}
which guarantees that $\widehat{\wbf} = \text{MSTLS}^{(1)}(\overline{\Gbf},\overline{\bbf})$ satisfies  $\supp{\widehat{\wbf}}\subset S^\star$. Next, using Theorem \ref{matrix_concentration} and Corollary \ref{matrix_corr}, we have that for all $t>0$ and sufficiently large $m$, it holds that
\[
\Pbb\left(\nrm{(\Gbf^{(m)},\bbf^{(m)}) - (\overline{\Gbf},\overline{\bbf})}_{\overrightarrow{\infty}} > t \right) \leq  2K\mathfrak{J}\exp\left(-\frac{c}{2}(mt)^{2/p_{\max}}\right).
\]
All that remains is to show that for sufficiently small $t$, $\nrm{(\Gbf^{(m)},\bbf^{(m)}) - (\overline{\Gbf},\overline{\bbf})}_{\overrightarrow{\infty}}<t$ leads to $\delta_1^{(m)},\delta_2^{(m)}>0$ where $\delta_1^{(m)},\delta_2^{(m)}$ are defined analogously to $\delta_1,\delta_2$ using $(\Gbf^{(m)},\bbf^{(m)})$. 

Indeed, assume that 
\begin{equation}\label{boundwitht}
\nrm{(\Gbf^{(m)},\bbf^{(m)}) - (\overline{\Gbf},\overline{\bbf})}_{\overrightarrow{\infty}}<t.
\end{equation}
If $t<\frac{\sigma_{\mathfrak{J}}(\overline{\Gbf})}{\sqrt{2K\mathfrak{J}}}$, where $\sigma_{\mathfrak{J}}(\overline{\Gbf})$ is the smallest singular value of $\overline{\Gbf}$, then by Lemma \ref{least_squares_cont}, we have that
\[\nrm{\overline{\wbf}^0-\wbf^{(m),0}}_\infty<Ct,\]
where $C=C(\overline{\Gbf},\overline{\wbf}^0)$ is defined in \eqref{xnormbnd}. This implies that
\[\delta_1^{(m)} := \min_{j\in S^\star}\left\vert \wbf^{(m),0}_j\right\vert-\max_{j\in (S^\star)^c}\left\vert \wbf^{(m),0}_j\right\vert \geq \delta_1 - 2Ct.\]
Hence, if $t < \min\left(\frac{\sigma_{\mathfrak{J}}(\overline{\Gbf})}{\sqrt{2K\mathfrak{J}}}, \frac{\delta_1}{2C}\right)$ and $m$ is sufficiently large, the probability that there exists $\widehat{\lambda}\in \pmb{\lambda}$ such that $\supp{\text{STLS}^{(1)}(\Gbf^{(m)},\bbf^{(m)},\widehat{\lambda})} = S^\star$ exceeds $1-2K\mathfrak{J}\exp\left(-\frac{c}{2}\left(mt\right)^{2/p_{\max}}\right)$. 
Next, we have
\begin{align*}
\delta_2^{(m)} &:= \frac{\nrm{\bbf^{(m)}}_2}{\mathfrak{J}}-\nrm{\Gbf^{(m)}(\wbf^{(m),0}-\wbf^{(m),\widehat{\lambda}})}_2 \\
&= \delta_2 + \frac{\nrm{\bbf^{(m)}}_2-\nrm{\overline{\bbf}}_2}{\mathfrak{J}}-\nrm{\Gbf^{(m)}(\wbf^{(m),0}-\wbf^{(m),\widehat{\lambda}})}_2 +\nrm{\overline{\Gbf}(\overline{\wbf}^0-\overline{\wbf}^{\widehat{\lambda}})}_2\\
&\geq \delta_2 - \frac{\nrm{\bbf^{(m)}-\overline{\bbf}}_2}{\mathfrak{J}} - \nrm{\Gbf^{(m)}-\overline{\Gbf}}_2\left(\nrm{\wbf^{(m),0}}_2+\nrm{\wbf^{(m),\widehat{\lambda}}}_2\right)\\
&\qquad-\nrm{\overline{\Gbf}}\left(\nrm{\wbf^{(m),0}-\overline{\wbf}^0}_2+\nrm{\wbf^{(m),\widehat{\lambda}}-\overline{\wbf}^{\widehat{\lambda}}}_2\right)\\
&\geq \delta_2 - \left(\frac{\sqrt{K}}{\mathfrak{J}} +\sqrt{K\mathfrak{J}}\left(\nrm{\overline{\wbf}^0}_2+\nrm{\overline{\wbf}^{\widehat{\lambda}}}_2+\delta_1\right)+2C\nrm{\overline{\Gbf}}_2\right)t\\
&=:\delta_2 - C't
\end{align*}
where we used the upper-bound $t<\frac{\delta_1}{2C}$ and the fact that 
\begin{equation}\label{subsetLSbound}
\nrm{\overline{\wbf}^\lambda-\wbf^{(m),\tilde{\lambda}}}_\infty<Ct
\end{equation}
for any $\lambda,\tilde{\lambda}$ that result in $\supp{\overline{\wbf}^\lambda} = \supp{\wbf^{(m),\tilde{\lambda}}}$, since \eqref{boundwitht} 
implies that $\nrm{\Gbf^{(m)}_{S'}-\overline{\Gbf}_{S'}}_{\overrightarrow{\infty}}<t$ on any subset $S'\subset\{1,\dots,\mathfrak{J}\}$. 
Hence, for sufficiently large $m$ and any 
\[ t < \min\left\{\frac{\sigma_{\mathfrak{J}}(\overline{\Gbf})}{\sqrt{2K\mathfrak{J}}},\ \frac{\delta_1}{2C},\ \frac{\delta_2}{C'}\right\}:=\tau,\]
we get that $\supp{\text{MSTLS}^{(1)}(\Gbf^{(m)},\bbf^{(m)})} \subset S^\star$ with probability exceeding $1-2K\mathfrak{J}\exp\left(-\frac{c}{2}\left(mt\right)^{2/p_{\max}}\right)$.

To prove the coefficient accuracy, if $\supp{\what^{(m)}} = S^\star$ (i.e.\ $\what^{(m)}=\wbf^{(m),\widehat{\lambda}}$), \eqref{subsetLSbound} implies that 
\[\nrm{\what^{(m)} - (\overline{\Gbf}_{S^\star})^\dagger\overline{\bbf}}_\infty <Ct\]
where $C=C(\overline{\Gbf}_S^\star, (\overline{\Gbf}_S^\star)^\dagger\overline{\bbf})$ is again defined in \eqref{xnormbnd}. Furthermore, properties of $\Abf$ imply that for $C''$ depending on the maximum polynomial degree $p$ and maximum frequency $\omega_{\max}$ appearing in the true model we have $\nrm{\overline{\Gbf}_{S^\star}-\Gbf_{S^\star}^\star}_\infty<C''\sigma^2$, which together with Lemma \ref{proj_cont} implies that
\[\nrm{(\overline{\Gbf}_{S^\star}^\dagger - \Gbf_{S^\star}^\dagger)\bbf^\star}_\infty<\tilde{C}\sigma^2.\]
Altogether, if $\supp{\what^{(m)}} = S^\star$, then with probability exceeding $1-2K\mathfrak{J}\exp\left(-\frac{c}{2}\left(mt\right)^{2/p_{\max}}\right)$ it holds that 
\[\nrm{\what^{(m)} - \wstar}_\infty\leq \nrm{\what^{(m)} - (\overline{\Gbf}_{S^\star})^\dagger\overline{\bbf}}_\infty+\nrm{(\overline{\Gbf}_{S^\star})^\dagger\overline{\bbf} - \wstar}_\infty\leq Ct + \tilde{C}\sigma^2 \leq C'(t+\sigma^2).\]
\end{proof}

\begin{rmrk}\label{promise_of_fullsupprec}
As alluded to at the beginning of this section, if we make an additional assumption on the noise-free continuous data, we can strengthen Lemma \ref{lemm:subset_support_rec} to ensure $\supp{\widehat{\wbf}}=\supp{\wstar}$ for all sufficiently small $\sigma$, which subsequently strengthens Theorem \ref{conv_theorem_nosmooth} to ensure $\supp{\widehat{\wbf}^{(m)}}=\supp{\wstar}$ with high probability for all sufficiently small $\sigma$ and sufficiently large $m$. This condition is the following: 
\begin{equation}\label{condforsupprec}
\mu^\star:= \min_{S\subsetneq S^\star}\frac{\nrm{\Pbf^\perp_{\Gbf^\star_{S^\star\setminus S}} \bbf^\star}}{\nrm{\bbf^\star}}-\frac{|S|+1}{\mathfrak{J}}>0.
\end{equation}
Using that $\Pbf^\perp_{\Gbf^\star_{S^\star\setminus S}} \bbf^\star = \Pbf^\perp_{\Gbf^\star_{S^\star\setminus S}}\Gbf^\star_S\wstar_S$, in words, this says that the contribution of each subset of true terms $\Gbf_S^\star\wstar_S$ that is orthogonal to the subspace spanned by the remaining terms ($\text{span}(\Gbf^\star_{S^\star\setminus S})$) cannot be arbitrarily small. Specifically, this orthogonal contribution must be at least $(\frac{|S|+1}{\mathfrak{J}})\nrm{\bbf^\star}$. In Appendix \ref{app:supprec} we show that this is sufficient to guarantee that the MSTLS loss $\CalL$ satisfies $\CalL(\widehat{\lambda})<\CalL(\tilde{\lambda})$ for all $\supp{\overline{\wbf}^{\widehat{\lambda}}} = S^\star$ and $\supp{\overline{\wbf}^{\tilde{\lambda}}}=\tilde{S}\subsetneq S^\star$. However, condition \eqref{condforsupprec} is unsatisfactory because it cannot be checked without knowledge of the true support $S^\star$. Nevertheless, with $|S|=1$, we can interpret \eqref{condforsupprec} as a modeling criteria: each true term in the model must provide a unique (orthogonal) contribution to the dynamics given in $\bbf^\star$ of at least $(2/\mathfrak{J})100\%$. 
\end{rmrk}

\subsection{Asymptotic consistency with filtering}\label{smooth}

In this section we show that filtering the data prior to building the linear system $(\Gbf^{(m)},\bbf^{(m)})$ results in unconditional convergence of estimators to the true weights, in other words $\what^{(m)}\to \wstar$ in probability, provided the filtering window scales appropriately.

We define the filtered data $\tilde{\Ubf}^{(m)}$ with respect to a discrete convolutional filter $\pmb{\nu}^{(m)}$ as
\[\tilde{\Ubf}^{(m)} = \pmb{\nu}^{(m)}*\Ubf^{(m)},\]
where $\pmb{\nu}^{(m)}\in \Rbb^{n_1^\nu\times\cdots\times n_{d+1}^\nu}$ satisfies $\nrm{\pmb{\nu}^{(m)}}_{\overrightarrow{1}}=1$ and $\pmb{\nu}^{(m)} > 0$ (i.e\ $\pmb{\nu}^{(m)}$ is a discrete probability distribution). For convenience we perform symmetric reflection of the data at the boundaries to maintain the same number of data points in $\tilde{\Ubf}^{(m)}$ and $\Ubf^{(m)}$. The filter $\pmb{\nu}^{(m)}$ is characterized by its {\it filter width} in total number of gridpoints at level $m$: 
\[m^{(\nu)}:= \prod_{i=1}^{d+1} n_i^\nu.\]

The class of $\pmb{\nu}^{(m)}$ resulting in convergence is large, however for simplicity we restrict our attention to the simple moving average filter, in which case each entry of $\pmb{\nu}^{(m)}$ is equal to $1/m^{(\nu)}$. In particular, we show that the moving average filter provides convergence so long as $m^{(\nu)}\gtrsim m^{\alpha}$ as $m\to \infty$ for some $\alpha<1$, and as before, $m$ is the number of points that the test function $\psi$ is support on at the grid resolution $(\Delta x^{(m)},\Delta t^{(m)})$. Since the simple moving average filter reduces variance by a factor of $1/m^{(\nu)}=\nrm{\pmb{\nu}^{(m)}}_{\overrightarrow{2}}^2$, we conjecture that the results below hold for other filters satisfying $\nrm{\pmb{\nu}^{(m)}}_{\overrightarrow{2}}^2\leq m^{-\alpha}$.

We denote the WSINDy linear system built using the filtered data $\tilde{\Ubf}^{(m)}$ as
\[(\tilde{\Gbf}^{(m)}, \tilde{\bbf}^{(m)}).\]
Our approach is to show that $(\tilde{\Gbf}^{(m)}, \tilde{\bbf}^{(m)})\to (\Gbf^\star, \bbf^\star)$ in probability, after which the full-rank assumption on $\Gbf^\star$ implies that the least-squares solution $\tilde{\wbf}^0 = (\tilde{\Gbf}^{(m)})^\dagger \tilde{\bbf}^{(m)}$ converges to $\wstar$. The main hurdle in proving this is the correlations that arise from filtering the data, since entries of $(\tilde{\Gbf}^{(m)}, \tilde{\bbf}^{(m)})$ are now sums of {\it correlated} heavy-tailed random variables, which falls outside of the analysis of works  \cite{vladimirova2020sub,bakhshizadeh2020sharp}. The proof can be found in Appendix \ref{app:conclemmas}.

\begin{thm}\label{smoothconc}
Let assumptions \ref{sec:regularity}, \ref{sec:sampling}, and \ref{sec:testfcn}
hold and let $\pmb{\nu}^{(m)}$ be a simple moving average filter with $m^{(\nu)}\gtrsim m^{\alpha}$ for some $\alpha\in(0,1)$. Further, let $\CalF$ be a library of locally Lipschitz functions such that $\forall f\in \CalF$ and $(x,y)\in \Rbb$,
\[|f(x)-f(y)|\leq C_\CalF|x-y|\left(1+|x-y|^{p_{\max}-1}\right)\]
for some $C_\CalF\geq 0$ and $p_{\max}\geq 1$. Then for $t> \mathfrak{t}(m)$, where $\mathfrak{t}(m) = \CalO(m^{-\alpha \left(\frac{k}{d+1}-\frac{1}{2}\right)})$ and $k$ is the Sobolev regularity of $u\in \CalH^{k,\infty}(D)$, and sufficiently large $m$, we have 
\[\Pbb \left(\|(\tilde{\Gbf}^{(m)},\tilde{\Gbf}^{(m)}) - (\Gbf^\star,\bbf^\star)\|_\infty>t\right) \leq 3K\mathfrak{J}\exp\left(-\frac{c}{2^{1+p_{\max}}}[m(t-\mathfrak{t}(m))]^{2/{p_{\max}}}\right).\]
where $c$ is the same rate from Theorem \ref{matrix_concentration}.
\end{thm}

\begin{rmrk}\label{increase_rate_rmrk}
A more careful analysis will lead to the rate $c$ replaced by a variable rate $c^*(m)$ (bounded below by $c$ from Theorem \ref{matrix_concentration}) that increases as $\sim m^{\alpha}$, since filtering transforms the sub-Gaussian noise distribution at each $m$ into a new distribution $\rho^{(m)}$ that has a reduced variance $\sigma^2/m^\alpha$. 
\end{rmrk}


\begin{rmrk}
As well as factors which serve to increase the rate $c$ in Theorem \ref{smoothconc}, convergence is largely dictated by how rapidly $\mathfrak{t}(m)$ decreases. This is determined by three deterministic effects: (1) the bias $|\tilde{\Ubf}^{(m),\star}-\Ubf^{(m),\star}|$ between the filtered and unfiltered noise-free data, (2) the distributional convergence of $\rho^{(m)}$ to $\delta_0$ as $m\to \infty$, and (3) the convergence of the trapezoidal rule on functions $\partial \psi^{\alpha^s}(\cdot) f(u(\cdot))$. In the case that $u$ is locally polynomial of degree $q+1$, higher-order filters can be used to reduce (1) to $\CalO(m^{-\frac{(q+1)(1-\alpha)}{d+1}})$. However, in most cases (2) still dominates, limiting $\mathfrak{t}(m)$ to $\CalO(m^{-\alpha/2})$. If $\CalF$ contains only analytic functions, then (2) acquires the rate $\CalO(m^{-\alpha})$. Higher-order filters are of no help in decreasing (2) since $|f*\rho^{(m)}(u)-f*\delta_0(u)|=\CalO(\nrm{\pmb{\nu}^{(m)}}_2^2)$, and $\nrm{\pmb{\nu}^{(m)}}_2^2$ is minimized at the simple moving average filter. For (3), if $u\in H^k$ with $k>(d+1)/2$ and $f_j$ are smooth, then the convergence rate of the trapezoidal rule is $\CalO(m^{-k+(d+1)/2})$, which leads to rapidly decreasing $\mathfrak{t}(m)$ with higher $k$, and thus faster entry into the region of exponential concentration. Ultimately, while filtering the data enables concentration to the noise-free problem, it is clear that more work is needed to make filtering useful in practice due to these intermediate biases.  
\end{rmrk}

Lastly, Lemma \ref{least_squares_cont} and Theorem \ref{smoothconc} directly imply the following coefficient accuracy.

\begin{corr}
Under the same conditions as Theorem \ref{smoothconc}, it holds that 
\[\nrm{\tilde{\wbf}^{(m),0}-\wstar}_\infty\leq Ct\]
 with probability exceeding $1 - 3K\mathfrak{J}\exp\left(-\frac{c}{2^{1+p_{\max}}}[m(t-\mathfrak{t}(m))]^{2/{p_{\max}}}\right)$,
where $C=C(\Gbf^\star,\wstar)$ is defined in Lemma \ref{least_squares_cont} and $\tilde{\wbf}^{(m),0} = (\tilde{\Gbf}^{(m)})^\dagger\tilde{\bbf}$ is the filtered least squares solution. 
\end{corr}

\begin{rmrk}
Since it is not advised to simply use the least-squares solution as in the previous corollary (which will in general not be sparse), we note that Theorem \ref{smoothconc} directly implies that $\supp{\what^{(m)}}\subset \supp{\wstar}$ with high probability as in Theorem \ref{conv_theorem_nosmooth}, where here $\what^{(m)} = \text{MSTLS}^{(1)}(\tilde{\Gbf}^{(m)},\tilde{\bbf}^{(m)})$ is the one-shot MSTLS solution on the filtered linear system $(\tilde{\Gbf}^{(m)},\tilde{\bbf}^{(m)})$. Moreover, in the case of filtered data the result is not restricted to trigonometric and polynomial libraries, but holds for {\it any locally Lipschitz library} $\CalF$. Finally, condition \eqref{condforsupprec} implies full support recovery $\supp{\what^{(m)}} = \supp{\wstar}$ with high probability from filtered data, using similar arguments as in Appendix \ref{app:supprec}.
\end{rmrk}

\section{Numerical Experiments}\label{exp}

We now test the theoretical results in the previous sections using four example problems: (1) the Lorenz system, (2) a hyper-diffusive Kuramoto-Sivashinsky (KS)-type equation, (3) a cubic oscillator, and (4) a nonlinear viscous Burgers-type model. Examples (1) and (2) are ordinary and partial differential equations, respectively, which {do not exhibit a critical noise (equivalently, $\sigma_c=\infty$). For these examples we demonstrate that recovery of the correct model occurs with probability approaching 1 as $m\to \infty$ across the noise spectrum. In contrast, systems (3) and (4) do exhibit a finite critical noise $\sigma_c$, and in these cases we show that a simple moving average filter is sufficient to enable convergence for $\sigma\geq \sigma_c$.

\begin{rmrk}
For cases of unconditional convergence (examples (1) and (2)) we do not report results for recovery from filtered data, as these were consistently worse than their unfiltered counterparts. This is supported by Theorems \ref{conv_theorem_nosmooth} and \ref{smoothconc}, as convergence for the filtered linear system is slower than for raw data. 
\end{rmrk}

\subsection{Data generation and WSINDy settings}

For each example we subsample a high-accuracy fine-grid simulation of the system in order to test the performance at different test function support size $m$. We then add mean-zero Gaussian white noise to every datapoint with specified variance $\sigma^2$. Where specified, we use a noise ratio $\sigma_{NR}$ to determine $\sigma$, and set
\[\sigma = \nrm{\Ubf^{(m),\star}}_{stdev'}\sigma_{NR}\]
where $\nrm{\Ubf^{(m),\star}}_{stdev'}$ indicates the standard deviation of the clean data $\Ubf^{(m),\star}$ stretched into a column vector.

Throughout we use the MSTLS algorithm outlined in equations \eqref{STLS}-\eqref{lossfcn} with sparsity thresholds $\log_{10}\pmb{\lambda} = \texttt{linspace}(-4,0,100)$. For simplicity, we fix the reference test function along each coordinate to be the $C^\infty_c(\Rbb)$ bump function
\[\phi(v) = \ind{(-1,1)}(v)\exp\left(\frac{9}{v^2-1}\right).\]
The model library and convolution query points vary by example as described below.

\subsection{Performance metrics}

We are mostly concerned with verifying the different asymptotic consistency results above for raw and filtered data over a range of noise levels. For each noise level $\sigma$, we run the WSINDy algorithm on 200 different independent noise instantiations at each test function support size $m$, producing 200 learned models $(\what^{(m),(i)})_{i=1,\dots,200}$ which we average results over. We record the probability of support recovery
\begin{equation}
\Pbb\left(S^{(m)}=S^\star\right)\approx \frac{1}{200}\sum_{i=1}^{200} \mathbbm{1}\left( \supp{\what^{(m),(i)}}=\supp{\wstar} \right)
\end{equation}
and where relevant the probability of support inclusion
\begin{equation}
\Pbb\left(S^{(m)}\subset S^\star\right)\approx \frac{1}{200}\sum_{i=1}^{200} \mathbbm{1}\left( \supp{\what^{(m),(i)}}\subset\supp{\wstar} \right).
\end{equation}
We also report the maximum relative coefficient error over true support set, defined as
\begin{equation}\label{eq:err_inf}
E_\infty(m) = \Ebb\left[\max_{j\in S^\star} \frac{|\what^{(m)}_j-\wstar_j|}{|\wstar_j|}\right]\approx  \frac{1}{200}\sum_{i=1}^{200}\max_{j\in S^\star} \frac{|\what^{(m),(i)}_j-\wstar_j|}{|\wstar_j|} .
\end{equation}

\subsection{Unconditional consistency} 

\begin{figure}[H]
\begin{center}
\begin{tabular}{cc}
\includegraphics[trim={40 60 60 15},clip,width=0.45\textwidth]{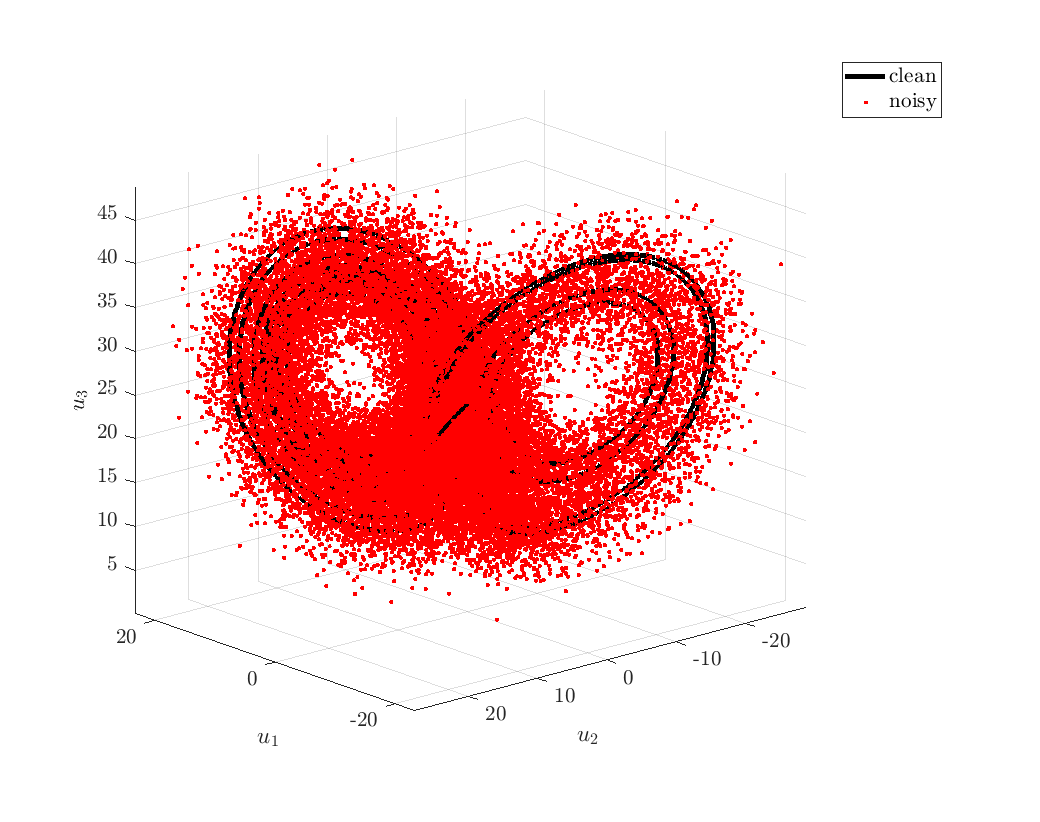} &
\includegraphics[trim={15 3 40 10},clip,width=0.45\textwidth]{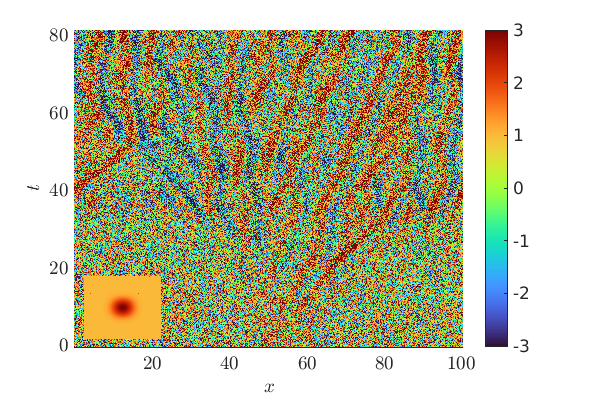} 
\end{tabular}
\end{center}
\caption{Example datasets for the Lorenz system (left) at 20\% noise and hyper-KS system (right) at $100\%$ noise.}\label{lorenzKSsol}
\end{figure}

\subsubsection{Lorenz System}
The true model equations for the Lorenz system are
\begin{equation}\label{eq:lorenz}
\begin{dcases} 
\frac{du_1}{dt} = (-10)u_1+(10)u_2 \\
\frac{du_2}{dt} = (28)u_1+(-1)u_2+(-1)u_1u_3 \\
\frac{du_3}{dt} = (-8/3)u_3+(1)u_1u_2.
\end{dcases}
\end{equation}
Since each term is either linear or bilinear, the associated continuum linear system is unbiased, in other words $\overline{\Gbf}_{S^\star} = \Gbf^\star_{S^\star}$, hence the system does not exhibit a critical noise. We simulate \eqref{eq:lorenz} for $t\in [0,10]$ using RK-45 with absolute tolerance $10^{-12}$ using $250,000$ equally-spaced points. This fine-grid solution is then subsampled by successive factors of two, leading to coarser data with approximate total numbers of points $\{2^{-k}250,000\ :\ k=0,\dots,9\}$. At each level of resolution we specify the test function width $m$ to be $2\%$ of the total timeseries, so that $|\supp{\psi}|/T=0.02$ for all $m$. We use a maximum of $K=1000$ equally-spaced convolution query points when constructing $\Gbf^{(m)}$. (When the number of timepoints $M$ is less than 1000 we use all possible query points, or $K = M-m+1$). We let $\CalF$ be the set of polynomials up to total degree $p_{\max}=6$ in the state variables $(u_1,u_2,u_3)$, leading to 84 library terms.

We examine noise ratios $\sigma_{NR}$ from $0.001$ to $1$, which translate to $\sigma$ from $0.013$ to $13$, due to $\nrm{\Ubf^{(m),\star}}_{stdev'}$ $\approx 13$. Figure \ref{lorenzKS} (left column) shows that across the noise spectrum we get asymptotic recovery of the correct system, with coefficient errors $E_\infty$ eventually entering a Monte-Carlo-type $\CalO(m^{-1/2})$ convergence.

\subsubsection{Hyper-KS}

In order to demonstrate the ability of WSINDy to recover high-order PDEs, we examine a hyper-diffusive, dispersive evolution equation that exhibits spatiotemporal chaos similar to the KS equation:
\begin{equation}\label{hyperKS}
\partial_t u = (1)\partial_{xxxx} u + (0.75)\partial_{xxxxxx} u + (-0.5)\partial_x(u^2) + (0.1)\partial_{xxx} u^2.
\end{equation}
Such models are of general interest for their potential to model challenging dynamics such as earthquakes and flame propagation. The dynamics of \eqref{hyperKS} are elaborated on in Appendix \ref{app:hyperKS}. Despite its complexity at face-value, the model \eqref{hyperKS} can be recovered using WSINDy as $m\to \infty$ with no restrictions on the noise level. This is due to the fact that the system is composed of only linear and quadratic terms, with no quadratic growth term, so falling into case $(i)$ for Lemmas \ref{lemm:bounds on sigma_c} and \ref{lemm:subset_support_rec}.

We simulate \eqref{hyperKS} using a ETDRK4 and Fourier-spectral collocation on a fine grid $(\Xbf^f,\tbf^f)\subset [0,32\pi]\times [0,82]$ containing $1024\times 1025$ points in space and time to mimic the continuum limit. We then examine a range of resolutions, with the coarsest grid have $64\times 65$ points, or $32$ times coarser than $(\Xbf^f,\tbf^f)$. We examine $\sigma_{NR}$ in the range $0.001$ to $1$, where $\sigma_{NR}\approx \sigma$ here since $\nrm{\Ubf^{(m),\star}}_{stdev'}\approx 1$. The reference test function is set so that $|\supp{\psi}|/|\Omega\times(0,T)| = 1/25$ (see inset plot of Figure \ref{lorenzKSsol} (right)). We fix the library at $\CalF = (u^q)_{q=0,\dots,8}$ and the differential operators $\partial^{\pmb{\alpha}} = (\partial_x^q)_{q=0,\dots,8}$, leading to a total of 73 library terms.

Similar to the Lorenz system, we observe support recovery of \eqref{hyperKS} in Figure \ref{lorenzKS} (right column) with probability approaching one for all $\sigma_{NR}$ examined. As well, the error scales asymptotically like $\CalO(m^{-1/2})$, where we note that $m=m_xm_t$. It is interesting to note that for lower resolutions, we still recover the correct system at fairly high probability. At $\sigma_{NR} = 0.01$, we recover the correct system with $>98\%$ probability from a grid of size $86\times 86$ (corresponding to $(m_x,m_t) = (15,15)$), and with average coefficient error $E_\infty<0.026$. This is surprising due to the dependence on $\partial^4_xu$ and $\partial^6_xu$. 

\begin{figure}
\begin{center}
\begin{tabular}{cc}
\includegraphics[trim={5 0 40 20},clip,width=0.49\textwidth]{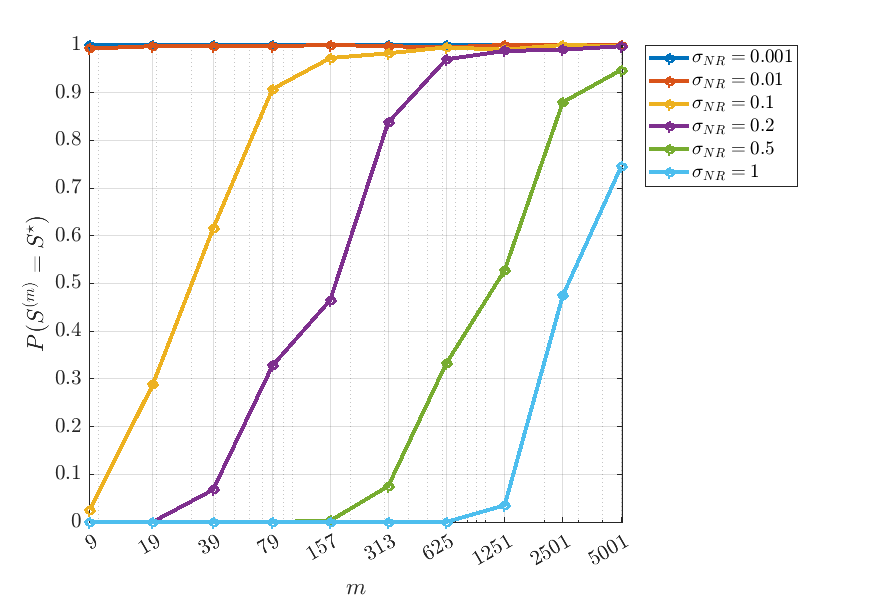} &
\includegraphics[trim={5 0 40 20},clip,width=0.49\textwidth]{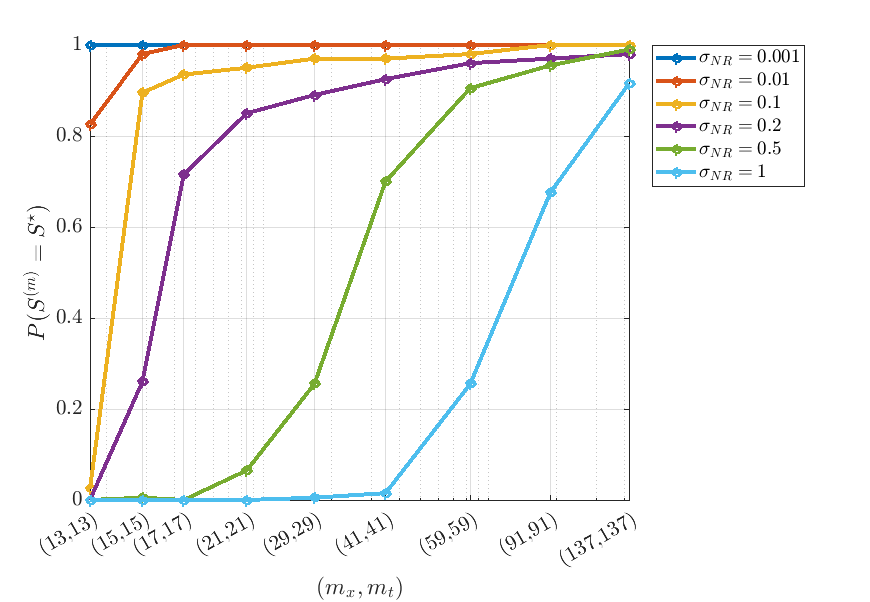} \\
\includegraphics[trim={5 0 40 20},clip,width=0.49\textwidth]{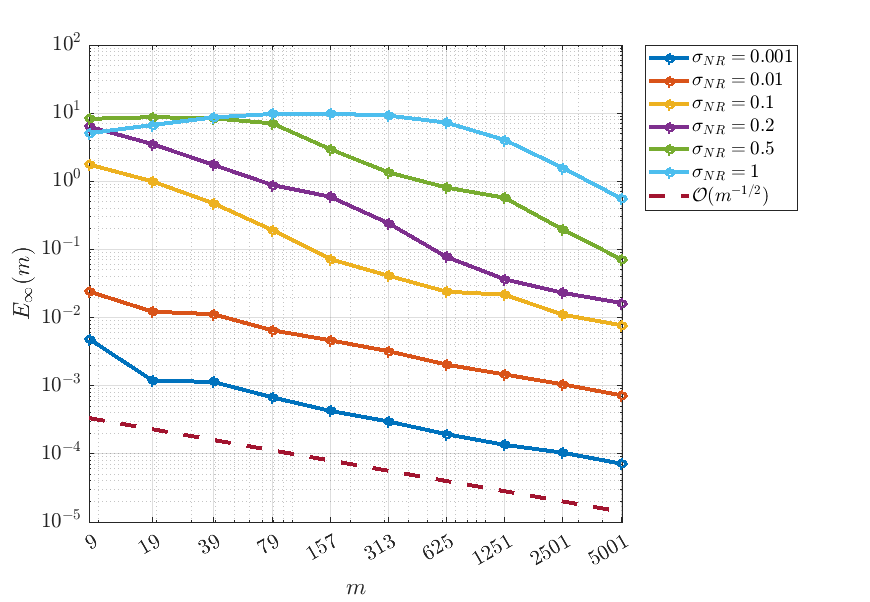} &
\includegraphics[trim={5 0 40 20},clip,width=0.49\textwidth]{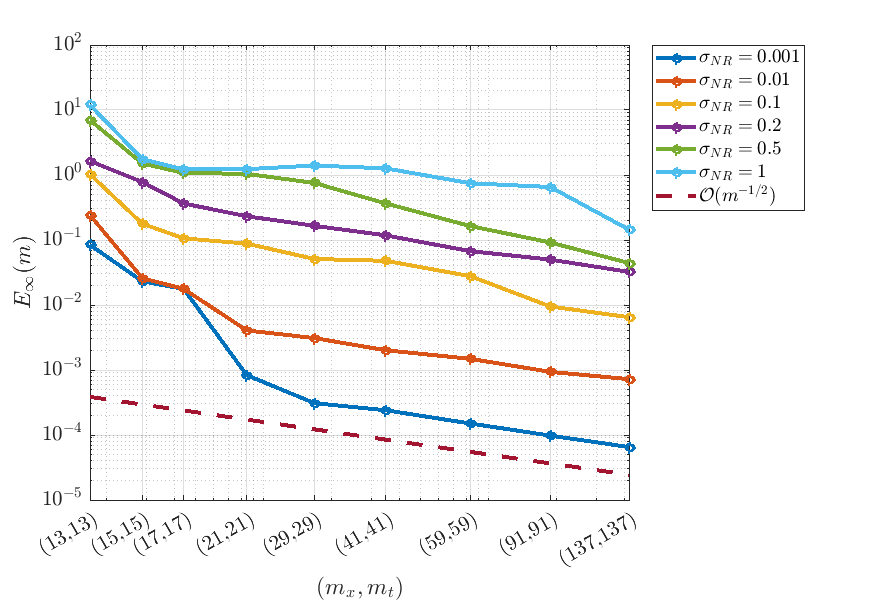} 
\end{tabular}
\end{center}
\caption{Recovery results for the Lorenz system (left) and hyper-KS equation (right). Both systems exhibit similar asymptotic consistency trends, where support recovery is achievable at any noise level if $m$ is taken large enough, and errors eventually decrease at a rate $\CalO(m^{-1/2})$ (note in the right column $m=m_xm_t$).}
\label{lorenzKS}
\end{figure}

\subsection{Conditional consistency and filtering}

\begin{figure}
\begin{center}
\begin{tabular}{cc}
\includegraphics[trim={15 3 45 10},clip,width=0.45\textwidth]{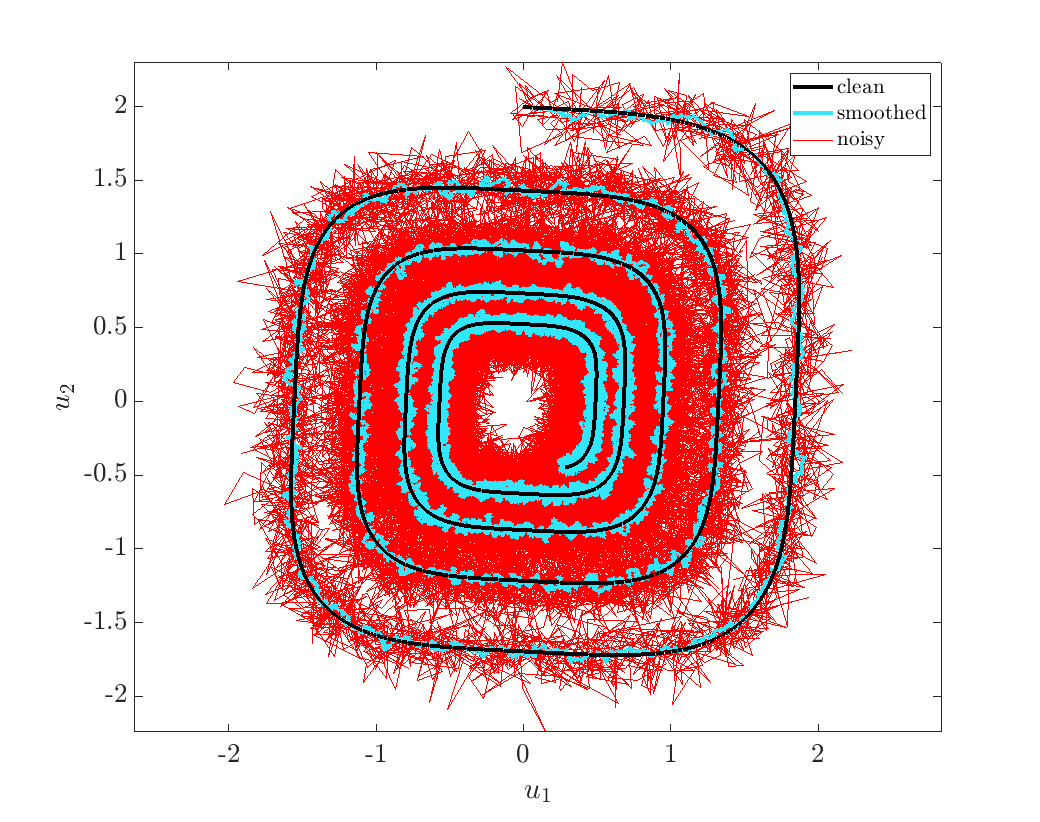} &
\includegraphics[trim={15 3 45 10},clip,width=0.5\textwidth]{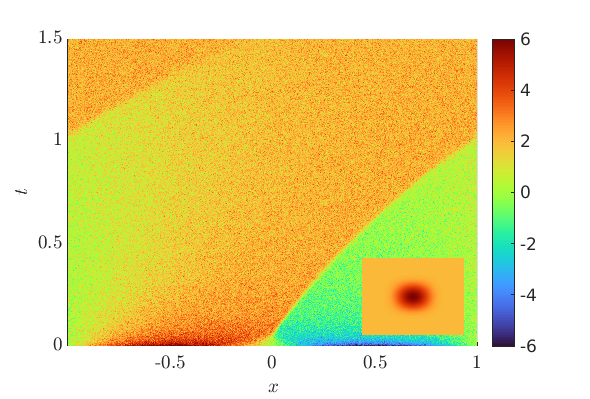}
\end{tabular}
\end{center}
\caption{Left: numerical solution to the cubic oscillator \eqref{cubicosc} together with noisy data at $\sigma=\sigma_c$ and corresponding filtered data. Right: numerical solution to the nonlinear viscous Burgers model \eqref{burgers} with noise $\sigma = 10\sigma_c$. }
\label{sol}
\end{figure}

The next two examples do exhibit a finite critical noise, and so we report the performance of WSINDy on the raw data as well as filtered data using a simple moving average filter. We design the filter width $m^{(\nu)}$ using the bounds \eqref{explicitboundsonsigmac}
in Lemma \ref{lemm:bounds on sigma_c}. First we use the data to calculate an estimate of $\sigma$, denoted by $\sigma_{\text{est}}$, using the method in Appendix \ref{app:sigma_est}. We then specify a prior $\tau^\star=0.01$ for the ratio $\min_{i\in S^\star}|\wstar_i|/\nrm{\wstar}_\infty$. Assuming $\tau^\star$ is well-specified, a desirable quality of the filter given the bounds \eqref{explicitboundsonsigmac} is to have 
\[\frac{1}{\tau^\star}{p_{\max} \choose 2}\sigma_\nu^2 < 1\]
where $\sigma_\nu^2 = \frac{\sigma^2}{m^{(\nu)}}$ is the reduced variance of the filtered data, given an initial noise variance of $\sigma^2$. We can then solve for $m^{(\nu)}$ using our variance estimate and our prior $\tau^\star$:
\[m^{(\nu)} > {p_{\max} \choose 2}\frac{\sigma_\text{est}^2}{\tau^\star}.\]
In both examples we use polynomials up to total degree $p_{\max} = 6$, and we set $\tau^\star = 0.01$, which translates into a filter width of 
\begin{equation}
m^{(\nu)}_1 > \left\lfloor (1500  \sigma_\text{est}^2)^{\frac{1}{d+1}} \right\rfloor
\end{equation}
points in each dimension. We use this together with a restriction on the filter width depending on the test function support size $m$, and set 
\[m^{(\nu)}_1 = \left\lfloor \min\left(2(1500  \sigma_\text{est}^2)^{\frac{1}{d+1}}, m^{\frac{1}{d+1}}/2\right)\right\rfloor.\]
In this way the variance is reduced enough (given assumptions on $\tau^\star$) to cancel spurious terms that arise, but the filter only covers at most $1/2^{d+1}$ of $|\supp{\psi}|$ to limit the resulting bias. We make no claim that this is the optimal way to choose the filter, but as we'll see below, it is sufficient to ensure support recovery beyond the critical noise level.

\begin{figure}
\begin{center}
\begin{tabular}{cc}
\includegraphics[trim={5 0 40 20},clip,width=0.46\textwidth]{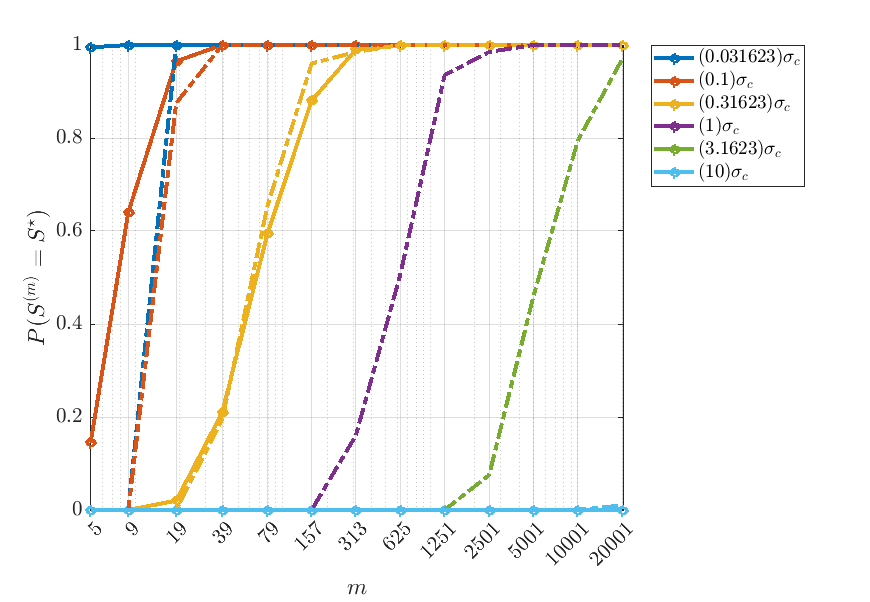} &
\includegraphics[trim={5 0 40 20},clip,width=0.46\textwidth]{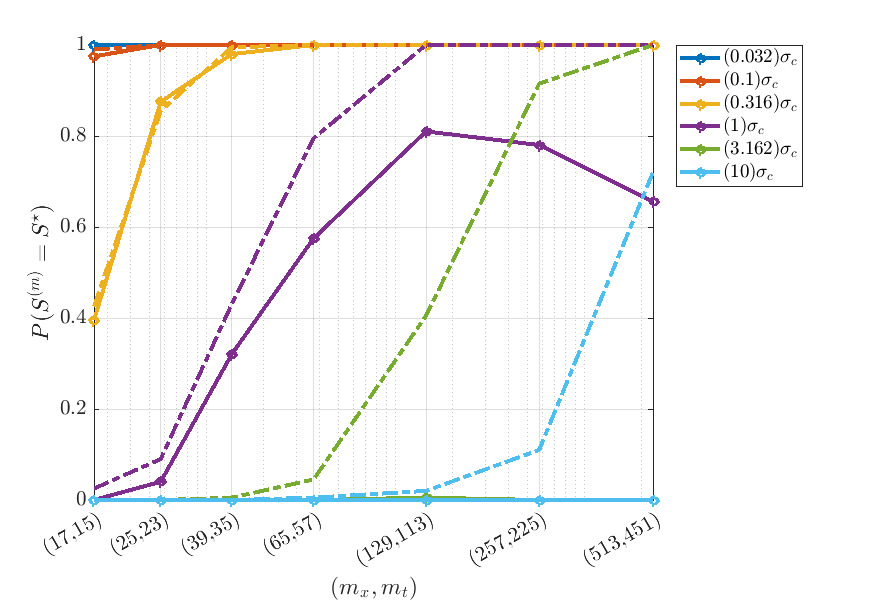} \\
\includegraphics[trim={5 0 40 20},clip,width=0.46\textwidth]{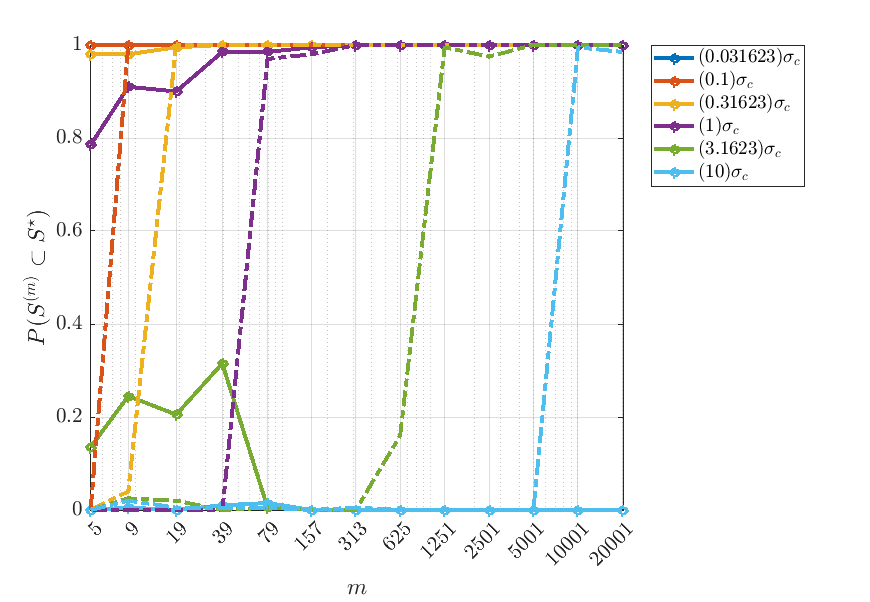} &
\includegraphics[trim={5 0 40 20},clip,width=0.46\textwidth]{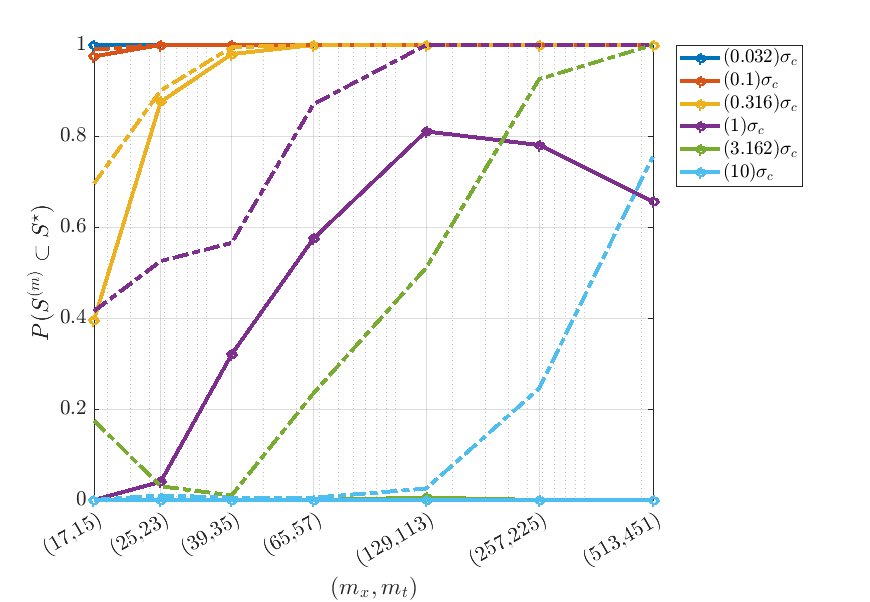} \\
\includegraphics[trim={5 0 40 20},clip,width=0.46\textwidth]{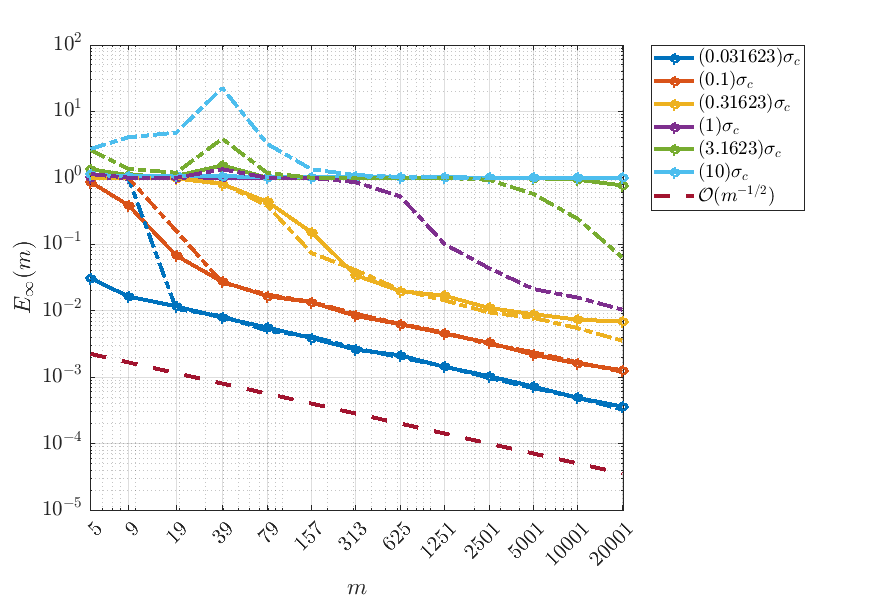} &
\includegraphics[trim={5 0 40 20},clip,width=0.46\textwidth]{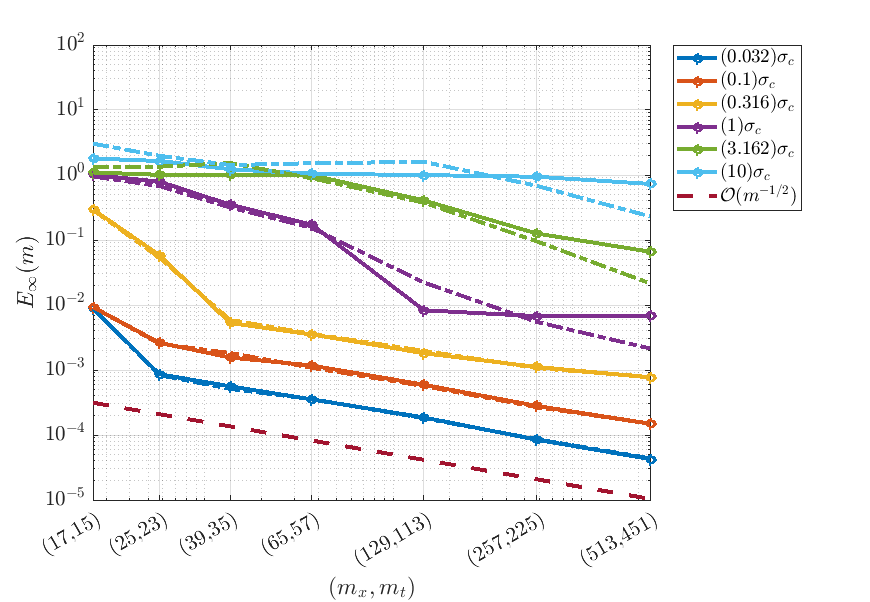} 
\end{tabular}
\end{center}
\caption{Recovery results for the Cubic Oscillator \eqref{cubicosc} (left) and nonlinear viscous Burgers equation \eqref{burgers} (right). Solid lines correspond to raw (unfiltered) data and dot-dashed lines correspond to filtered data. For both systems WSINDy applied to the raw data only recovers the correct model for $\sigma<\sigma_c$, whereas support recovery is achievable at higher noise levels from the filtered data if $m$ is taken large enough. Moreover, solutions converge on the filtered data at a higher rate than $\CalO(m^{-1/2})$ (see dot-dashed lines in the bottom right plot). (Note that in the right column $m=m_xm_t$).}
\label{cubicNLVB}
\end{figure}

\subsubsection{Cubic Oscillator}

\begin{equation}\label{cubicosc}
\begin{dcases} 
\frac{du_1}{dt} = (2)u_1^3+(-0.1)u_2^3 \\
\frac{du_2}{dt} = (-0.1)u_1^3+(-2)u_2^3
\end{dcases}
\end{equation}

The cubic oscillator has origins in quantum mechanics and represents a particle in an anharmonic potential well. The cubic terms each generate spurious linear terms in the limit, which persist in the learned equations for noise levels above the critical noise $\sigma_c = \sqrt{0.1/6}\approx 0.13$. A fine-grid solution is obtained on $10^6$ equally-spaced timepoints $\tbf \subset [0,25]$. The data is then subsampled to obtain coarser resolutions. As with the Lorenz system, we fix $|\supp{\psi}|=T/50$ and we use a maximum of $K=1000$ equally-spaced query points (meaning, as before, that $K=1000$ for total number of timepoints $M>1000$, and otherwise $K=M-m+1$). We examine noise levels $\sigma$ in the range $10^{-2}\sigma_c\approx 0.0013$ to $10\sigma_c \approx 1.3$. In terms of the noise ratio $\sigma_{NR}$, given that $\nrm{\Ubf^{(m),\star}}_{stdev'}\approx 0.6253$, this corresponds to $\sigma_{NR}$ in the range $0.0021$ to $2.064$.

In Figure \ref{cubicNLVB} (left) we report the asymptotic recovery trends for \eqref{cubicosc}. For noise levels $\sigma < \sigma_c$, we recover the correct support as $m\to\infty$ for both raw and filtered datasets, with WSINDy performing better overall on the raw data. For $\sigma\geq \sigma_c$, the raw data fails to lead to support recovery, but the method works successfully on the filtered data, as predicted. However, it is clear that this particular problem exhibits slow convergence. We observe rapid subset support recovery (left middle panel of Figure \ref{cubicNLVB}) using filtered data, corresponding to dropping one or both of the smaller terms in \eqref{cubicosc}, but full support recovery requires impractical amounts of data using the chosen hyperparameters. Nevertheless the asymptotic theoretical results are captured. 

\subsubsection{Nonlinear Viscous Burgers}

In our last example we examine the model
\begin{equation}\label{burgers}
\partial_t u = (0.01)\partial_{xx} u + (-0.5)\partial_x(u^2) + (-1)u^3+(2)u^2+(1)u^0,
\end{equation}
which allows for explicit computation of the critical noise level $\sigma_c$ similar to the cubic oscillator. Specifically, the $(-1)u^3$ term generates a spurious term $(-3\sigma^2) u^1$ in the continuum limit, and the term $(2)u^2$ generates a term $(2\sigma^2)u^0$, which corrupts the existing coefficient of $u^0$. With a diffusion coefficient of $0.01$, the true model will not be identified for noise levels above the critical noise $\sigma_c:=\sqrt{0.01/3}\approx 0.058$.  We examine noise levels $\sigma$ in the range $10^{-2}\sigma_c\approx 0.0006$ to $10\sigma_c \approx 0.5774$. In terms of the noise ratio $\sigma_{NR}$, given that $\nrm{\Ubf^{(m),\star}}_{stdev'}\approx 1.16$, this corresponds to $\sigma_{NR}$ in the range $0.0005$ to $0.5$.

We simulate \eqref{burgers} on a fine grid $(\Xbf^f,\tbf^f)\subset [-1,1]\times [0,1.5]$ containing (2048,1801) points in space and time to mimic the continuum limit. We examine a range of resolutions, with the coarsest grid have $(64,57)$ points, or $32$ times coarser than $(\Xbf^f,\tbf^f)$. The raw data without noise is depicted in Figure \ref{sol} at the finest resolution with $\sigma=10\sigma_c$ along with the reference test function $\psi$ overlaid. We fix the reference test function such that $|\supp{\psi}|/|\Omega\times(0,T)| = 1/16$ (see Figure \ref{sol}), and we fix the differential operators $\partial^{\pmb{\alpha}} = (\partial_x^q)_{q=0,\dots,6}$ and the nonlinearities $\CalF = (u^q)_{q=0,\dots,6}$.

Figure \ref{cubicNLVB} (right column) shows recovery results as $m\to \infty$ for various noise levels, comparing raw data (solid lines) with filtered data (dot-dashed lines). As expected, for $\sigma<\sigma_c$ we see rapid recovery of the true support using both raw and filtered data. For $\sigma = \sigma_c$ with raw data (purple solid line), we see an initial increase in the probability of recovery, followed by a slow decline, which is to be expected since the continuum problem is not expected to yield correct recovery in this case. All noise levels $\sigma>\sigma_c$ using the raw data yield zero probability of correct support recovery. On the other hand, filtering the data successfully enables support recovery at all higher noise levels.

At noise levels $\sigma<\sigma_c$ the coefficient error $E_\infty(m)$ approximately decreases with Monte-Carlo rate of $\CalO(m^{-1/2})$ (bottom right panel of Figure \ref{cubicNLVB}). For higher noise levels, the $E_\infty$ decreases at an increased rate on filtered data, which is a reflection of the observation in Remark \ref{increase_rate_rmrk}.

\section{Discussion}\label{discussion}

In this work we have provided an in-depth analysis of the WSINDy algorithm in the limit of continuum data. We have presented results in as general a framework as possible, in particular showing convergence for non-smooth solutions of differential equations. This analysis includes several key results that can be used to improve algorithms for equation learning in general. In particular, we identify that WSINDy models may in  general be only conditionally consistent with the ground truth, and that spurious terms may arise in the limit if the noise level is above a critical threshold. Included in this result is the identification of a large class of models for which WSINDy {\it does} recover the correct model in the limit for any noise level, which explains the previously report  robustness of weak-form methods.

We also examine a filtered WSINDy approach, and identify that filtering the data leads to unconditional convergence for a wide range of systems and noise distributions. We propose that optimal filtering be a major priority for future algorithms, as the theoretical advantages of filtering are undeniable, yet is well-known that filtering can lead to much worse results if disadvantageous filters are use (see e.g.\ \cite{de2019filters}).

As far as other numerical results, we demonstrated that the theoretical findings are in fact exhibited in practice. In particular, filtering the data enables recovery of systems beyond the critical noise level. Also, crucially, for system which {\it do not} exhibit a critical noise level, we do observe unconditional convergence, no matter how high the orders of spatial derivatives are in true model (as in the Hyper-KS case).

We now offer several ideas for extensions and combination of our results with existing results.


\begin{enumerate}[label=(\Roman*)]
\item We have examined only models that can be integrated-by-parts to put all derivatives onto the test function. For other such terms, it is possible to use local-polynomial differentiation as in \cite{he2021asymptotic} and then apply the weak form to get convergence. We conjecture that local polynomial filtering may also provide advantages over the simple moving average filter employed here for general data filtering, hence a hybrid weak-form / local-polynomial discretization appears advantageous both theoretically and in practice.  

\item In this work we have restricted our analyses to cases where the true model is contained in the model library. A recent work \cite{russo2022convergence} has proven the convergence of WSINDy surrogate models using terms that are not restricted to the true support set. A fruitful future direction would be to merge these results and our current findings to extend results in \cite{russo2022convergence} to the case of discrete and noisy data. 

\item While results are presented here for sub-Gaussian noise distributions, extensions to heavier-tailed noise distributions are of course possible. For instance, the results above carry over immediately to sub-exponential noise albeit with slower concentration rates.

\item We have analyzed the MSTLS algorithm from \cite{messenger2020weakpde} due to its practical performance and speed. We also conjecture that thresholding-based sparse regression routines may offer an advantage in the case highly-correlated linear system. However, alternative sparse optimization techniques as employed in \cite{cortiella2020sparse} and \cite{bertsimas2022learning} may offer additional advantages with regard to the biases introduced by nonlinear functions in $\CalF$. On the other hand, suitably rescaling the data and coordinates as introduced in \cite{messenger2020weakpde} has the ability to increase the critical noise threshold $\sigma_c$ (e.g.\ by increasing the lower bound in equation \eqref{explicitboundsonsigmac}). Altogether, we believe that our results may inform future developments in sparse regression algorithms for equation learning.

\end{enumerate}

\section{Acknowledgments}
This research was supported in part by the NSF Mathematical Biology MODULUS grant 2054085, in part by the NSF/NIH Joint DMS/NIGMS Mathematical Biology Initiative grant R01GM126559, and in part by the NSF Computing and Communications Foundations grant 1815983. This work also utilized resources from the University of Colorado Boulder Research Computing Group, which is supported by the National Science Foundation (awards ACI-1532235 and ACI-1532236), the University of Colorado Boulder, and Colorado State University. The authors would also like to thank Prof.\,Vanja Duki\'c (University of Colorado at Boulder, Department of Applied Mathematics) for insightful discussions about the statistical aspects of this work and in particular the results concerning the bias.

\appendix

\section{Notation}\label{app:notation}
In Table \ref{notationtable} we include relevant notations used throughout along with text references.
\begin{table}
\begin{center}
{\footnotesize
\begin{tabular}{|c|c|c|c|}
\hline Symbol & Domain & Meaning & First Mention \\ 
\hline $\Omega$ & $\Rbb^d$ & spatial domain of true solution & Sec. \ref{sec:wsindyoverview} \\ 
\hline $(0,T)$ & $\Rbb$ & temporal domain of true solution & Sec. \ref{sec:wsindyoverview} \\ 
\hline $\CalH^{k,p}(D)$ & --- & spaces considered for true solutions & eq. \eqref{eq:Hspaces} \\
\hline $u$ & $\CalH^{k,\infty}(\Omega\times (0,T))$ & true solution & Sec. \ref{sec:wsindyoverview} \\
\hline $\CalF:=(f_j)_{j\in [J]}$ & $C(\Rbb^n,\Rbb)$ & trial functions & Sec. \ref{sec:wsindyoverview}\\
\hline $F^{\pmb{\omega}} = \{\exp(i\omega^Tu\}_{\omega \in \pmb{\omega}}$ & $C(\Rbb^n,\Rbb)$ & finite trigonometric basis & Sec. \ref{sec:modellibrary}\\
\hline $P^{(q)}$ & $C(\Rbb^n,\Rbb)$ & polynomials of degree at most $q$ over $\Rbb^n$ & Sec. \ref{sec:modellibrary}\\
\hline $\pmb{\alpha}:=(\alpha^s)_{s\in [S]}$ & $\Nbb^{d+1}$ & trial differential operators as multi-indices & Sec. \ref{sec:wsindyoverview}\\
\hline $S$ & $\Nbb$ & number of differential operators in library & Sec. \ref{sec:wsindyoverview}\\
\hline $J$ & $\Nbb$ & number of trial functions in library & Sec. \ref{sec:wsindyoverview}\\
\hline $\mathfrak{J}=SJ$ & $\Nbb$ & total number of library terms & Sec. \ref{sec:wsindyoverview}\\
\hline $\wstar$ & $\Rbb^{\mathfrak{J}\times n}$ & true model coefficients & Sec. \ref{sec:intro}\\
\hline $\psi$ & $C^{\nrm{\pmb{\alpha}}_\infty}_0(\Omega\times (0,T))$ & reference test function & Sec. \ref{sec:wsindyoverview} \\
\hline $\CalQ := \{(\xbf_k,t_k)\}_{k\in[K]}$ & $\Omega\times (0,T)$ & Query points used to evaluate convolutions with $\psi$ & Sec. \ref{sec:wsindyoverview} \\


\hline STLS & --- & sequential thresholding least squares & equation \eqref{STLS} \\

\hline $\CalL$ & $BV(0,\infty)$ & MSTLS loss function & eq. \eqref{lossfcn} \\

\hline $(\Delta x^{(m)},\Delta t^{(m)})$ & $\Rbb^2_+$ & spatiotemporal resolution at discretization level $m$ & Sec. \ref{sec:contprob} \\
\hline $(\Xbf^{(m)},\tbf^{(m)})$ & $\Rbb^{n_1^{(m)}\times \cdots \times n_d^{(m)}}\times \Rbb^{n_{d+1}^{(m)}}$ & spatiotemporal grid at discretization level $m$ & Sec. \ref{sec:contprob} \\
\hline $m$ & $\Nbb$ &  $\supp{\psi} \cap(\Xbf^{(m)},\tbf^{(m)})$ & Sec. \ref{sec:contprob} \\
\hline $\what^{(m)}$ & $\Rbb^{\mathfrak{J}\times n}$ & learned weight vector at discretization level $m$ & Sec. \ref{sec:intro}\\
\hline $\rho$ & $\CalP(\Rbb)$ & noise distribution & Sec. \ref{sec:sampling} \\
\hline $\Ubf^{(m)}$ & $(\Rbb^{n_1^{(m)}\times \cdots \times n_{d+1}^{(m)}})^n$ & noisy data evaluated at $(\Xbf^{(m)},\tbf^{(m)})$ & Sec. \ref{sec:contprob} \\
\hline $\ep$ & - & i.i.d.\ noise with distribution $\rho$ & Sec. \ref{sec:sampling} \\
\hline $(\Gbf^{(m)},\bbf^{(m)})$ & $\Rbb^{K\times\mathfrak{J}}\times \Rbb^{K\times n}$ & WSINDy linear system at discretization level $m$ & Sec. \ref{sec:contprob} \\
\hline $(\overline{\Gbf},\overline{\bbf})$ & $\Rbb^{K\times\mathfrak{J}}\times \Rbb^{K\times n}$ & WSINDy linear system in the limit $m\to \infty$ & Sec. \ref{sec:contprob} \\
\hline $(\Gbf^\star,\bbf^\star)$ & $\Rbb^{K\times\mathfrak{J}}\times \Rbb^{K\times n}$ & $(\overline{\Gbf},\overline{\bbf})$ in the case of noise-free data & Sec. \ref{sec:contprob} \\
\hline $\wbf^{(m),\lambda}$ & $\Rbb^{\mathfrak{J}\times n}$ & output of STLS$(\Gbf^{(m)},\bbf^{(m)},\lambda)$ & Sec. \ref{sec:contprob}\\
\hline $\what^{(m),\lambda}$ & $\Rbb^{\mathfrak{J}\times n}$ & output of MSTLS$(\Gbf^{(m)},\bbf^{(m)})$ & Sec. \ref{sec:contprob}\\
\hline $\overline{\wbf}^\lambda$ & $\Rbb^{\mathfrak{J}\times n}$ & output of STLS$(\overline{\Gbf},\overline{\bbf},\lambda)$ & Sec. \ref{sec:contprob}\\
\hline $\widehat{\overline{\wbf}}^\lambda$ & $\Rbb^{\mathfrak{J}\times n}$ & output of MSTLS$(\overline{\Gbf},\overline{\bbf})$ & Sec. \ref{sec:contprob}\\
\hline
\end{tabular}}
\end{center}
\caption{Notations used throughout.}
\label{notationtable}
\end{table}

\section{Convergence of the trapezoidal rule for $\CalH^{k,\infty}$}\label{app:trapconv}

Let $U\subset \Rbb^{n}$ be an open subset and $k>n/2$. Define $r = \left\lceil k-\frac{n}{2}\right\rceil-1$ and $\alpha = k-\frac{n}{2}-r$. By the Sobolev embedding theorem,  we have
\[H^k(U)\subset C^{r, \gamma}(U)\]
where $\gamma = \alpha$ if $\alpha\in (0,1)$, and if $\alpha=1$ the embedding holds for any $\gamma \in (0,1)$. The H\"older space $C^{r,\gamma}(U)$ is defined
\[C^{r,\gamma}(U) = \left\{f\in C(U)\ :\ L:=\max_{|\beta|\leq r}\sup_{x,y\in U}\frac{|\partial^\beta f(x)-\partial^{\beta}f(y)|}{|x-y|^\gamma}<\infty\right\}.\]
The following lemma verifies that the trapezoidal rule converges for functions $\partial^{\alpha^s}\psi(\cdot) f(u(\cdot))$ appearing in the WSINDy linear system, so long as $u \in \CalH^{k,\infty}(\Omega\times(0,T))$ and $\psi$, $f$ are smooth. For smoother functions the convergence rate is greater.
\begin{lemm}\label{trapconv}
Let $g\in \CalH^{k,\infty}(D)$ have compact support where $D\subset \Rbb^{n}$ is open and bounded. Then the multidimensional trapezoidal rule approximation of $\int_D g(x)dx$ over a grid $\Xbf\subset \overline{D}$ with uniform mesh width $h$ takes the form 
\[I_h(g) = h^n \sum_{x_i\in \Xbf} g(x_i).\]
Then we have the following cases:
\begin{enumerate}[label=(\roman*)]
\item If $k > \frac{n}{2}$ and we also have $g\in H^k(D)$, then for $C$ independent of $h$ it holds that
\[\left\vert I_h(g) - \int_D g(x) dx \right\vert \leq Ch^{k-\frac{n}{2}}.\]
\item If $\frac{n}{2} < k \leq \frac{n}{2}+1$, then for $C$ independent of $h$ it holds that
\[\left\vert I_h(g) - \int_D g(x) dx \right\vert \leq Ch^{k-\frac{n}{2}},\]
\item If $k > \frac{n}{2}+1$, then for $C$ independent of $h$ it holds that
\[\left\vert I_h(g) - \int_D g(x) dx \right\vert \leq Ch.\]
\end{enumerate}
\end{lemm}
\begin{proof}
Case $(i)$ follows from the compact support of $g$ using Fourier analysis. We now examine cases $(ii)$ and $(iii)$. By the definition of $\CalH^{k,\infty}(D)$, there exists a finite partition $D = D_1\cup\cdots \cup D_\ell$ where each subdomain $D_i$ has Lipschitz boundary. Define the interior boundary
\[\CalS_{\text{int}} = \left(\bigcup_{i=1}^L \partial D_i\right)\setminus \partial D\]
and the $\varepsilon$-tube of $\CalS_{\text{int}}$
\[\CalS_\varepsilon = \CalS_{\text{int}} +B_n(0,\varepsilon)\]
where $B_n(0,\varepsilon)$ is the ball of radius $\varepsilon$ in $\Rbb^n$ centered at $0$. Due to the regularity of the boundaries, the surface measure of $\CalS_{\text{int}}$ is finite, hence the volume of $\CalS_\varepsilon$ satisfies $\lim_{\varepsilon\to 0}|\CalS_\varepsilon| = 0$. Fix $h>0$ so that $D$ is partitioned into $M = |D|/h^n$ hypercubes and choose $\varepsilon < h/2$ so that $|\CalS_\varepsilon|\leq |\CalS_{\text{int}}|h$. Any given hypercube $K$ that intersects $\CalS_{\text{int}}$ contributes a worst-case error of 
\[\left\vert\frac{h^n}{2^n}\sum_{j=1}^{2^n}g_j - \int_{K}g(x)dx \right\vert \leq 2\nrm{g}_\infty h^n,\]
where $g_j=g(x^{(j)})$ for vertices $(x^{(j)})_{j=1}^{2^n}$ of $K$. There are fewer than $M\frac{h|\CalS_{\text{int}}|}{|D|}$ such hypercubes. 

For case $(ii)$, using the mean value theorem and the H\"older regularity, we have that integration over hypercubes $K$ with $K\cap \CalS_{\text{int}}=\emptyset$ satisfies
\[\left\vert\frac{h^n}{2^n}\sum_{j=1}^{2^n}g_j - \int_{K}g(x)dx \right\vert  = \left\vert\frac{h^n}{2^n}\sum_{j=1}^{2^n}\left(g_j - g(\tilde{x})\right) \right\vert  \leq \frac{Lh^n}{2^n}\sum_{j=1}^{2^n} |x^{(j)}-\tilde{x}|^\gamma \leq L\sqrt{n}^{k-\frac{n}{2}}h^{k+\frac{n}{2}}\]
where $g(\tilde{x}) = \frac{1}{|K|}\int_Kg(x)dx$. Here $L$ is the H\"older constant for $g$ in $C^{0,k-\frac{n}{2}}(D)$. This gives us the error
\[E \leq \left(M\frac{h|\CalS_{\text{int}}|}{|D|}\right) 2\nrm{g}_\infty h^n + M L\sqrt{n}^{k-\frac{n}{2}}h^{k+\frac{n}{2}}\]
\[\leq 2\nrm{g}_\infty|\CalS_\text{int}|h + L|D|\sqrt{n}^{k-\frac{n}{2}}h^{k-\frac{n}{2}} \leq C h^{k-\frac{n}{2}}.\]
For case $(iii)$, we have that $g\big\vert_{D_i}\in C^{1,\gamma}$ where $\gamma = \min(1,k-\frac{n}{2}-1)$ for $i=1,\dots,\ell$, hence we can exploit smoothness in the 1st derivative. In the 1D case ($n=1$), using that 
\[g(x) = g(0) + \int_0^xg'(y)dy = g(h) -\int_x^hg'(y)dy\]
we have the following:
\begin{align*}
\left\vert\frac{h}{2}(g(0)+g(h))-\int_0^hg(x)dx\right\vert &=\frac{1}{2}\left\vert\int_0^h\left(\int_0^xg'(y)dy-\int_x^hg'(y)dy\right)dx\right\vert\\
\text{(switch order of integration)}\qquad &=\frac{1}{2}\left\vert\int_0^h\int_0^y(g'(x)-g'(y))dxdy\right\vert \\
&\leq \frac{L}{2}\int_0^h\int_0^y(y-x)^{\gamma}dxdy \\
&= \frac{L}{2(\gamma+1)(\gamma+2)}h^{2+\gamma}.
\end{align*}
We can extend this to the trapezoidal rule in $n$ dimensions as follows. Let $K$ be a hypercube satisfying $K\cap \CalS_{\text{int}}=\emptyset$ with vertices $(x^{(j)})_{j=1}^{2^n}$. Then
\[\left\vert\frac{h^n}{2^n}\sum_{j=1}^{2^n}g_j - \int_{K}g(x)dx \right\vert = \frac{1}{2^n}\left\vert \int_{K} \left(\sum_{j=1}^{2^n} \int_0^1 \frac{d}{dt}g(z_j(t))dt \right)dx\right\vert \]
where each $z_j$ is a connected curve satisfying $z_j(0)=x^{(j)}$ and $z_j(1)=x$, using that
\[g(x) = g(x^{(j)}) + \int_0^1 \frac{d}{dt}g(z_j(t))dt, \qquad j=1,\dots,2^n.\]
Since this holds for any absolutely continuous curves $z_j$, we can select curves to cancel analogously to the 1D case. Let each $z_j$ take the path from $x_j$ to $x$ along the coordinate axes. For example, let $x^{(1)}=(0,\dots,0)$, $x^{(2)}=(h,0,\dots,0)$ and consider the paths
\[z_1: x^{(1)}\to (x_1,0,\dots,0) \to (x_1,x_2,0,\dots,0) \to \cdots \to x.\]
\[z_2: x^{(2)}\to (x_1,0,\dots,0) \to (x_1,x_2,0,\dots,0) \to \cdots \to x.\]
Then combining the first edges of each path in the sum above leads to the term 
\[I = \int_K\left(\int_0^{x_1} \partial_{x_1} g(y_1,0,\dots,0)dy_1 - \int_{x_1}^h\partial_{x_1}g(y_1,0,\dots,0)dy_1\right)dx \]
upon which we can use the same trick as before to swap the order of integration, leading to 
\[I\leq \frac{L}{(\gamma+1)(\gamma+2)}h^{2+\gamma + (n-1)}.\]
Each of the $2^n$ paths $z_j$ have $n$ edges, each of which gets paired with another edge to yield a similar bound, leading to an overall $n2^n / 2$ such error terms, or
\[\left\vert\frac{h^n}{2^n}\sum_{j=1}^{2^n}g_j - \int_{K}g(x)dx \right\vert \leq \frac{n 2^{n-1}}{2^n} \frac{L}{(\gamma+1)(\gamma+2)}h^{2+\gamma + (n-1)} = \frac{n}{2}\frac{L}{(\gamma+1)(\gamma+2)}h^{1+\gamma + n}.\]
This holds for all $K \cap \CalS_{int} = \emptyset$, and the overall error satisfies
\[E \leq \left(M\frac{h|\CalS_{\text{int}}|}{|D|}\right) 2\nrm{g}_\infty h^n + M \frac{n}{2}\frac{L}{(\gamma+1)(\gamma+2)}h^{1+\gamma + n} \]
\[\leq 2\nrm{g}_\infty|\CalS_\text{int}|h + \frac{n}{2}\frac{L|D|}{(\gamma+1)(\gamma+2)}h^{1+\gamma} \leq Ch.\]
\end{proof}

\section{Moment matrix lemmas}\label{app:comblemms}
\begin{lemm}\label{Ainv}
Define the Gaussian moment matrix
\begin{equation}\label{gaussmoment}
\Abf^{(p)}_{i,j} = \begin{dcases} 1, & i=j\\ {j\choose i}(j-i-1)!!\sigma^{j-i}, & j>i, \ (j-i)\ \text{even}  \\ 0, &\text{otherwise.}\end{dcases}
\end{equation}
Then the inverse of $\Abf^{(p)}$ is given by   
\[(\Abf^{(p)})^{-1}_{i,j} := (-1)^{\frac{j-i}{2}}\Abf^{(p)}_{i,j}.\]
\end{lemm}
\begin{proof}
Define $\Bbf^{(p)}_{i,j} := (-1)^{\frac{j-i}{2}}\Abf^{(p)}_{i,j}$ and consider an entry of the product $\Cbf_{i,j} = \left(\Abf^{(p)}\Bbf^{(p)}\right)_{i,j}$. We first note that $\Cbf_{i,j} = 0$ if either $j<i$ or $j-i$ is odd, the former due to the fact that $\Cbf$ is upper triangular, the latter because both $\Abf^{(p)}$ and $\Bbf^{(p)}$ have a ``checkerboard'' sparsity pattern with zeros on odd superdiagonals, a pattern which is inherited by $\Cbf$. Next since both $\Abf^{(p)}$ and $\Bbf^{(p)}$ have 1's along the diagonal, $\Cbf$ will also have 1's along the diagonal. For the remaining entries, in other words $\Cbf_{i,j}$ with $j-i$ an even positive integer,
\begin{align}
\Cbf_{i,j} &= \sigma^{j-i}\sum_{\substack{k=i \\k-i \text{ even}}}^j {k \choose k-i}{j \choose j-k}(k-i-1)!!(j-k-1)!!(-1)^\frac{j-k}{2}\\ 
\label{a11} &=\sigma^{j-i}{j \choose j-i}(j-i)!\sum_{\substack{k=i \\k-i \text{ even}}}^j \frac{(-1)^\frac{j-k}{2}}{(k-i)!!(j-k)!!}\\
&=\sigma^{j-i}{j \choose j-i}(j-i)!\sum_{\ell=0}^{\frac{j-i}{2}} \frac{(-1)^\frac{j-i-2\ell}{2}}{(2\ell)!!(j-i-2\ell)!!}\\
&=\left(-\frac{\sigma^{j-i}{j \choose j-i}(j-i)!}{2^\frac{j-i}{2}\left(\frac{j-i}{2}\right)!}(-1)^\frac{j-i}{2}\right)\left(\sum_{\ell=0}^{\frac{j-i}{2}} (-1)^\ell { \frac{j-i}{2} \choose \ell}\right)\\
&=0,
\end{align}
where we used the identities
\[{p \choose q}{p-q \choose j} = {p \choose q+j}{q+j \choose j}, \qquad p-q\geq j\geq 0\]
and
\[\sum_{\ell=0}^{p} (-1)^\ell { p \choose \ell}=0, \qquad p\geq 1.\]
This shows that $\Cbf$ is the identity, so that $\Bbf^{(p)} = (\Abf^{(p)})^{-1}$.
\end{proof}

\begin{lemm}\label{app:lemm_bnd_on_AL}
Let $\Abf^{(p)} = \Ibf + \Lbf^{(p)}$ be the Gaussian moment matrix defined in \eqref{gaussmoment}. Then it holds that 
\[\max\{\nrm{\Abf^{(p)}}_\infty,\nrm{\Abf^{(p)}}_2\}\leq \nrm{\Abf^{(p)}}_1 \leq \exp\left(\sigma^2{p \choose 2}\right)\]
and 
\[\max\{\nrm{\Lbf^{(p)}}_\infty,\nrm{\Lbf^{(p)}}_2\}\leq \nrm{\Lbf^{(p)}}_1 \leq \sigma^2{p \choose 2}\exp\left(\sigma^2{p \choose 2}\right).\]
\end{lemm}

\begin{proof}
Set $\alpha = \sigma^2{p \choose 2}$. First, note that the maximum entry along any superdiagonal of $\Abf^{(p)}$ occurs in the $p$th column, i.e.
\[\max_{j-i = k} \Abf^{(p)}_{ij} = \Abf^{(p)}_{p-k,p}:=y_k \] 
where for $k$ even,
\[y_k = {p \choose p-k}(k-1)!!\sigma^{k} = \frac{{p \choose p-k}(k-1)!!}{{p \choose 2}^{\frac{k}{2}}}\alpha^{k/2}\]
and for $k$ odd, we have $y_k = 0$. This implies that the $p$th column achieves the maximum column sum, so that $\nrm{\Abf^{(p)}}_1 = \nrm{\ybf}_1$. This also implies that we can upper-bound any row sum by $\nrm{\ybf}_1$, hence $\nrm{\Abf^{(p)}}_\infty\leq \nrm{\Abf^{(p)}}_1$. Now let $\Tbf(\ybf)$ be the Toeplitz matrix formed by the vector $\ybf$:
\[\Tbf(\ybf) = \begin{bmatrix} y_0 & y_1 & \cdots & y_p \\ y_p & y_0 & \cdots & y_{p-1} \\  \vdots & \ddots & \ddots & \vdots \\ y_1 & y_2 & \cdots & y_0 \end{bmatrix}.\]
Since the entries of $\Abf^{(p)}$ are positive, we have
\[\nrm{\Abf^{(p)}}_2\leq \nrm{\Tbf(\ybf)}_2 = \nrm{\ybf}_1 = \nrm{\Abf^{(p)}}_1,\]
where $\nrm{\Tbf(\ybf)}_2 = \nrm{\ybf}_1$ follows from the fact that the product $\Tbf(\ybf)\xbf$ is the convolution of $\xbf$ with $\ybf_-\ =\ \ybf$ in reverse order, so by Young's inequality for convolutions have 
\[\nrm{\Tbf(\ybf)}_2 =\max_{\nrm{\xbf}_2=1}\nrm{\xbf*\ybf_-}_2\leq \nrm{\ybf}_1,\]
and $\nrm{\Tbf(\ybf)\xbf}_2 = \nrm{\ybf}_1$ is achieved by $\xbf = (1,\dots,1)^T$. Revisiting the entries of $\ybf$, we then have for $k$ even,
\[y_k \leq \frac{2^{\frac{k}{2}}}{(k)!!}\alpha^{k/2} = \frac{\alpha^{k/2}}{\left(\frac{k}{2}\right)!}\]
so that  
\[\nrm{\ybf}_1 \leq \sum_{\substack{k=0 \\ k\text{ even}}}^p \frac{\alpha^{k/2}}{\left(\frac{k}{2}\right)!} = \sum_{\ell=0}^{\left\lfloor p/2\right\rfloor}\frac{\alpha^\ell}{\ell!}\leq \exp(\alpha).\]
The inequalities for $\Lbf^{(p)}$ follow similarly, noting that 
\[\max\{\nrm{\Lbf^{(p)}}_\infty,\nrm{\Lbf^{(p)}}_2\}\leq \nrm{\Lbf^{(p)}}_1 = \sum_{k=2}^py_k = \alpha \sum_{\ell=0}^{\left\lfloor p/2\right\rfloor}\frac{\alpha^\ell}{\ell!}\frac{1}{\ell+1}\leq \alpha\exp(\alpha).\]
\end{proof}
For general moment matrices we have the following.

\begin{lemm}\label{generalmomentmatrix}
Let $\Abf^{(p)}$ be the moment matrix of order $p$ for probability distribution $\rho$ such that the monomials $P^{(p)} = \{1,x,\dots,x^p\}$ transform under cross-correlation with $\rho$ according to $P^{(p)} \star \rho = P^{(p)}\Abf^{(p)}$. Then the inverse of $\Abf^{(p)}$ is given by $(\Abf^{(p)})^{-1}_{ij} = f(j-i)\Abf^{(p)}_{ij}$ where $f:\Nbb\cup\{0\}\to \Rbb$ obeys the recurrence
\begin{equation}\label{recrel}
f(q) = -\sum_{\ell=0}^{q-1}{q\choose \ell} \left(\frac{M_{q-\ell}M_\ell}{M_q}\right) f(\ell), \qquad f(0) = 1
\end{equation}
for moments $M_k := \int_\Rbb x^k d\rho(x)$. Moreover, if $\rho$ is symmetric and sub-Gaussian,
then with $\tilde{\Lbf}^{(p)} := (\Abf^{(p)})^{-1} - \Ibf^{(p)}$ it holds that 
\[\|\tilde{\Lbf}^{(p)}\| =\CalO(\|\rho\|^2_\text{SG}).\]
\end{lemm}

\begin{proof}
$\Abf^{(p)}$ is defined by
\begin{equation}
\Abf^{(p)}_{ij} = \begin{dcases} M_{j-i}{j \choose i}, &j\geq i \\0, & \text{otherwise}.\end{dcases}
\end{equation}
Defining $(\Bbf^{(p)})_{ij} := f(j-i)\Abf^{(p)}_{ij}$, we see that solving for $f$ such that
\begin{align*}
\delta_{ij} &= (\Abf^{(p)}\Bbf^{(p)})_{ij} \\
&=\sum_{k=i}^j {k \choose i}{j \choose k}M_{k-i}M_{j-k}f(j-k)\\
&= {j \choose i}\sum_{\ell=0}^q{q\choose \ell}M_{q-\ell}M_\ell f(\ell)
\end{align*}
leads directly to \eqref{recrel}, noting that only $j\geq i$ needs to be considered. 

Now consider $\rho$ to be sub-Gaussian. It then holds that
\[|\tilde{\Lbf}^{(p)}_{ij}| = \begin{dcases} |f(j-i)|M_{j-i}{j\choose i}&, j-i \text{ even}, j\geq i+2. \\ 0&,\text{otherwise}.\end{dcases}\]
If we now assume that $\|\rho\|_\text{SG}\leq B$ for any $B>0$, it holds through sub-Gaussianity that 
\[M_q\leq \frac{q!!}{2^{{q}/{2}-1}}\|\rho\|_\text{SG}^q\leq C'\|\rho\|_\text{SG}^2, \qquad q\geq 2\]
where $C'$ depends only on $B$ and $q$. Using the recurrence \eqref{recrel} and Jensen's inequality we can bound $f$ independently of $\rho$,
\[|f(q)|\leq \sum_{\ell=0}^{q-1}{q\choose \ell}|f(\ell)|,\]
hence we see that 
\[|\tilde{\Lbf}^{(p)}_{ij}| \leq C\|\rho\|_\text{SG}^2\]
where $C$ depends only on $B$ and $p$. In this way we have that 
\[\|\tilde{\Lbf}^{(p)}\| = \CalO(\|\rho\|_\text{SG}^2)\]
for any norm $\|\cdot\|$.
\end{proof}

\section{Concentration results}\label{app:conclemmas}

\begin{lemm}[Lemma 3.1]
Let $f:\Rbb\to \Rbb$ satisfy
\begin{equation}
|f(x)|\leq C_f\left(1+|x|^p\right),
\end{equation}
for some $p\geq 1$ and $C_f>0$. For bounded sequences of real numbers $(\alpha_i)_{i\in \Nbb}$ and $(u_i)_{i\in \Nbb}$, and for $\ep_i\sim \rho$ i.i.d.\, $i\in \Nbb$, define the random variables $Y_i = \alpha_i f(u_i +\ep_i)$. Then for any $\kappa\geq \nrm{\rho}_{\text{SG}}$ and $t>0$, it holds that the right tails of $Y_i$ are captured by a common rate function $I(t)$,
\[\Pbb\left(Y_i > t\right)\leq \exp (-I(t)),\]
where 
\begin{equation}\label{eq:rate_fcn_app}
I(t):=\ind{(t^*,\infty)}(t)\left[\frac{1}{\kappa^2}\left(\left(\frac{t}{C_f\alpha^*}-1\right)^{1/p}-u^*\right)^2-\log(2)\right] = \frac{t^{2/p}}{\kappa^2(C_f\alpha^*)^{2/p}}I_0(t),
\end{equation}
for $\alpha^* = \sup_i|\alpha_i|$, $u^* = \sup_i|u_i|$, and $t^*:= C_f\alpha^*\left(1+\left(u^*+ \kappa\sqrt{\log(2)}\right)^p\right)$. Moreover, $I_0(t)$ is monotonically increasing from $0$ to $1$ over $t\in(t^*,\infty)$, and is defined in the proof.
\end{lemm}
\begin{proof}
Using \eqref{eq:growthcondition} and the symmetry of $\rho$, we get for each $i$ and every $t>t^*$,
\[\Pbb\left(Y_i >t\right)\leq \Pbb\left(C_f\alpha^*\left(1+|u_i+\ep|^p\right)> t \right) \leq 2\Pbb\left(\ep>\left(\frac{t}{C_f\alpha^*}-1\right)^{1/p}-u^*\right)\]
from which sub-Gaussianity implies 
\[\Pbb\left(Y_i >t\right) \leq \exp(-I(t))\]
for $I(t)$ defined as above. This leads to
\[I_0(t) = \ind{(t^*,\infty)}(t)\left[\left(\left(1-\frac{C_f\alpha^*}{t}\right)^{1/p}-u^*\left(\frac{C_f\alpha^*}{t}\right)^{1/p}\right)^2 -\frac{\kappa^2(C_f\alpha^*)^{2/p}\log(2)}{t^{2/p}}\right]\]
where the quantity in brackets is positive over $t>t^*$ due to the definition of $t^*$, and approaches 1 as $t\to \infty$. It is monotonically increasing as it contains only sums and compositions of functions that preserve monotonicity over $t>t^*$ (namely, $t\to A-Bt^{-\theta}$ is monotonically increasing for all $A,B,\theta \geq 0$, and sums of monotonically increasing functions are monotonically increasing).  
\end{proof}

\begin{lemm}[Lemma 3.2]
Let $Y_i$ be defined under the same conditions as Lemma \ref{lemm1} and choose $\beta\in(0,1)$. Then there exists $\overline{v}(\beta) < \infty$ such that the sum $S_m = \sum_{i=1}^{m} Y_i$ satisfies
\begin{equation}\label{eq:convergence_of_sums_app}
\Pbb\left(|S_m-\Ebb S_m|> mt\right) \leq  \begin{dcases} 2 \exp\left(-\frac{\beta}{2} I(mt)\right) + 2m\exp\left(-I(mt)\right), & t \geq t_m(\beta)\\
2 \exp\left(-\frac{mt^2}{2\overline{v}(\beta)}\right) + 2m\exp\left(-\frac{mt_m(\beta)^2}{\overline{v}(\beta)}\right), & 0\leq t < t_m(\beta), \end{dcases}
\end{equation}
where $t_m (\beta):= \sup\{t\geq 0\ :\ t\leq \beta\overline{v}(\beta)\frac{I(mt)}{mt}\}$.
\end{lemm}
\begin{proof}
This follows from a simple modification of Theorem 1 in \cite{bakhshizadeh2020sharp} to the case of independent but not identically distributed random variables. As in \cite{bakhshizadeh2020sharp} we define 
\[Y_i^L= Y_i\ind{Y_i\leq L}\]
and use that (\cite[Lemma 1]{bakhshizadeh2020sharp})
\[\Ebb\left[\exp\left(\lambda(Y_i^L-\Ebb[Y_i])\right)\right]\leq \frac{k_i(L,\lambda)}{2}\lambda^2\]
for all $\lambda,L>0$, where
\[k_i(L,\lambda) := \Ebb\left[\left(Y_i^L-\Ebb[Y_i]\right)^2\ind{Y_i^L\leq \Ebb[Y_i]} + \left(Y_i^L-\Ebb[Y_i]\right)^2\exp\left(\lambda\left(Y_i^L-\Ebb[Y_i]\right)\right)\ind{Y_i^L> \Ebb[Y_i]}\right].\]
We then note that $\sup_{i\in \Nbb} k_i(L,\lambda)=: \overline{k}(L,\lambda)$ is bounded, for instance by 
\[\overline{k}(L,\lambda)\leq \sup_i\Vbb[Y_i]+(L+\mu^*)^2\exp(\lambda(L+\mu^*))),\] 
where $\mu^*:=\sup_i|\Ebb[Y_i]|$. We then modify the proof of \cite[Theorem 1]{bakhshizadeh2020sharp} as follows,
\begin{align*}
\Pbb\left(S_m-\Ebb[S_m]>mt\right)&\leq \Pbb\left(\sum_{i=1}^m Y_i^L-\Ebb[S_m] > mt\right)+\Pbb\left(\exists i\ \text{ s.t. }\ Y_i> L\right)\\
&\leq \exp(-\lambda m t)\prod_{i=1}^m\Ebb\left[\exp(\lambda(Y_i^L-\Ebb[Y_i]))\right] + m\Pbb\left(Y_i>L\right)\\
&\leq \exp\left(m\left(-\lambda t + \frac{\overline{k}(L,\lambda)}{2}\lambda^2\right)\right) + m\exp\left(-I(L)\right).
\end{align*}
Choosing $\lambda=\beta\frac{I(mt)}{mt}$, $L=mt$ for $t>t_m(\beta)$ and $\lambda = \frac{t}{\overline{v}(\beta)}$, $L=mt_m(\beta)$ for $t\leq t_m(\beta)$ gives \eqref{eq:convergence_of_sums_app}, where $\overline{v}(\beta) = \sup_{L>0}\overline{k}(L,\beta\frac{I(L)}{L})$. We also used that $1-\frac{\beta\overline{v}(\beta)I(mt)}{2mt^2}\in[\frac{1}{2},1]$ for $t\geq t_m(\beta)$.\\
Finally, the bound for $\overline{v}(\beta)$ can be obtained using \cite[Lemma 3]{bakhshizadeh2020sharp}, which asserts that
\[k_i\left(L,\beta\frac{I(L)}{L}\right) \leq \sup_i\Vbb[Y_i]+\exp\left(-\beta \Ebb[Y_i]\frac{I(L)}{L}\right)\int_0^{L-\Ebb[Y_i]}\exp(-(1-\beta)I(t+\Ebb[Y_i]))(2t+\beta tI(t))dt.\]
Letting $C_1 = \sup_i\Vbb[Y_i]$, $C_2 = \sup_{i\in \Nbb,L>0}\exp\left(-\beta \Ebb[Y_i]\frac{I(L)}{L}\right)$, and introducing a parameter $\gamma>1$, we can then use the definition of $I(t)$ to get
\begin{align*}
k_i\left(L,\beta\frac{I(L)}{L}\right) &\leq C_1+C_2\int_0^\infty\exp(-(1-\beta)I(t-\mu^*))(2t+\beta tI(t))dt\\
&\leq C_1+C_2\int_0^{\gamma t^*+\mu^*}(2t+\beta tI(t))dt+\int_{\gamma t^*+\mu^*}^\infty \exp(-(1-\beta)I(t-\mu^*))(2t+\beta tI(t))dt\\
&\leq C'+C_2\int_{\gamma t^*}^\infty \exp\left(-\frac{(1-\beta)I_0(\gamma t^*)}{K^2(C_f\alpha^*)^{2/p}}s^{2/p}\right)r(s)ds \\
&\leq C'+C_2\int_0^\infty \exp\left(-A s^{2/p}\right)r(s)ds,
\end{align*}
where 
\[r(s) = 2(s+\mu^*)+\frac{1}{\kappa^2(C_f\alpha^*)^{2/p}}(s+\mu^*)^{2/p+1} \leq r_0+r_1 s +r_2 s^{2/p+1}.\]
We see from this that $\overline{v}(\beta)=\sup_{i\in \Nbb,L>0} k_i(L,\beta\frac{I(L)}{L})$ is finite, since $\gamma>1 $ and $\beta\in(0,1)$ imply $A>0$, and $s\to \exp(-As^{2/p})$ has finite moments of all order over $s\in[0,\infty)$.
\end{proof}

\begin{thm}[Theorem 3.1]
Suppose each function in library $\CalF$ satisfies the growth bound \eqref{eq:growthcondition} for some $p:=p_{\max}$. Then it holds that for every $t>\overline{t}(m)$, where $\overline{t}(m)\to 0$ as $m\to \infty$, we have the concentration rates
\begin{equation}\label{eq:convergence_of_G_app}
\Pbb\left(\nrm{(\Gbf^{(m)},\bbf^{(m)}) - (\overline{\Gbf},\overline{\bbf})}_{\overrightarrow{\infty}}> t \right) \leq  \begin{dcases} K\mathfrak{J}\exp\left(-\frac{c}{2} (mt)^{2/p_{\max}}\right) + K\mathfrak{J}m\exp\left(-c (mt)^{2/p_{\max}}\right), & t \geq t_m\\
K\mathfrak{J} \exp\left(-\frac{mt^2}{2\overline{v}}\right) + K\mathfrak{J}m\exp\left(-\frac{mt_m^2}{\overline{v}}\right), & 0\leq t < t_m, \end{dcases}
\end{equation}
where the rate factor $c$ depends on $\nrm{u}_\infty$, $|\Omega\times (0,T)|$, $\pmb{\alpha}$, $\psi$, $\CalF$, and $\nrm{\rho}_{\text{SG}}^2$.
\end{thm}
\begin{proof}
For convenience we have chosen $\beta=0.5$ from Lemma \ref{lemm2}, so that $t_m  =t_m(0.5)$, $\overline{v}=\overline{v}(0.5)$. 

It suffices to consider only the concentration of $\Gbf^{(m)}$ to $\overline{\Gbf}$ since entries of $\bbf^{(m)}$ are of the same type\footnote{In addition, when $\bbf^{(m)}$ is linear in the data (e.g.\ when the left-hand side of \eqref{diffform} is a linear differential operator), $\bbf^{(m)}$ concentrates at a Gaussian rate to $\overline{\bbf}$, and $\overline{\bbf}=\bbf^\star$}. First notice that 
\[\nrm{\Gbf^{(m)}-\overline{\Gbf}}_{\overrightarrow{\infty}}\leq \underbrace{\nrm{\Gbf^{(m)}-\Ebb\Gbf^{(m)}}_{\overrightarrow{\infty}}}_{\text{deviation from mean}}+\underbrace{\nrm{\Ebb\Gbf^{(m)}-\overline{\Gbf}}_{\overrightarrow{\infty}}}_{\text{integration error}}\]
and similarly for $\bbf^{(m)}$. The entries of $\Gbf^{(m)}$ satisfy
\[\Gbf^{(m)}_{k,(s-1)J+j} = \sum_{\substack{(\xbf_\ell,t_\ell)\in \\ \supp{\psi(\xbf_k-\cdot,t_k-\cdot)}\cap (\Xbf^{(m)},\tbf^{(m)})}} \partial^{\alpha^s}\psi(\xbf_k-\xbf_\ell,t_k-t_\ell)f_j(\Ubf^{(m)}(\xbf_\ell,t_\ell))(\Delta x^{(m)})^d\Delta t^{(m)}\]
\[= \frac{1}{m}\sum_{\ell=1}^{m} \alpha_{s,k,\ell} f_j(u_\ell+\ep_\ell)\] 
where
\[\alpha_{s,k,\ell} = |\supp{\psi}|\partial^{\alpha^s}\psi(\xbf_k-\xbf_\ell,t_k-t_\ell), \quad u_\ell = u(\xbf_\ell,t_\ell),\] 
$m$ is the number of points $\psi$ is supported on at the resolution $(\Delta x^{(m)},\Delta t^{(m)})$. By the smoothness and compact support of $\psi$, together with the regularity $u\in \CalH^{k,\infty}(\Omega\times (0,T))$, we have that $\alpha^*:=\sup_{s,k,\ell} |\alpha_{s,k,\ell}|<\infty$ and $u^*:=\sup_\ell |u_\ell|<\infty$. Hence, each entry $\Gbf^{(m)}_{k,(s-1)J+j}$ concentrates to $\Ebb \Gbf^{(m)}_{k,(s-1)J+j}$ according to Lemma \ref{lemm2}, and in particular its concentration can be modelled by a common rate function $I(t)$ given by \eqref{eq:rate_fcn_app}. To make use of the desired $\CalO(\exp(-(mt)^{2/p_{\max}}))$ concentration, we must take $mt>t^*$. Further, $t$ must be larger than the integration error, which is at-most $\CalO(m^{-\min(1,k-(d+1)/2)})$, hence depends on the smoothness of $u$ relative to the spatiotemporal dimension $d+1$. Denote the integration error by $e_\text{int}(m)$. Then for some $\gamma>1$, taking 
\begin{equation}
\overline{t}(m) = \max\left(\frac{\gamma t^*}{m},e_\text{int}(m)\right)
\end{equation}
we arrive at the desired concentation rates for $t>\overline{t}(m)$, and with an at-worst rate factor
\[c = \frac{I_0(\gamma t^*)}{4\nrm{\rho}_{\text{SG}}^2(C_f \alpha^*)^{2/p_{\max}}}\]
Finally, each $\Gbf^{(m)}$ has $K\mathfrak{J}$ entries, hence a union bound provides the desired concentration result. 
\end{proof}

\begin{thm}[Theorem 4.2]
Let the above assumption hold with $m^{(\nu)}\gtrsim m^{\alpha}$ for some $\alpha\in(0,1)$ and $\CalF$ locally Lipschitz with polynomial growth, i.e. $\forall f\in \CalF$ and $(x,y)\in \Rbb$,
\[|f(x)-f(y)|\leq C_\CalF|x-y|\left(1+|x-y|^{p_{\max}-1}\right)\]
for some $C_\CalF\geq 0$ and $p_{\max}\geq 1$, and recall that $u\in \CalH^{k,\infty}(D)$. Then for $t> \mathfrak{t}(m)$, where $\mathfrak{t}(m) = \CalO((m^{-\alpha \left(\frac{k}{d+1}-\frac{1}{2}\right)})$, and sufficiently large $m$, we have 
\[\Pbb \left(\|(\tilde{\Gbf}^{(m)},\tilde{\Gbf}^{(m)}) - (\Gbf^\star,\bbf^\star)\|_{\overrightarrow{\infty}}>t\right) \leq 3K\mathfrak{J}\exp\left(-\frac{c}{2^{1+p_{\max}}}[m(t-\mathfrak{t}(m))]^{2/{p_{\max}}}\right).\]
where $c$ is the same rate from Theorem \ref{matrix_concentration}.
\end{thm}

\begin{proof}
As in Theorem \ref{matrix_concentration}, it suffices to consider concentration of $\tilde{\Gbf}^{(m)}$. Let $\pmb{\nu}^{(m)}$ be a simple moving average filter with $n^\nu_q$ points in dimension $q$ and filter width $m^{(\nu)} := \prod_{q=1}^{d+1}n^\nu_q$ satisfying $m^{(\nu)} = m^\alpha$ for some $\alpha \in (0,1)$. Denote the filtered data, used to build $\tilde{\Gbf}^{(m)}$, by $\tilde{\Ubf}^{(m)} = \pmb{\nu}^{(m)}* \Ubf^{(m)} = \tilde{\Ubf}^{(m),\star}+\tilde{\ep}^{(m)}$, where $\tilde{\Ubf}^{(m),\star} = \pmb{\nu}^{(m)}*\Ubf^{(m),\star}$ is the filtered clean data and the filtered noise $\tilde{\ep}^{(m)}$ is mean zero and correlated accorded to 
\[\Ebb\tilde{\ep}^{(m)}_i\tilde{\ep}^{(m)}_j = \frac{\sigma^2}{m^\alpha}\Sigma(i,j), \qquad \Sigma(i,j) = \prod_{q=1}^{d+1}\max\left(1-\frac{|i_q-j_q|}{n^\nu_q},\ 0\right)\in[0,1],\]
where we've treated the indices $i$ and $j$ as vectors in $\Rbb^{d+1}$. We will couple $\tilde{\Ubf}^{(m)}$ to a filtered, {\it uncorrelated} dataset $\widehat{\Ubf}^{(m)}$ which we inject as an intermediary, defined as $\widehat{\Ubf}^{(m)} = \tilde{\Ubf}^{(m),\star}+\widehat{\ep}^{(m)}$ where $\widehat{\ep}$ and $\tilde{\ep}$ are all identically distributed according to $\rho^{(m)}$ with variance $\frac{\sigma^2}{m^\alpha}$, yet $\widehat{\ep}$ satisfy
\[\Ebb\widehat{\ep}^{(m)}_i\widehat{\ep}^{(m)}_j = \delta_{ij}\frac{\sigma^2}{m^\alpha}\]
and 
\[\Ebb\widehat{\ep}^{(m)}_i\tilde{\ep}^{(m)}_j = 0 \quad \forall i,j\in \Rbb^{d+1}.\]
In other words, entries $\widehat{\ep}_i$ are independent, and are independent from all $\tilde{\ep}_j$. 

We will split the error into the following terms:
\[\|\tilde{\Gbf}^{(m)} - \Gbf^\star\| \leq \|\tilde{\Gbf}^{(m)} - \widehat{\Gbf}^{(m)}\| + \|\widehat{\Gbf}^{(m)} - \overline{\Gbf}^{(m)}\| +\|\overline{\Gbf}^{(m)} - \overline{\Gbf}^{(m),\star}\|+\|\overline{\Gbf}^{(m),\star} - \Gbf^{(m),\star}\|+\|\Gbf^{(m),\star} - \Gbf^\star\|\]
\[ = I_1+I_2+I_3+I_4+I_5.\]
Using subscript $m$ to denote trapezoidal rule approximation of inner products on the grid $(\Xbf^{(m)}, \tbf^{(m)})$, the different intermediate matrices are defined by
\begin{itemize}
\item $\tilde{\Gbf}^{(m)}_{ksj} = \lan \partial^{\alpha^s}\psi_k, f_j(\tilde{\Ubf}^{(m),\star}+\tilde{\ep}^{(m)})\ran_m$ \\
\item $\widehat{\Gbf}^{(m)}_{ksj} = \lan \partial^{\alpha^s}\psi_k, f_j(\tilde{\Ubf}^{(m),\star}+\widehat{\ep}^{(m)})\ran_m$ \\
\item $\overline{\Gbf}^{(m)}_{ksj} = \lan \partial^{\alpha^s}\psi_k, \rho^{(m)}\star f_j(\tilde{\Ubf}^{(m),\star})\ran_m$ \\
\item $\overline{\Gbf}^{(m),\star}_{ksj} = \lan \partial^{\alpha^s}\psi_k, f_j(\tilde{\Ubf}^{(m),\star})\ran_m$ \\
\item $\Gbf^{(m),\star}_{ksj} = \lan \partial^{\alpha^s}\psi_k, f_j(\Ubf^{(m),\star})\ran_m$ \\
\item $\Gbf^\star_{ksj} = \lan \partial^{\alpha^s}\psi_k, f_j(u)\ran$.
\end{itemize}

$I_4$: First, we note that $I_4=\|\overline{\Gbf}^{(m),\star} - \Gbf^{(m),\star}\|$ is determined by how fast the locally averaged clean data $\tilde{\Ubf}^{(m),\star}$ converges to the clean data $\Ubf^{(m),\star}$. The locally Lipschitz assumption on the library and the smoothness of $u\in \CalH^{k,\infty}$ with $k>\frac{d+1}{2}$ implies that there exists $\tilde{C}$ depending on $C_\CalF$, $\psi$, and $\nrm{u}_\infty$ such that $m$ large enough implies 
\[I_4 \leq \tilde{C} \|\tilde{\Ubf}^{(m),\star}-\Ubf^{(m),\star}\| \leq \tilde{C} m^{-\alpha \min\{\frac{1}{d+1},\frac{k}{d+1}-\frac{1}{2}\}}.\]
Further, this is a worst-case bound reached only when $u$ is both discontinuous and not continuously differentiable in the subregions of $D=\Omega\times(0,T)$ where $u$ is continuous. If $u\in C^{1,\gamma}(D)$ then the rate changes to $\CalO(m^{-\frac{\alpha}{d+1}(1+\gamma)})$.

$I_5$: Asymptotic behavior of $I_5=\|\Gbf^{(m),\star} - \Gbf^\star\|$ is similar to $I_4$ as it is determined by how quickly the trapezoidal rule converges on $\partial^{\alpha^s}\psi(\cdot) f_j(u(\cdot))$. The worst-case rate is given by 
\[I_5\leq \tilde{C} m^{-\min\{\frac{1}{d+1},\frac{k}{d+1}-\frac{1}{2}\}},\]
which similary improves with the smoothness of $u$. 

$I_3$: For $I_3=\|\widehat{\Gbf}^{(m)} - \overline{\Gbf}^{(m)}\|$, we can consider the pointwise rate of convergence
\[\rho^{(m)}\star f_j(x)\to f_j(x),\]
which for $f_j$ locally Lipschitz and with polynomial growth, we get 
\[|\rho^{(m)}\star f_j(x)-f_j(x)|\leq C_\CalF\left(\sqrt{M_2(\rho^{(m)})} + M_{p_{\max}}(\rho^{(m)})\right)\leq C_\CalF C_\sigma\frac{\sigma}{m^{\alpha/2}}\]
where we used Jensen's inequality and sub-Gaussianity of $\rho^{(m)}$ to relate higher moments to the variance, and $C_\sigma$ depends on $p_{max}$ and $\rho$ but not $m$. It should be noted that this rate also increases to $\CalO(m^{-\alpha})$ for polynomial and trigonometric libraries $\CalF$. 

Hence for $m$ large enough, it suffices to consider 
\[\Pbb\left(I_1+I_2 > t'\right) \leq \Pbb\left(I_1 > t'/2\right)+\Pbb\left(I_2 > t'/2\right) := E_1+E_2,\]
where $t' = t-\sum_{i=3}^5I_i$, and for $k>1+\frac{d+1}{2}$, $t' = t - \tilde{C}m^{-\frac{\alpha}{d+1}}$. For $E_2$ we are in the same situation as in Theorem \ref{matrix_concentration}, hence for $t'>0$ we can use Corollary \ref{matrix_corr} to get the exponential concentration:
\[E_2 = \Pbb\left( \|\widehat{\Gbf}^{(m)} - \overline{\Gbf}^{(m)}\|_{\overrightarrow{\infty}} > t'/2\right) \leq 2K\mathfrak{J}\exp\left(-\frac{c}{2^{1+2/{p_{\max}}}}(mt')^{2/{p_{\max}}}\right).\]
This leaves $E_1$, which we now show also yields exponential concentration, but at a sub-Gaussian rate. Indeed,

\begin{align*}
\Pbb\left(|\tilde{\Gbf}^{(m)}_{ksj} - \widehat{\Gbf}_{ksj}^{(m)}|> t'/2\right) &= \Pbb\left(|\lan \partial^{\alpha^s}\psi_k, f_j(\tilde{\Ubf}^{(m),\star}+\tilde{\ep}^{(m)})-f_j(\tilde{\Ubf}^{(m),\star}+\widehat{\ep}^{(m)})\ran_m| > t'/2\right) \\ 
&\leq \Pbb\left(\nrm{\partial^{\alpha^s}\psi_k}_\infty\nrm{f_j(\tilde{\Ubf}^{(m),\star}+\tilde{\ep}^{(m)})-f_j(\tilde{\Ubf}^{(m),\star}+\widehat{\ep}^{(m)})}_{\overrightarrow{1}} > mt'/(2\gamma T|\Omega|)\right)\\
&\leq \sum_{i=1}^m \Pbb\left(|f_j(\tilde{\Ubf}^{(m),\star}_i+\tilde{\ep}^{(m)}_i)-f_j(\tilde{\Ubf}^{(m),\star}_i+\widehat{\ep}_i^{(m)})| > Ct'\right)
\end{align*}
where $C = (2\gamma T|\Omega| \nrm{\partial^{\alpha^s}\psi_k}_\infty)^{-1}$. (Here we replaced the volume element $(\Delta x)^d\Delta t$ with the equivalent expression $\gamma\frac{T|\Omega|}{m}$ using that $|\supp{\psi}| = \gamma T|\Omega|$, for some $\gamma<1$.) Now, define the sets
\[A_i = \{|f_j(\tilde{\Ubf}^{(m),\star}_i+\tilde{\ep}^{(m)}_i)-f_j(\tilde{\Ubf}^{(m),\star}_i+\widehat{\ep}^{(m)}_i)| > Ct'\},\]
\[B_i = \{|\tilde{\ep}^{(m)}_i - \widehat{\ep}_i^{(m)}|>1\}.\]
We get 
\begin{align*}
\Pbb\left(|f_j(\tilde{\Ubf}^{(m),\star}_i+\tilde{\ep}^{(m)}_i)-f_j(\tilde{\Ubf}^{(m),\star}_i+\widehat{\ep}_i^{(m)})| > Ct'\right)&=\Pbb(A_i\cap B_i)+\Pbb(A_i\cap B_i^c)
\end{align*}
where, using sub-Gaussianity,
\[\Pbb(A_i\cap B_i) \leq \Pbb(B_i)\leq \tilde{C} e^{-\frac{m^\alpha}{2\sigma^2}}\]
and using the locally Lipschitz condition on $f_j$, together with the independent of $\tilde{\ep}^{(m)}_i$ and $\widehat{\ep}^{(m)}_i$, 
\begin{align*}
\Pbb(A_i\cap B_i) &\leq \Pbb\left(\{C_\CalF(|\tilde{\ep}^{(m)}_i - \widehat{\ep}_i^{(m)}|+|\tilde{\ep}^{(m)}_i - \widehat{\ep}_i^{(m)}|^{p_{\max}})>Ct'\}\cap\{|\tilde{\ep}^{(m)}_i - \widehat{\ep}_i^{(m)}|<1\}\right)\\
&\leq \Pbb\left(|\tilde{\ep}^{(m)}_i - \widehat{\ep}_i^{(m)}|>\frac{C}{2C_\CalF}t'\right) \\ 
&\leq \tilde{C} e^{-\frac{m^\alpha}{2\sigma^2}\left(\frac{Ct'}{2C_\CalF}\right)^2}.
\end{align*}
Hence, summing over all $i$, we get 
\[\Pbb\left(|\tilde{\Gbf}^{(m)}_{ksj} - \widehat{\Gbf}_{ksj}^{(m)}|> t'/2\right) \leq \tilde{C}me^{-c_1 m^\alpha (t')^2}\]
for some fixed $c_1>0$, which shows that the term $E_2$ asymptotically dominates $E_1$.

Overall this shows that $\tilde{\Gbf}^{(m)}$ converges to $\Gbf^\star$ elementwise, and since the number of elements is fixed, we can take a union bound to conclude that $\|\tilde{\Gbf}^{(m)}-\Gbf^\star\|_\infty\to 0$ in probability at the same rate. We can characterize this convergence as follows. (1) For each $m$ we acrue a deterministic error $\mathfrak{t}(m)$ due to the bias introduced when averaging the clean data, as well as the integration errors. The worst-case asymptotic order of this error is $\CalO(m^{-\alpha \min\{\frac{1}{d+1},\frac{k}{d+1}-\frac{1}{2}\}})$ and arises from $I_4$, the difference between $\overline{\Gbf}^{(m),\star}$ which uses the filtered clean data $\tilde{\Ubf}^{(m),\star}$ and $\Gbf^{(m),\star}$ which uses the clean data $\Ubf^{(m),\star}$ (in effect, $I_4$ is the combined bias and integration error). (2) For $t>\mathfrak{t}(m)$ and large enough $m$ we have 
\[\Pbb \left(\|\tilde{\Gbf}^{(m)} - \Gbf^\star\|_\infty>t\right) \leq 3K\mathfrak{J}\exp\left(-\frac{c}{2^{1+2/{p_{\max}}}}[m(t-\mathfrak{t}(m))]^{2/{p_{\max}}}\right),\]
using the worst-case concentration rate arising from the error $I_2$.
\end{proof}

\section{Full support recovery}\label{app:supprec}

As indicated in Remark \ref{promise_of_fullsupprec}, in this section we extend the results \eqref{subsetsupprec_lemm} of Lemma \ref{lemm:subset_support_rec} and \eqref{subsetsupprec_thm} of Theorem \ref{conv_theorem_nosmooth} to subset equality
\[\supp{\widehat{\wbf}} = \supp{\wstar} \qquad \& \qquad\supp{\widehat{\wbf}^{(m)}}=\supp{\wstar},\]
respectively, with the addition of a constraint on the linear system $(\Gbf^\star,\bbf^\star)$. First we need another stability lemma similar to Lemma \ref{least_squares_cont}. 

\begin{lemm}\label{proj_cont}
Let $\Abf \in \Rbb^{m\times n}$ have rank $n$ and let $\tilde{\Abf} = \Abf+\Ebf$ be a perturbed system satisfying $\nrm{\Ebf}_{2}<\varepsilon$. For sufficiently small $\varepsilon$, there exist constants $C,C'>0$ depending only on $\Abf$ such that the following stability holds for the pseudoinverse $\Abf^\dagger := (\Abf^T\Abf)^{-1}\Abf^T$ and projection operator $\Pbf_\Abf := \Abf\Abf^\dagger$:
\begin{equation}
\nrm{\Abf^\dagger - \tilde{\Abf}^\dagger}_2<C\varepsilon \qquad \& \qquad \nrm{\Pbf_\Abf-\Pbf_{\tilde{\Abf}}}_2<C'\varepsilon.
\end{equation}
\end{lemm}
\begin{proof}
The pseudoinverse of $\tilde{\Abf}$ is given by 
\[\tilde{\Abf}^\dagger = (\tilde{\Abf}^T\tilde{\Abf})^{-1}\tilde{\Abf}^T = (\Ibf+\Bbf)^{-1}(\Abf^\dagger+(\Abf^T\Abf)^{-1}\Ebf^T),\]
provided $(\Ibf+\Bbf)^{-1}$ exists, where 
\[\Bbf = \Abf^\dagger\Ebf+(\Abf^T\Abf)^{-1}\Ebf^T(\Abf + \Ebf).\]
A sufficient condition for invertibility of $\Ibf+\Bbf$ is 
\[\nrm{\Bbf}_2<1,\]
which is guaranteed for sufficiently small $\varepsilon$ since $\nrm{\Bbf}_2 = \CalO(\varepsilon)$. In this case using the Von Neumann series for $(\Ibf+\Bbf)^{-1}$, we get that 
\[\tilde{\Abf}^\dagger = \Abf^\dagger + \tilde{\Ebf},\]
where 
\[\tilde{\Ebf} = (\Abf^T\Abf)^{-1}\Ebf^T+\sum_{k=1}^\infty(-\Bbf)^k\left(\Abf^\dagger+(\Abf^T\Abf)^{-1}\Ebf^T\right).\]
Since $\Ebf=\CalO(\varepsilon)$ and $\Bbf = \CalO(\varepsilon)$, we have $\tilde{\Ebf} = \CalO(\varepsilon)$, hence
\[\nrm{\Abf^\dagger - \tilde{\Abf}^\dagger}_2 = \nrm{\tilde{\Ebf}}_2 = \CalO(\varepsilon).\]
The same stability readily applies to $\Pbf_{\tilde{\Abf}} = \Pbf_{\Abf} + \Abf\tilde{\Ebf}+\Ebf\tilde{\Abf}^\dagger$.
\end{proof}

We can now use this stability to extend Lemma \ref{lemm:subset_support_rec} to yield full support recovery on the continuum problem for small enough $\sigma^2$, provided condition \eqref{app:condforsupprec} is satisfied.

\begin{lemm}\label{lemm:fullsupprec}
Let $\rho$ be a Gaussian mean-zero noise distribution and let $p$ be the maximum polynomial degree appearing in the true model. Assume that $(\Gbf^\star,\bbf^\star)$ satisfy
\begin{equation}\label{app:condforsupprec}
\mu^\star:= \min_{S\subsetneq S^\star}\frac{\nrm{\Pbf^\perp_{\Gbf^\star_{S^\star\setminus S}} \bbf^\star}}{\nrm{\bbf^\star}}-\frac{|S|+1}{\mathfrak{J}}>0.
\end{equation}
Then there exists a critical noise level $\sigma_c$ such that for any $\sigma\leq \sigma_c$ the estimator $\what = \text{MSTLS}^{(1)}(\overline{\Gbf},\overline{\bbf})$ satisfies 
\begin{equation}\label{subsetsupprec_lemm_app}
\supp{\what}=\supp{\wstar}.
\end{equation}
\end{lemm}

\begin{proof}
In Lemma \ref{lemm:subset_support_rec} we already showed that 
\begin{enumerate}[label=(\alph*)]
\item For sufficiently small $\sigma$ there exists $\widehat{\lambda}\in \pmb{\lambda}$ such that $\supp{\overline{\wbf}^{\widehat{\lambda}}} = \supp{\wstar}$, where $\overline{\wbf}^{\widehat{\lambda}} = \text{STLS}^{(1)}(\overline{\Gbf},\overline{\bbf},\widehat{\lambda})$.
\item For sufficiently small $\sigma$ we have 
\[\frac{\nrm{\overline{\Gbf}\left(\overline{\wbf}^0 -\overline{\wbf}^{\widehat{\lambda}}\right)}_2}{\nrm{\overline{\bbf}}_2} <\frac{1}{\mathfrak{J}}\]
which implies that for $\lambda\in \pmb{\lambda}$ with $\lambda<\widehat{\lambda}$, we have 
\[\CalL(\widehat{\lambda})<\CalL(\lambda).\]
\end{enumerate}
Assume that $\sigma$ is small enough so that (a) and (b) hold, and consider $\lambda\in \pmb{\lambda}$ with $\lambda>\widehat{\lambda}$. All we need to show is that for sufficiently small $\sigma$, we have $\CalL(\widehat{\lambda})<\CalL(\lambda)$ in this case as well. Indeed, by construction of $\pmb{\lambda}$, for $\lambda>\widehat{\lambda}$ we have
\[\supp{\overline{\wbf}^{\lambda}} \subsetneq S^\star.\]
Setting $S = S^\star \setminus \supp{\overline{\wbf}^{\lambda}}$, we see that
\[\CalL(\widehat{\lambda}) - \CalL(\lambda) = \frac{\nrm{\overline{\Gbf}\left(\overline{\wbf}^0 -\overline{\wbf}^{\widehat{\lambda}}\right)}_2}{\nrm{\overline{\bbf}}_2} - \frac{\nrm{\overline{\Gbf}\left(\overline{\wbf}^0 -\overline{\wbf}^{\lambda}\right)}_2}{\nrm{\overline{\bbf}}_2}+\frac{|S|}{\mathfrak{J}}.\]
Using (b), and rewriting using the projection operator, we have
\[\CalL(\widehat{\lambda}) - \CalL(\lambda) < - \frac{\nrm{\overline{\Gbf}\left(\overline{\wbf}^0 -\overline{\wbf}^{\lambda}\right)}_2}{\nrm{\overline{\bbf}}_2}+\frac{|S|+1}{\mathfrak{J}} = -\frac{\nrm{\Pbf^\perp_{\overline{\Gbf}_{S^\star \setminus S}}\overline{\bbf}}_2}{\nrm{\overline{\bbf}}_2}+\frac{|S|+1}{\mathfrak{J}}.\]
Now, for small enough $\sigma$, using stability of the projection operator, we have 
\[\max_{S\subsetneq S^\star}\nrm{\Pbf^\perp_{\overline{\Gbf}_{S^\star \setminus S}} - \Pbf^\perp_{\Gbf^\star_{S^\star \setminus S}} }_2 < C\sigma^2,\]
since properties of $\Dbf$ and $\Lbf$ imply
\[\overline{\Gbf}_{S^\star \setminus S} -\Gbf^\star_{S^\star \setminus S} = \Gbf^\star (\Abf_{S^\star \setminus S} -\Ibf_{S^\star \setminus S}) = \Gbf^\star(\Dbf-\Ibf)_{S^\star \setminus S}+ \Gbf^\star \Lbf_{S^\star \setminus S} =  \CalO(\sigma^2).\]
Hence, using the reverse triangle inequality we get 
\begin{align*}
\CalL(\widehat{\lambda}) - \CalL(\lambda) &< -\frac{\nrm{\Pbf^\perp_{\Gbf^\star_{S^\star \setminus S}}\overline{\bbf}}_2}{\nrm{\overline{\bbf}}_2} + \frac{\nrm{\left(\Pbf^\perp_{\Gbf^\star_{S^\star \setminus S}}- \Pbf^\perp_{\Gbf^\star_{S^\star \setminus S}}\right)\overline{\bbf}}_2}{\nrm{\overline{\bbf}}_2} +\frac{|S|+1}{\mathfrak{J}} \\
&< -\mu^\star + C\sigma^2
\end{align*}
which is negative for all $\sigma^2 <\mu^\star/C$. Thus, to recover the full support on systems satisfying \eqref{app:condforsupprec}, in addition to other constraints derived on $\sigma_c$, it holds that $\sigma_c^2\geq \mu^\star/C$. 
\end{proof}

The previous Lemma, together with results shown in Section \ref{nosmooth}, directly imply the following, using similar techniques as Theorem \ref{conv_theorem_nosmooth}.

\begin{thm}\label{full_theorem_nosmooth}
Let Assumptions \ref{sec:regularity}-\ref{sec:conditioning} hold with $\rho$ a mean-zero Gaussian distribution. In addition assume that the condition \eqref{app:condforsupprec} holds. Then exists a critical noise $\sigma_c>0$ and a stability tolerance $\tau$, both independent of $m$, such that for all $\sigma<\sigma_c$ and $t< \tau$, and for sufficiently large $m$, it holds that 
\begin{equation}\label{supprec_thm}
\supp{\widehat{\wbf}^{(m)}}=\supp{\wstar} \qquad \text{and} \qquad \nrm{\widehat{\wbf}^{(m)}-\wstar}_\infty < C'(t+\sigma^2)
\end{equation}
with probability exceeding $1-2K\mathfrak{J}\exp\left(-\frac{c}{2}\left(mt\right)^{2/p_{\max}}\right)$, where $\widehat{\wbf}^{(m)} = \text{MSTLS}^{(1)}(\Gbf^{(m)},\bbf^{(m)})$, $c$ is from Theorem \ref{matrix_concentration}, and $C'$ depends on only $(\Gbf^\star_{S^\star},\bbf^\star)$. 
\end{thm}

\section{Hyper-KS dynamics}\label{app:hyperKS}

\begin{equation}\label{hyperKSapp}
\partial_t u = (1)\partial_{xxxx} u + (0.75)\partial_{xxxxxx} u + (-0.5)\partial_x(u^2) + (0.1)\partial_{xxx} u^2.
\end{equation}

Motivation for equation \ref{hyperKSapp} comes from the fact that variants of this equation arise when performing equation learning on corrupted Kuramoto-Sivashinsky data using a disadvantageous test function, if derivatives up to order 6 are included in the library. (See for example Figure 16 of \cite{tang2022weakident}.) Two effects dominate the dynamics. At low wavelengths $|k|\leq k^*$ below some threshold $k^*$ the system is unstable. Then at all wavelengths sufficiently separated from a critical wavenumber $\overline{k}$, we have dispersive mixing of wavemodes. At the critical wavenumber $\overline{k}$, dispersive and transport effects cancel out, such that $\overline{k}$ is entirely subject to the growth and decay dynamics of the diffusive terms. This implies that if $|\overline{k}|< k^*$, the system can become unstable.

\section{Variance estimation}\label{app:sigma_est}

We estimate the noise variance using the following procedure. Let $\fbf\in \Rbb^{2L+1}$ satisfy the discrete moment conditions 
\[M_k(\fbf) = \sum_{j=-L}^L  j^k\fbf_f = 0, \quad k=0,\dots,M_{\max}.\]
and 
\[\nrm{\fbf}_2=1.\]
Now consider data $\Ubf = u(\tbf)+\ep$ such that $u$ is locally a polynomial of degree $M_{\max}$ over intervals $(t-L\Delta t, t+L\Delta t)$ and $\ep$ is mean zero i.i.d.\ noise with variance $\sigma^2$. Then the data $\Ubf$ will satisfy
\[\Ebb\left[|(\fbf*\Ubf)^2_i\right] = \sigma^2, \quad i\in \Zbb,\]
hence 
\[\nrm{\fbf*\Ubf}_{rms}\approx \sigma\]
to very high accuracy. We use this approach and take $\fbf$ to be 6th-order finite difference weights normalized so that $\nrm{\fbf}_2=1$.

\bibliographystyle{plain}
\bibliography{researchCU}

\end{document}

%% file: ex_shared.tex
\usepackage{lipsum}
\usepackage{amsfonts}
\usepackage{graphicx}
\usepackage{epstopdf}
\usepackage{subcaption}
\usepackage{algorithmic}
\ifpdf
  \DeclareGraphicsExtensions{.eps,.pdf,.png,.jpg}
\else
  \DeclareGraphicsExtensions{.eps}
\fi

\newsiamremark{remark}{Remark}
\newsiamremark{hypothesis}{Hypothesis}
\crefname{hypothesis}{Hypothesis}{Hypotheses}
\newsiamthm{claim}{Claim}

\headers{Asymptotic consistency of the WSINDy algorithm}{D. A. Messenger and D. M. Bortz}

\title{Asymptotic consistency of the WSINDy algorithm in the limit of continuum data}

\author{Daniel A. Messenger\thanks{Department of Applied Mathematics, University of Colorado, Boulder, CO 80309-0526, USA. 
	(\email{daniel.messenger@colorado.edu}, \email{dmbortz@colorado.edu}).}
\and David M. Bortz\footnotemark[1]}

\usepackage[latin9]{inputenc}
\usepackage{geometry}
\PassOptionsToPackage{normalem}{ulem}

\usepackage{
amsopn,
amsmath,
amssymb,
appendix,
bm,
bbm,
caption,
color,
epsf,
enumitem,
float,
graphicx,
hyperref,
listings,
mathtools,
mathrsfs,
pdfcolmk,
subcaption,
tikz,
titletoc,
url,
ulem,
xcolor
}

\hypersetup{
    colorlinks=true, 
    linktoc=all,     
    linkcolor=blue,  
}

\numberwithin{equation}{section}

\makeatletter

\providecolor{lyxadded}{rgb}{0,0,1}
\providecolor{lyxdeleted}{rgb}{1,0,0}

\DeclareRobustCommand{\lyxsout}[1]{\ifx\\#1\else\sout{#1}\fi}

\input{notation}

\usepackage{stmaryrd}

\makeatother


%% file: notation.tex
\newcommand{\Ebb}{\mathbb{E}}
\newcommand{\Nbb}{\mathbb{N}}
\newcommand{\Pbb}{\mathbb{P}}
\newcommand{\Rbb}{\mathbb{R}}

\newcommand{\Vbb}{\mathbb{V}}
\newcommand{\Zbb}{\mathbb{Z}}

\newcommand{\Abf}{\mathbf{A}}
\newcommand{\bbf}{\mathbf{b}}
\newcommand{\Bbf}{\mathbf{B}}

\newcommand{\Cbf}{\mathbf{C}}

\newcommand{\Dbf}{\mathbf{D}}

\newcommand{\Ebf}{\mathbf{E}}
\newcommand{\fbf}{\mathbf{f}}

\newcommand{\Gbf}{\mathbf{G}}

\newcommand{\Ibf}{\mathbf{I}}

\newcommand{\Lbf}{\mathbf{L}}

\newcommand{\Pbf}{\mathbf{P}}

\newcommand{\tbf}{\mathbf{t}}
\newcommand{\Tbf}{\mathbf{T}}

\newcommand{\Ubf}{\mathbf{U}}

\newcommand{\wbf}{\mathbf{w}}

\newcommand{\xbf}{\mathbf{x}}
\newcommand{\Xbf}{\mathbf{X}}
\newcommand{\ybf}{\mathbf{y}}
\newcommand{\Ybf}{\mathbf{Y}}

\newcommand{\ep}{\epsilon}

\newcommand{\CalF}{{\mathcal{F}}}

\newcommand{\CalH}{{\mathcal{H}}}

\newcommand{\CalL}{{\mathcal{L}}}

\newcommand{\CalO}{{\mathcal{O}}}
\newcommand{\CalP}{{\mathcal{P}}}
\newcommand{\CalQ}{{\mathcal{Q}}}

\newcommand{\CalS}{{\mathcal{S}}}

\newcommand{\nrm}[1]{\left\Vert {#1} \right\Vert}

\newcommand{\supp}[1]{\text{supp}\left(#1\right)}

\newcommand{\ind}[1]{\mathbbm{1}_{#1}}

\renewcommand{\tilde}[1]{\widetilde{#1}}

\DeclareMathOperator*{\argmin}{argmin}

\newcommand{\lan}{\left\langle}
\newcommand{\ran}{\right\rangle}

\newcommand{\wstar}{\wbf^\star}
\newcommand{\what}{{\widehat{\wbf}}}

\newtheorem{thm}{Theorem}[section]
\newtheorem{lemm}{Lemma}[section]

\newtheorem{corr}{Corollary}[section]
\newtheorem{rmrk}{Remark}[section]